\documentclass[1p,fleqn,10pt]{elsarticle}

\usepackage{graphicx} 
\usepackage{amsfonts,amsmath,amssymb,amsthm}
\usepackage{pgfplots}
\usepackage{xcolor}
\usepackage{etoolbox}
\usepackage{subcaption}
\usepackage{caption}
\usepackage{cancel}
\usepackage{nicefrac}
\usepackage[framemethod=TikZ]{mdframed}
\usetikzlibrary{spy}
\usetikzlibrary{calc}
\usepackage{booktabs}
\usepackage{mathtools}

\definecolor{GREEN}{HTML}{1b9e77}
\definecolor{ORANGE}{HTML}{d95f02}
\definecolor{MAGENTA}{HTML}{7570b3}

\newtheorem{theorem}{Theorem}[section]
\newtheorem{lemma}{Lemma}[section]
\newtheorem{proposition}{Proposition}[section]
\newtheorem{corollary}{Corollary}[section]

\newtheorem{assumption}{Assumption}[section]

\newtheorem{formulation}{Formulation}[section]

\newtheorem{remark}{Remark}[section]
\begin{document}

\newcommand{\curl}[1]{\mathrm{curl}\ifstrempty{#1}{}{\left(#1\right)}}
\newcommand{\rot}[1]{\mathrm{rot}\ifstrempty{#1}{}{\left(#1\right)}}
\renewcommand{\div}[1]{\mathrm{div}\ifstrempty{#1}{}{\left(#1\right)}}
\newcommand{\grad}[1]{\mathrm{grad}\ifstrempty{#1}{}{\left(#1\right)}}

\renewcommand{\deg}[0]{k}

\renewcommand{\vec}[1]{\boldsymbol{#1}}

\newcommand{\Ci}[0]{\mathrm{C_{inv}}}
\newcommand{\Ce}[0]{\mathrm{C_{est}}}
\newcommand{\tm}[0]{\mathrm{\tau_{M}^{-1}}}
\newcommand{\tc}[0]{\mathrm{\tau_{c}^{-1}}}

\newcommand{\PicIt}[0]{m}

\newcommand{\closure}[1]{\mathrm{clos}\left(#1\right)}


\newenvironment{assumpBox}
{
	\begin{mdframed}[
		outerlinewidth=0pt,
		roundcorner=5pt,
		innertopmargin=\baselineskip,
		innerbottommargin=\baselineskip,
		innerrightmargin=14pt,
		innerleftmargin=14pt,
		backgroundcolor=gray!14!white,
		fontcolor=black,
		align=center,
		skipabove=6pt,
		skipbelow=6pt]
	\begin{assumption}
}%
{
	\end{assumption}
	\end{mdframed}
}

\newenvironment{methodBox}
{
	\begin{mdframed}[
		outerlinewidth=0pt,
		roundcorner=5pt,
		innertopmargin=5pt,
		innerbottommargin=5pt,
		innerrightmargin=14pt,
		innerleftmargin=14pt,
		backgroundcolor=gray!14!white,
		fontcolor=black,
		align=center,
		skipabove=6pt,
		skipbelow=6pt]
	\begin{formulation}
}%
{
	\end{formulation}
	\end{mdframed}
}

\author[label1]{Kevin Dijkstra}
\ead{kwdijkstra@tudelft.nl}
\author[label1]{Deepesh Toshniwal}
\ead{d.toshniwal@tudelft.nl}

\address[label1]{Delft Institute of Applied Mathematics, Delft University of Technology - The Netherlands}

\title{Structure-preserving Variational Multiscale Stabilization of the Incompressible Navier--Stokes Equations}
\begin{abstract}
    This paper introduces a Variational Multiscale Stabilization (VMS) formulation of the incompressible Navier--Stokes equations that utilizes the Finite Element Exterior Calculus (FEEC) framework.
    The FEEC framework preserves the geometric and topological structure of continuous spaces and PDEs in the discrete spaces and model, and helps build stable and convergent discretizations.
    For the Navier--Stokes equations, this structure is encoded in the de Rham complex.
    In this work, we consider the vorticity-velocity-pressure formulation discretized within the FEEC framework.
    We model the effect of the unresolved scales on the finite-dimensional solution by introducing appropriate fine-scale governing equations, which we also discretize using the FEEC approach.
    This preserves the structure of the continuous problem in both the coarse- and fine-scale solutions; for instance, both the coarse- and fine-scale velocities are pointwise incompressible.
    We demonstrate that the resulting formulation is residual-based, energetically stable, and optimally convergent.
    Moreover, our fine-scale model allows for an efficient computational approach -- by decoupling fine-scale problems on one element from another, it is possible to solve them in parallel.
    In fact, the fine-scale equations can be entirely eliminated during matrix assembly, thereby leading to a VMS formulation where the problem size is governed by the coarse-scale discretization only.
    Finally, the proposed formulation applies to both the lowest regularity discretizations of the de Rham complex and high-regularity isogeometric discretizations.
    We validate our theoretical results through numerical experiments, simulating both steady-, unsteady-, viscous-, and inviscid-flow problems.
    These tests show that the stabilized solutions are qualitatively better than the unstabilized ones, converge at optimal rates, and, as the mesh is refined, the stabilization is asymptotically turned off.
\end{abstract}

\begin{keyword}
    Incompressible Navier--Stokes equations \sep Vorticity-velocity-pressure formulation \sep Variational multiscale analysis \sep The de Rham complex \sep Finite element exterior calculus \sep Divergence-conforming discretizations
\end{keyword}

\maketitle

\section{Introduction}
Structure-preserving methods offer a robust and accurate framework for discretizing partial differential equations (PDEs), including the incompressible Navier--Stokes equations. 
These methods aim to conserve essential physical invariants, such as mass, momentum, or kinetic energy, in the discrete setting. This focus enhances the numerical properties compared to standard discretizations.
For instance, \cite{arakawa_computational_1966} demonstrated that conservation of energy and enstrophy was crucial for improving long-time simulations of the incompressible Navier--Stokes equations.

As a result, constructing structure-preserving formulations of the incompressible Navier--Stokes Equations that preserve invariants has been an active area of research.
For instance, in \cite{perot_conservation_2000}, it was noted that discretizations of the vorticity formulation can conserve kinetic energy up to machine precision, and the divergence formulation can conserve kinetic energy and momentum.
Several such structure-preserving methods for fluid flows have since been proposed within the frameworks of Discrete Exterior Calculus (DEC) \cite{desbrun_discrete_2005,elcott_stable_2007,mohamed_discrete_2016,bothe_upwind_2017}, Finite Element Exterior Calculus (FEEC) \cite{evans_isogeometric_2013-1,evans_isogeometric_2013, evans_variational_2020, hanot_arbitrary_2023,tonnon_semi-lagrangian_2024,carlier_mass_2025}, mimetic methods (mass, energy, enstrophy and vorticity conserving)  \cite{palha_mass_2017,zhang_mass-_2022} and others, such as \cite{rebholz_energy-_2007, olshanskii_velocityvorticityhelicity_2010}.

In this work, we will follow the FEEC framework \cite{arnold_finite_2018}, where the essential geometric and topological structure embedded in the continuous spaces and PDE is preserved, to obtain stable and convergent finite element discretizations.
In the case of incompressible Navier--Stokes, this structure is encoded in the so-called two- or three-dimensional de Rham Hilbert complex. For example, in three dimensions, this complex is given by
\begin{equation}\label{eq:deRhamComplex3d}
    H^1(\Omega) \xrightarrow{\mathrm{grad}} H(\curl{};\Omega)\xrightarrow{\mathrm{curl{}}} H(\div{};\Omega) \xrightarrow{\mathrm{div}} L^2(\Omega)\;.
\end{equation}
Finite Element Exterior Calculus has been utilized in various ways for the incompressible Navier--Stokes equations.
The velocity-pressure formulation of the incompressible Navier--Stokes equations was explored in \cite{evans_isogeometric_2013-1,evans_isogeometric_2013}, where only the second half of \eqref{eq:deRhamComplex3d}, specifically $H(\div{};\Omega) \xrightarrow{\mathrm{div}} L^2(\Omega)$, was discretized with divergence-conforming B-spline spaces.
This approach ensures that the finite element discretization satisfies pointwise incompressibility and admits discrete balance laws for momentum, angular momentum, energy, vorticity, enstrophy, and helicity.
By only considering the second half of the complex, the tangential no-slip boundary condition had to be enforced weakly via Nitsche's method.
More importantly, these works use the symmetric gradient operator of the incompressible Navier--Stokes equations.
Instead, in Finite Element Exterior Calculus \cite{arnold_finite_2018}, the Hodge Laplacian plays a crucial role in the structure of the de Rham complex; hence, treating the diffusion operator as a Hodge Laplacian is the more natural approach.
However, the Hodge Laplacian relies on the adjoint of the discrete differential operators, which are global and deteriorate the sparsity of the system \cite{hanot_arbitrary_2023}.
While the issue of the adjoint operators can be resolved, for example through broken FEEC \cite{carlier_mass_2025}), the adjoint operators can also be eliminated by the introduction of an auxiliary variable $\omega = \curl{\vec{u}}$ (the vorticity), leading to the vorticity-velocity-pressure formulation; see \cite{hiemstra_high_2014} for an example discretization based on FEEC and isogeometric analysis.
Examples of methods for the vorticity-velocity-pressure formulation are \cite{anaya_analysis_2019, hanot_arbitrary_2023}, where the linearizations of the incompressible Navier--Stokes equations are studied, resulting in pointwise incompressibility of the velocity field.
Compared to \cite{evans_isogeometric_2013}, the resulting systems are larger to solve, due to the need to solve for the auxiliary variable.
Alternatively, in \cite{palha_mass_2017}, the velocity-pressure formulation is extended by introducing the vorticity and supplementing it with a vorticity transport equation, leading to a coupled system of nonlinear equations. 
For computational efficiency, a leapfrog temporal discretization scheme was used, resulting in a coupled system of linear equations.
Additionally, the vorticity-velocity-pressure formulation allows the boundary conditions to be incorporated without the need for additional tools (e.g., Nitsche's method). 

Structure-preserving methods create discretizations that maintain physical invariants along with geometric and topological structures. However, when transitioning to a discrete space, we inherently overlook the effects of the scales that these discrete spaces cannot capture. This is not just an issue for structure-preserving methods but affects all discrete numerical methods.
The variational multiscale stabilization (VMS) framework \cite{hughes_multiscale_1995, hughes_variational_1998} was developed to incorporate the missing, uncaptured effects (the so-called fine scales) into the equations governing the captured scales (the so-called coarse scales).
Here, the governing equations are separated into coarse- and fine-scale equations with respect to a projector.
If the Green's function is known, the solution of the fine-scale governing equations can be found explicitly  \cite{hughes_variational_2007}. 
Alternatively, the solution needs to be approximated either via a subgrid model or a discretization of the fine-scale equations.
This fine-scale solution is then incorporated into the coarse-scale governing equations, which were found to have a stabilizing effect. 
Indeed, under the right conditions, this technique reproduces the well-known Streamline Upwind/Petrov-Galerkin (SUPG) stabilization \cite{brooks_streamline_1982,hughes_variational_2007}. 

The VMS framework has been used extensively to develop stabilized formulations of the incompressible Navier--Stokes equations \cite{hughes_large_2000,bazilevs_variational_2007,bazilevs_isogeometric_2010,ten_eikelder_correct_2018,ahmed_review_2015}.
It is important to note that, in general, VMS formulations do not adhere to the de Rham complex without careful consideration.
For example, consider the residual-based VMS formulation introduced in \cite{bazilevs_variational_2007}.
In this formulation, it is assumed that the strong residual, when multiplied by a stabilization parameter, approximates the fine-scale solutions effectively.
As a result, only the coarse-scale fields need to be computed explicitly, after substitution. 
However, since the strong residual is in general not divergence-free, nor is the resulting fine-scale solution, conservation of mass is not guaranteed by the method.
Instead, in \cite{van_opstal_isogeometric_2017,evans_variational_2020}, where the velocity-pressure formulation is considered, both the coarse- and fine-scale solutions are constructed to be pointwise divergence free (conservation of mass).
In \cite{van_opstal_isogeometric_2017}, the authors introduce an elementwise fine-scale model, based on the Darcy equations, that conforms to an elementwise version of the de Rham complex \eqref{eq:deRhamComplex3d}.
Multiple boundary conditions are discussed, such as homogeneous Dirichlet, homogeneous Neumann, and periodic boundary conditions.
This approach can be considered a divergence-conforming extension of the residual-free bubble technique \cite{brezzi_choosing_1994,brezzi_bint_1997}.
Although this results in a significantly larger number of unknowns, the elementwise nature of the model makes it efficiently solvable.
Alternatively, in \cite{evans_variational_2020}, a global fine-scale model is constructed where the fine-scale velocity is weakly divergence-free. However, both the coarse- and fine-scale pressure must be $H^1$-conforming and explicitly computed to solve the system.
Furthermore, the higher regularity of the coarse-scale pressure makes this method unsuitable for low-order discretizations.
More recently, VMS was used within FEEC \cite{shrestha_optimal_2025} with the vorticity-velocity-pressure formulation.
Here, the fine-scale computation requires the solution of the symmetric part of the PDE with discrete spaces defined on a finer mesh and with a higher polynomial degree.
As a consequence, large coupled nonlinear problems must be solved simultaneously.
The authors numerically demonstrate the convergence of the method and apply it to the 2D MEEVC formulation.

In this work, we address the velocity-vorticity-pressure formulation of the incompressible Navier--Stokes Equations with VMS and FEEC.
Here, we follow the work of \cite{van_opstal_isogeometric_2017}, where the fine-scale model is defined elementwise and discretized with bubble functions that conform to the complete de Rham complex \eqref{eq:deRhamComplex3d}, and where appropriate Dirichlet boundary conditions are imposed on the skeleton of the mesh.
As in \cite{van_opstal_isogeometric_2017,evans_variational_2020}, we will consider fine-scale governing equations based on the Darcy equations with a modelling parameter; however, for our analysis, we augment them with an additional diffusion term so that the fine-scale governing equations might be better characterized as the Darcy--Stokes equations.
In addition, based on the work of \cite{ten_eikelder_correct_2018}, we introduce additional terms that imply an explicit equation for the kinetic energy evolution.
In order not to impose any additional regularity on the coarse spaces, we use the weak residual rather than the strong residual to drive the fine-scale governing equations, which allows our method to be applicable to both low regularity discretizations of the de Rham complex, such as the mimetic finite elements \cite{palha_physics-compatible_2014} or high regularity isogeometric discretizations \cite{Buffa2011, buffa_isogeometric_2011, evans_isogeometric_2013,hiemstra_high_2014}.
We show convergence and stability of our method under suitable assumptions, which we numerically validate.
Here, we qualitatively observe a reduction in oscillations in under-resolved solutions and additional dissipation of kinetic energy due to the introduction of the fine-scale solution.
Finally, the behavior of both coarse and fine-scale solutions can be adjusted using the stabilization parameter to fit the model to validated data.

In Section \ref{sec:background}, the basic notions of structure-preserving techniques for the incompressible Navier--Stokes equations are introduced. Section \ref{sec:VMS-method} introduces our VMS formulation, which is shown to be energetically stable. In Section \ref{sec:convergence}, our formulation is applied to the Oseen equation, and uniqueness, stability, and optimal convergence rates are shown. Lastly, the formulation is numerically validated in Section \ref{sec:numerics}.

\section{Structure-preserving discretizations for incompressible Navier--Stokes}
\label{sec:background}

\subsection{Function spaces}
Let $L^2(\Omega)$ denote the space of square integrable functions
\begin{equation}
    L^2(\Omega) := \left\{~f~:~\left(f,f\right):=\int_{\Omega}f^2dV <\infty~ \right\}\;.
\end{equation}
In addition, we introduce the following Hilbert spaces where the solutions of the continuous problem live:
\begin{align}
    H^1(\Omega) &:= \left\{~f\in L^2(\Omega)~:~\grad{f}\in\left(L^2(\Omega)\right)^d~ \right\}\;,\\
    H(\curl{};\Omega) &:= \left\{~\vec{f}\in \left(L^2(\Omega)\right)^d~:~\curl{\vec{f}}\in\left(L^2(\Omega)\right)^d~ \right\}\;,\\
    H(\div{};\Omega) &:= \left\{~\vec{f}\in \left(L^2(\Omega)\right)^d~:~\div{\vec{f}}\in L^2(\Omega)~ \right\}\;.
\end{align}
The trace operators, which restrict fields to $\Gamma \subseteq \partial\Omega$ are defined as:
\begin{equation}
    \mathcal{T}f = \left.f\right|_{\Gamma},~\mathcal{T}_{\perp}\vec{f} = \left.\vec{f}\cdot\vec{n}\right|_{\Gamma},~\mathcal{T}_{\parallel}\vec{f} = \left.\vec{f}\times\vec{n}\right|_{\Gamma},~\mathcal{T}_{\vec{n}}\vec{f} = \left.\vec{n}\times\left(\vec{f}\cdot\vec{n}\right)\right|_{\Gamma}.
\end{equation}
The corresponding spaces with boundary conditions are defined as:
\begin{align}
    H^1_{\phi}(\Omega,\Gamma) &:= \left\{~f\in H^1(\Omega)~:~\mathcal{T}f = \phi\text{ on }\Gamma~ \right\}\;,\\
    H_{\vec{\phi}}(\curl{};\Omega, \Gamma) &:= \left\{~\vec{f}\in H(\curl{};\Omega)~:~\mathcal{T}_{\vec{n}}\vec{f} = \vec{\phi}\text{ on }\Gamma~ \right\}\;,\\
    H_{{\phi}}(\div{};\Omega, \Gamma) &:= \left\{~\vec{f}\in H(\div{};\Omega)~:~\mathcal{T}_{\perp}\vec{f} = {\phi}\text{ on }\Gamma~ \right\}\;,\\
    L^2_{\phi}(\Omega) &:= \left\{~f\in L^2(\Omega)~:~(f,1) = \phi~ \right\}\;.
\end{align}
We will use the natural norms for all of these spaces and, to simplify notation, we will denote the $L^2$-norm as $\|\cdot\|$.
Finally, we will denote with $H^s(\Omega)$, $s \in \mathbb{N}$, the space of $L^2$ functions such that their derivatives up to order $s$ are also in $L^2(\Omega)$.
Similarly, we will denote with $H^s(\curl{};\Omega)$ (respectively, $H^s(\div{};\Omega)$) the space of $L^2$ vector fields such that their curl (respectively, divergence) is in $(H^{s}(\Omega))^d$ (respectively, $H^{s}(\Omega)$).

\subsection{Continuous weak formulation}
The strong form of the incompressible Navier--Stokes equations in the velocity-vorticity-pressure form is:
Find $\vec{u}:\Omega\times[0,T]\rightarrow\mathbb{R}^d,\vec{\omega}:\Omega\times(0,T]\rightarrow\mathbb{R}^d$ and $p:\Omega\times(0,T]\rightarrow \mathbb{R}$ such that,
\begin{subequations}\label{eq:ss-ns-rotational-from}
    \begin{align}
    \partial_t \vec{u} + \vec{\omega} \times \vec{u} - \mathrm{Re}^{-1} \curl{\vec{\omega}} +  \grad{p} &= \vec{f}\;,\\
    \div{\vec{u}} &= 0\;,\\
    \vec{\omega} - \curl{\vec{u}} &= \vec{0}\;,
\end{align}
\end{subequations}
Here $0<\mathrm{Re}<\infty$ is the Reynolds number and $\vec{f}$ is a source term. This equation is closed with the boundary conditions:
\begin{equation}\label{eq:intro-boundary-conditions}
    \begin{cases}
        p = \hat{p} &\text{on }\Gamma_{\hat{p}}\times[0,T]\\
        \vec{u}\cdot\vec{n} = \hat{u} &\text{on }\Gamma_{\hat{u}}\times[0,T]\\
    \end{cases}\qquad\text{and}\qquad 
    \begin{cases}
        \vec{u}\times\vec{n} = \hat{\vec{u}} &\text{on }\Gamma_{\hat{\vec{u}}}\times[0,T]\\
        \vec{n}\times\left(\vec{\omega}\times\vec{n}\right) = \hat{\vec{\omega}} &\text{on }\Gamma_{\hat{\vec{\omega}}}\times[0,T]\\
    \end{cases}
\end{equation}
where $\vec{n}$ is the outward unit normal vector, $\left\{\Gamma_{\hat{p}},\Gamma_{\hat{u}}\right\}$ and $\left\{\Gamma_{\hat{\vec{u}}},\Gamma_{\hat{\vec{\omega}}}\right\}$ partition $\partial \Omega$, such that $\Gamma_{\hat{p}}\cup \Gamma_{\hat{u}} = \partial \Omega$, $\Gamma_{\hat{p}}\cap \Gamma_{\hat{u}}=\emptyset$, and $\Gamma_{\hat{\vec{u}}}\cup \Gamma_{\hat{\vec{\omega}}} = \partial \Omega$, $\Gamma_{\hat{\vec{u}}}\cap \Gamma_{\hat{\vec{\omega}}}=\emptyset$.

The following weak form coresponding to \eqref{eq:ss-ns-rotational-from} can be derived: 
given the initial velocity field $\vec{u}_0 \in H_{\hat{u}}(\div{};\Omega,\Gamma_{\hat{u}})$, find $\vec{\omega}\in H_{\hat{\vec{\omega}}}(\curl{};\Omega,\Gamma_{\hat{\vec{\omega}}})$, $\vec{u}\in L^2((0,T);H_{\hat{u}}(\div{};\Omega,\Gamma_{\hat{u}}))\cap H^1((0,T);L^2(\Omega)^d)$ and $p\in L^2(\Omega)$, such that for all $\vec{\tau}\in H_{\vec{0}}(\curl{};\Omega,\Gamma_{\hat{\vec{\omega}}})$, $\vec{v}\in H_0(\div{};\Omega,\Gamma_{\hat{u}})$ and $q\in L^2(\Omega)$, it holds that
\begin{subequations}
\label{eq:intro-weak-NS}
    \begin{align}
        (\partial_t \vec{u},\vec{v})+(\vec{\omega}\times \vec{u},\vec{v}) + \mathrm{Re}^{-1}(\curl{\vec{\omega}},\vec{v})- (p,\div{\vec{v}}) &= (\vec{f},\vec{v}) - (\hat{p},\mathcal{T}_{\perp}\vec{v})_{\Gamma_{\hat{p}}}\;,\\
        (q,\div{\vec{u}}) &= 0\;,\\
        (\vec{\omega},\vec{\tau}) -(\vec{u},\curl{\vec{\tau}}) &= - (\hat{\vec{u}},\mathcal{T}_{\parallel}\vec{\tau})_{\Gamma_{\hat{\vec{u}}}}\;.
    \end{align}
\end{subequations}

\subsection{The de Rham complex}
For our discretization, we will seek solutions in finite-dimensional spaces which discretize the continuous de Rham complex \eqref{eq:deRhamComplex3d}.
This sequence is called a complex as each differential operator ($\grad{},\curl{},\div{}$) maps the previous function space into the next one and, additionally, the application of any two successive differential operators in the sequence is zero, i.e.,
\begin{equation}
    \div{\curl{\vec{\omega}}} = 0\;,~\forall \vec{\omega}\in H(\curl{};\Omega)\;,\quad \curl{\grad{\phi}} = \vec{0}\;,~\forall \phi\in H^1(\Omega)\;.
\end{equation}
This sequence encodes the geometric and topological properties of PDEs such as fluid flows and electromagnetism, and a stable and accurate discrete formulation for such PDEs can be constructed by using finite-dimensional spaces that form a structure-preserving subcomplex of \eqref{eq:deRhamComplex3d}, see \cite{arnold_finite_2018}.
Consider the finite-dimensional spaces $\mathbb{V}_h^0\subset H^1(\Omega),\mathbb{V}_{h}^1\subset H(\curl{};\Omega),\mathbb{V}_{h}^2\subset H(\div{};\Omega),\mathbb{V}_h^3\subset L^2(\Omega)$ such that they form such a structure-preserving subcomplex: 
\begin{equation}\label{eq:intro-2D-hilbert-complex-finite}
    0 \xrightarrow{} \mathbb{V}_h^0 \xrightarrow{\grad{}} \mathbb{V}_{h}^1 \xrightarrow{\curl{}} \mathbb{V}_{h}^2 \xrightarrow[]{\div{}}\mathbb{V}_h^3 \rightarrow 0\;.
\end{equation}
For example, one could use the complex formed by Lagrange elements of degree $\deg$, N\'ed\'elec (first kind) of degree $\deg-1$, the Raviart-Thomas elements of degree $\deg-1$, and discontinuous Lagrange elements of degree $\deg-1$ \cite{arnold_finite_2006}.
Alternatively, on (adaptively-refined) tensor product meshes, one could use multivariate B-spline spaces of varying polynomial degrees and regularities; see, for instance, \cite{buffa2011isogeometric,evans_divergence-free_2011,buffa2014isogeometric,evans2020hierarchical,shepherd2024locally}.

To simplify the formulation and analysis of the stabilized discrete formulation, we will also consider the problem with homogeneous boundary conditions.
The corresponding de Rham complex with homogeneous boundary conditions on $\partial \Omega$ is given by:
\begin{equation}\label{eq:intro-3D-hilbert-complex-homogeneous-bc}
        0 \xrightarrow{} H^1_0(\Omega,\partial \Omega) \xrightarrow{\grad{}} H_{0}(\curl{};\Omega,\partial \Omega) \xrightarrow{\curl{}} H_{0}(\div{};\Omega,\partial \Omega) \xrightarrow{\div{}} L^2_0(\Omega) \rightarrow 0\;,
\end{equation}
This complex can likewise be discretized using spaces $\mathbb{V}_{h,0}^0 = \mathbb{V}_h^0 \cap H^1_{}(\Omega,\partial \Omega)$, $\mathbb{V}_{h,\vec{0}}^1 = \mathbb{V}_h^1 \cap H_{\vec{0}}(\curl{};\Omega,\partial \Omega)$,
$\mathbb{V}_{h,0}^2 = \mathbb{V}_h^2 \cap H_{0}(\div{};\Omega,\partial \Omega)$, and $\mathbb{V}_{h,0}^3  = \mathbb{V}_h^3 \cap L^2_0(\Omega)$ that form a structure-preserving subcomplex:
\begin{equation}\label{eq:intro-2D-hilbert-complex-finite-hom-bc}
    0 \xrightarrow{} \mathbb{V}_{h,0}^0 \xrightarrow{\grad{}} \mathbb{V}_{h,\vec{0}}^1 \xrightarrow{\curl{}} \mathbb{V}_{h,0}^2 \xrightarrow[]{\div{}}\mathbb{V}_{h,0}^3 \rightarrow 0\;.
\end{equation}

\section{Variational Multiscale Stabilization via discontinuous fine-scales}\label{sec:VMS-method}
Proceeding with homogeneous boundary conditions, we now introduce a Variational Multiscale Stabilization (VMS) formulation for the incompressible Navier--Stokes equations in the vorticity-velocity-pressure form.
Within the framework of VMS, for a given projector $\Pi$:
\begin{equation*}
    \Pi : H_{\vec{0}}(\curl{};\Omega,\partial\Omega)\times H_{{0}}(\div{};\Omega,\partial\Omega)\times L^2_0(\Omega) \rightarrow \mathbb{V}_{h,\vec{0}}^1\times\mathbb{V}_{h,0}^2\times\mathbb{V}_{h,0}^3\;,
\end{equation*}
and the continuous solution $\left(\vec{\omega},\vec{u},p\right)\in H_{\vec{0}}(\curl{};\Omega,\partial\Omega)\times H_{{0}}(\div{};\Omega,\partial\Omega)\times L^2_0(\Omega)$, the goal is to develop a numerical method whose solution is the projected continuous solution, i.e., $\left(\vec{\omega}^h,\vec{u}^h,p^h\right) := \Pi \left(\vec{\omega},\vec{u},p\right)$.
For this, we consider the Stokes projector for $\Pi$, which is defined such that for all $(\vec{\tau}^h,\vec{v}^h,q^h)\in \mathbb{V}_{h,\vec{0}}^1\times\mathbb{V}_{h,0}^2\times\mathbb{V}_{h,0}^3$:
\begin{multline}\label{eq:stokes-projector}
      \mathrm{Re}^{-1}(\curl{\vec{\omega}^h},\vec{v}^h) - (p^h,\div{\vec{v}^h}) + (q^h,\div{\vec{u}^h}) + (\vec{\omega}^h,\vec{\tau}^h) - (\vec{u}^h,\curl{\vec{\tau}^h}) \\
     =\mathrm{Re}^{-1}(\curl{\vec{\omega}},\vec{v}^h) - (p,\div{\vec{v}^h}) + (q^h,\div{\vec{u}}) + (\vec{\omega},\vec{\tau}^h) - (\vec{u},\curl{\vec{\tau}^h})\;.
\end{multline}
This projector orthogonally decomposes any triple $(\vec{\tau},\vec{v},q)\in H_{\vec{0}}(\curl{};\Omega,\partial\Omega)\times H_{{0}}(\div{};\Omega,\partial\Omega)\times L^2_0(\Omega)$ into a finite-dimensional component $(\vec{\tau}^h,\vec{v}^h,q^h) = \Pi (\vec{\tau},\vec{v},q)$ and an infinite-dimensional component $(\vec{\tau}',\vec{v}',q') = (I - \Pi)(\vec{\tau},\vec{v},q)$, which we will refer to as the coarse-scale and fine-scale components.
We denote the orthogonal complement as:
\begin{equation}
    \left(\mathbb{V}_{h,\vec{0}}^1\times\mathbb{V}_{h,0}^2\times\mathbb{V}_{h,0}^3\right)^\perp := \left\{(\vec{\tau},\vec{v},q)\in H_{\vec{0}}(\curl{};\Omega,\partial\Omega)\times H_{{0}}(\div{};\Omega,\partial\Omega)\times L^2_0(\Omega)~:~\Pi(\vec{\tau},\vec{v},q)=0\right\}\;.
\end{equation}

From this orthogonal decomposition, the weak formulation \eqref{eq:intro-weak-NS} is split into a coarse-scale and a fine-scale problem, which must be solved simultaneously. 
The coarse-scale part of the problem is given as: find $(\vec{\omega}^h, \vec{u}^h,p^h) \in \mathbb{V}_{h,\vec{0}}^1\times\mathbb{V}_{h,0}^2\times\mathbb{V}_{h,0}^3$ such that for all $(\vec{\tau}^h,\vec{v}^h,q^h) \in \mathbb{V}_{h,\vec{0}}^1\times\mathbb{V}_{h,0}^2\times\mathbb{V}_{h,0}^3$:
\begin{subequations}\label{eq:coarse-oseen-equation}
\begin{align}
    \nonumber \left(\partial_t \vec{u}^h + \partial_t \vec{u}',\vec{v}^h\right) + (\left[\vec{\omega}^h+\vec{\omega}'\right]\times\left[ \vec{u}^h + \vec{u}'\right],\vec{v}^h)+& & \\
    \mathrm{Re}^{-1}(\curl{\vec{\omega}^h},\vec{v}^h) - (p^h,\div{\vec{v}^h}) 
    &= (\vec{f},\vec{v}^h) \;,\\
    (q^h,\div{\vec{u}^h}) &= 0\;,\\
    (\vec{\omega}^h,\vec{\tau}^h) - (\vec{u}^h,\curl{\vec{\tau}^h})
    &= 0\;.
\end{align}
\end{subequations}
Note that by the orthogonality of \eqref{eq:stokes-projector}, the impact of the fine scales in the above problem is solely expressed through time-evolution and the advection terms.

For the fine-scale part of the equations, we do not make use of the orthogonality of $\Pi$. This is because, to simplify their discretization, we will shortly introduce modelling assumptions on the structure of the space of fine-scales, and then the orthogonality may no longer hold.
Here, for all $(\vec{\tau}',\vec{v}',q') \in \left(\mathbb{V}_{h,\vec{0}}^1\times\mathbb{V}_{h,0}^2\times\mathbb{V}_{h,0}^3\right)^\perp$, we look for $(\vec{\omega}',\vec{u}',p') \in \left(\mathbb{V}_{h,\vec{0}}^1\times\mathbb{V}_{h,0}^2\times\mathbb{V}_{h,0}^3\right)^\perp$ such that:
\begin{subequations}\label{eq:fine-oseen-equation}
\begin{align}
    \nonumber \left(\partial_t \vec{u}^h + \partial_t \vec{u}',\vec{v}'\right) +(\left[\vec{\omega}^h+ \vec{\omega}'\right]\times\left[\vec{u}^h+\vec{u}'\right],\vec{v}') + &\\
    \mathrm{Re}^{-1}(\curl{\vec{\omega}^h+\vec{\omega}'},\vec{v}')-(p^h+p',\div{\vec{v}'})
    &= (\vec{f},\vec{v}')\;, \\
    (q',\div{\vec{u}'}) + (q',\div{\vec{u}^h}) &= 0\;,\\
    (\vec{\omega}',\vec{\tau}') - (\vec{u}',\curl{\vec{\tau}}') + (\vec{\omega}^h,\vec{\tau}') - (\vec{u}^h,\curl{\vec{\tau}}')
    &= 0\;.
\end{align}
\end{subequations}

\subsection{Modelling assumptions}
\label{sec:modelling-choices}
To develop a practical discretization method for the coupled system of coarse- and fine-scale problems introduced above, we make the following modelling assumptions:
\begin{itemize}
    \item We ignore the higher order fine-scale terms $(\vec{\omega}'\times\vec{u}',\vec{v}^h)$ and $(\vec{\omega}'\times\vec{u}',\vec{v}')$.
    \item We assume that $(q',\div{\vec{u}^h}) = 0$ and $(\vec{\omega}^h,\vec{\tau}') = (\vec{u}^h,\curl{\vec{\tau}'})$.
    Since our computed coarse-scale and fine-scale velocities will be pointwise divergence-free, the first assumption is automatically satisfied.
    \item We model the advection-diffusion terms in the fine-scale momentum balance as:
    \begin{align}\label{eq:advection-diffusion-modeling-choice}
        (\vec{\omega}'\times\vec{u}^h,v') + \mathrm{Re}^{-1}(\curl{\vec{\omega}'},\vec{v}') &\approx (\tm  \vec{u}',\vec{v}') + \frac{\mathrm{Re}^{-1}}{2}(\curl{\vec{\omega}'},\vec{v}')\;,
    \end{align}
    where $\tm$ is the stabilisation parameter.
    \item For energy stability, inspired by \cite{ten_eikelder_correct_2018}, we add the term \(\left( \vec{u}', \curl{\vec{\tau}^h}\right)\) to the coarse-scale vorticity equation.
    \item We assume that the fine-scale function space is the product of the following bubble-function spaces, which form a de Rham complex on each mesh element:
    \begin{equation}
        \bigcup_{e\in\mathcal{M}} H_{\vec{0}}(\curl{};\Omega^e,\partial \Omega^e) \xrightarrow[]{\curl{}} \bigcup_{e\in\mathcal{M}} H_0(\div{};\Omega^e,\partial\Omega^e) \xrightarrow[]{\div{}} \bigcup_{e\in\mathcal{M}} L^2_{0}(\Omega^e)\;,
    \end{equation}
    and where $\mathcal{M} =\{e~:~\Omega^e\text{ is a mesh element}\}$ is the set of all mesh elements.
    
\end{itemize}

\subsection{Stabilized semi-discrete formulation}
Given the above assumptions, we will discretize the coarse-scale problem by choosing the finite-dimensional spaces $\mathbb{V}_{h,\vec{0}}^1$, $\mathbb{V}_{h,0}^2$, and $\mathbb{V}_{h,0}^3$ introduced in \eqref{eq:intro-2D-hilbert-complex-finite-hom-bc}.
For the fine-scale problem, we will discretize the modified fine-scale product-space with finite-dimensional bubble-function spaces which form a structure-preserving subcomplex on each mesh element:
\begin{subequations}
    \begin{align*}
        \mathbb{W}_{h,\vec{0}}^1&\subset\bigcup_{e\in\mathcal{M}} H_{\vec{0}}(\curl{};\Omega^e,\partial \Omega^e)\;, &
        \mathbb{W}_{h,0}^2&\subset\bigcup_{e\in\mathcal{M}} H_0(\div{};\Omega^e,\partial \Omega^e)\;, &
        \mathbb{W}_{h,0}^3&\subset\bigcup_{e\in\mathcal{M}} L^2_0(\Omega^e)\;.
    \end{align*}
\end{subequations}  
Both the coarse- and fine-scale spaces will be chosen as piecewise-polynomial function spaces defined on the same partition, and the degree of the fine-scale spaces will be selected to be at least as large as the degree of the coarse-scale space.
To simplify the notation, we will define the product space:
\begin{equation}
    \mathbb{S}_h := \mathbb{V}_{h,\vec{0}}^1\times\mathbb{V}_{h,0}^2\times\mathbb{V}_{h,0}^3\times \mathbb{W}_{h,\vec{0}}^1\times\mathbb{W}_{h,0}^2\times\mathbb{W}_{h,0}^3\;.
\end{equation}
The following box presents our stabilized semi-discrete formulation for the incompressible Navier--Stokes equations in the velocity-vorticity-pressure form.
\begin{methodBox}\label{method:sd-ns-stabilization}
Find $(\vec{\omega}^h,\vec{u}^h,p^h,\vec{\omega}',\vec{u}',p') \in\mathbb{S}_h$ such that, for all $(\vec{\tau}^h,\vec{v}^h,q^h,\vec{\tau}',\vec{v}',q') \in\mathbb{S}_{h}$:
\begin{subequations}\label{eq:sd-ns-stabilized-formulation}
    \begin{align}
        \nonumber \left(\partial_t\vec{u}^h + \partial_t\vec{u}',\vec{v}^h\right) + (\vec{\omega}^h\times\left[\vec{u}^h+\vec{u}'\right],\vec{v}^h)+ &\\
    (\vec{\omega}'\times\vec{u}^h,\vec{v}^h)+ \mathrm{Re}^{-1}(\curl{\vec{\omega}^h},\vec{v}^h) - (p^h,\div{\vec{v}^h}) 
    &= (\vec{f},\vec{v}^h)\;,\\
    (q^h,\div{\vec{u}^h}) &= 0\;,\\
    (\vec{\omega}^h,\vec{\tau}^h) - (\vec{u}^h,\curl{\vec{\tau}^h}) - (\vec{u}',\curl{\vec{\tau}^h})
    &= 0\;,\\
    \nonumber \left(\partial_t\vec{u}^h + \partial_t\vec{u}',\vec{v}'\right) + (\vec{\omega}^h\times \vec{u}',\vec{v}')+(\tm \vec{u}',\vec{v}') +& \\
     \frac{\mathrm{Re}^{-1}}{2}( \curl{\vec{\omega}'},{\vec{v}'}) - ({p'} ,\div{\vec{v}'})-(p^h,\div{\vec{v}'}) &= (\mathrm{res_M},\vec{v}')\\
    ({q'},\div{\vec{u}'})  &= 0\;,\\
    (\vec{\omega}',\vec{\tau}') - ({\vec{u}'},\curl{\vec{\tau}'})  &= 0\;.
    \end{align}
\end{subequations}
where the momentum residual is defined as:
\begin{equation}
    \mathrm{res_M} := \vec{f} - \vec{\omega}^{h}\times \vec{u}^{h} - \mathrm{Re}^{-1}\curl{\vec{\omega}^{h}}.
\end{equation}
\end{methodBox}

\begin{remark}\label{rem:fine-scale-boundary-conditions}
    By our choice of fine-scale bubble spaces, we impose homogeneous boundary conditions $\vec{n}\cdot\vec{u}'=0$ and $\vec{\omega}'=\vec{0}$ on the entire mesh skeleton.
    Consequently, the fine-scale spaces also form a discrete de Rham complex on the entire domain $\Omega$:
    \begin{equation}
        0 \xrightarrow[]{} \mathbb{W}_{h,\vec{0}}^1 \xrightarrow[]{\curl{}} \mathbb{W}_{h,0}^2 \xrightarrow[]{\div{}} \mathbb{W}_{h,0}^3 \xrightarrow[]{} 0\;,
    \end{equation}
\end{remark}

\begin{remark}\label{rem:elementwise-fine-scales}
    Despite the way that the fine-scale problem is presented in \eqref{eq:sd-ns-stabilized-formulation}, due to our choice of fine-scale spaces, this problem decouples into elementwise problems which can be solved in parallel.
    That is, for each mesh element $\Omega^e$, we look for $(\vec{\omega}',\vec{u}',p') \in \mathbb{W}_{h,\vec{0}}^1(\Omega^e)\times\mathbb{W}_{h,0}^2(\Omega^e)\times\mathbb{W}_{h,0}^3(\Omega^e)$ such that, for all $(\vec{\tau}',\vec{v}',q') \in \mathbb{W}_{h,\vec{0}}^1(\Omega^e)\times\mathbb{W}_{h,0}^2(\Omega^e)\times\mathbb{W}_{h,0}^3(\Omega^e)$:
    \begin{subequations}\label{eq:sd-formulation-fine-scale}
    \begin{align}
        \nonumber \left(\partial_t\vec{u}^h + \partial_t\vec{u}',\vec{v}'\right)_{\Omega^e} + (\vec{\omega}^h\times \vec{u}',\vec{v}')_{\Omega^e}+(\tm \vec{u}',\vec{v}')_{\Omega^e} +& \\
         \frac{\mathrm{Re}^{-1}}{2}( \curl{\vec{\omega}'},{\vec{v}'})_{\Omega^e} - ({p'} ,\div{\vec{v}'})_{\Omega^e}-(p^h,\div{\vec{v}'})_{\Omega^e} &= (\mathrm{res_M},\vec{v}')_{\Omega^e}\;,\\
        ({q'},\div{\vec{u}'})_{\Omega^e}  &= 0\;,\\
        (\vec{\omega}',\vec{\tau}')_{\Omega^e} - ({\vec{u}'},\curl{\vec{\tau}'})_{\Omega^e}  &= 0\;.
    \end{align}
\end{subequations}
\end{remark}

\subsection{Energetically-stable fully-discrete formulation}
To construct an energetically-stable method, we temporally discretize Formulation \ref{method:sd-ns-stabilization} with Crank--Nicolson time integration. Specifically, for time step $\Delta t > 0$, we evaluate each quantity at $N = T/\Delta t$ equally-spaced time point and use the following approximations:
\begin{equation}
    \partial_t[\cdot] := \frac{[\cdot]^{n+1} - [\cdot]^n}{\Delta t}\;,\quad [\cdot]^{n+\nicefrac{1}{2}} := \frac{[\cdot]^{n+1} + [\cdot]^n}{2}\;.
\end{equation}
To address the nonlinearity in the Navier--Stokes equations at time step \(n\), we use Picard iterations. We denote quantities at the \( \PicIt \) iteration with a superscript of \( \PicIt \) and linearize the convection term, leading to the following formulation:
\begin{methodBox}\label{method:ns-stabilization}
Given the solution at time step $n$, $(\vec{\omega}^{h,n},\vec{u}^{h,n},p^{h,n},\vec{\omega}'^{n},\vec{u}'^{n},p'^{n}) \in\mathbb{S}_h$ and the $\PicIt$-th Picard iterate $(\vec{\omega}^{h,\PicIt,n+1},\vec{u}^{h,\PicIt,n+1},p^{h,\PicIt,n+1},\vec{\omega}'^{\PicIt,n+1},\vec{u}'^{\PicIt,n+1},p'^{\PicIt,n+1}) \in\mathbb{S}_h$, find the next Picard iterate $(\vec{\omega}^{h,\PicIt+1,n+1},\vec{u}^{h,\PicIt+1,n+1},p^{h,\PicIt+1,n+1},\vec{\omega}'^{\PicIt+1,n+1},\vec{u}'^{\PicIt+1,n+1},p'^{\PicIt+1,n+1}) \in\mathbb{S}_h$, such that, for all $(\vec{\tau}^h,\vec{v}^h,q^h,\vec{\tau}',\vec{v}',q') \in\mathbb{S}_{h}$:
\begin{subequations}
    \begin{align}
        \nonumber \left(\frac{\vec{u}^{h,\PicIt+1,n+1}-\vec{u}^{h,n}}{\Delta t}+\frac{\vec{u}'^{\PicIt+1,n+1}-\vec{u}'^{n}}{\Delta t},\vec{v}^h\right) + \\
        \nonumber (\vec{\omega}^{h,\PicIt+1,n+\nicefrac{1}{2}}\times\left[\vec{u}^{h,\PicIt,n+\nicefrac{1}{2}}+
        \vec{u}'^{\PicIt,n+\nicefrac{1}{2}}\right],\vec{v}^h)+&\\
        \nonumber (\vec{\omega}'^{\PicIt+1,n+\nicefrac{1}{2}}\times\vec{u}^{h,\PicIt,n+\nicefrac{1}{2}},\vec{v}^h) 
     +     \mathrm{Re}^{-1}(\curl{\vec{\omega}^{h,\PicIt+1,n+\nicefrac{1}{2}}},\vec{v}^h) - &\\
    (p^{h,\PicIt+1,n+\nicefrac{1}{2}},\div{\vec{v}^h})
    &= 
    (\vec{f},\vec{v}^h)\;,\\
    (q^h,\div{\vec{u}^{h,\PicIt+1,n+\nicefrac{1}{2}}}) &= 0\;,\\
    (\vec{\omega}^{h,\PicIt+1,n+\nicefrac{1}{2}},\vec{\tau}^h) - (\vec{u}^{h,\PicIt+1,n+\nicefrac{1}{2}}
    +\vec{u}'^{\PicIt+1,n+\nicefrac{1}{2}},\curl{\vec{\tau}^h})
    &= 0\;,\\
    \nonumber \left(\frac{\vec{u}^{h,\PicIt+1,n+1}-\vec{u}^{h,n}}{\Delta t}+\frac{\vec{u}'^{\PicIt+1,n+1}-\vec{u}'^{n}}{\Delta t},\vec{v}'\right) +\\
    \nonumber (\vec{\omega}^{h,\PicIt,n+\nicefrac{1}{2}}\times \vec{u}'^{\PicIt+1,n+\nicefrac{1}{2}},\vec{v}') 
     +(\tm \vec{u}'^{\PicIt+1,n+\nicefrac{1}{2}},\vec{v}') +&\\
    \nonumber + \frac{\mathrm{Re}^{-1}}{2}(\curl{\vec{\omega}'^{\PicIt+1,n+\nicefrac{1}{2}}},{\vec{v}'}) 
       - ({p'^{\PicIt+1,n+\nicefrac{1}{2}}} ,\div{\vec{v}'}) -&\\
       ({p^{h,\PicIt+1,n+\nicefrac{1}{2}}} ,\div{\vec{v}'}) &= (\mathrm{res_M},\vec{v}')\;,\\
    ({q'},\div{\vec{u}'^{\PicIt+1,n+\nicefrac{1}{2}}})  &= 0\;,\\
    (\vec{\omega}'^{\PicIt+1,n+\nicefrac{1}{2}},\vec{\tau}') - ({\vec{u}'^{\PicIt+1,n+\nicefrac{1}{2}}},\curl{\vec{\tau}'})  &= 0\;,
    \end{align}
\end{subequations}
where
\begin{align}
    \mathrm{res_M} &:= \vec{f} - \vec{\omega}^{h,\PicIt+1,n+\nicefrac{1}{2}}\times \vec{u}^{h,\PicIt,n+\nicefrac{1}{2}} - \mathrm{Re}^{-1}\curl{\vec{\omega}^{h,\PicIt+1,n+\nicefrac{1}{2}}}\;.
\end{align}
The solution at time step $n+1$ is defined as the limit of the Picard iterations $\PicIt\rightarrow\infty$.
\end{methodBox}
This method is energetically stable, as the following result shows:
\begin{lemma}\label{lem:energy-evolution}
    The solutions to Method \ref{method:ns-stabilization}, evolve over time as:
    \begin{multline}\label{eq:ns-stabilization-energy-estimate}
    \frac{1}{\Delta t}\left( \frac{1}{2}\Vert \vec{u}^{h,n+1}  + \vec{u}'^{n+1}\Vert^2 - \frac{1}{2} \Vert \vec{u}^{h,n} + \vec{u}'^n \Vert^2 \right) = \\
    (\vec{f},\vec{u}^{h,n+\nicefrac{1}{2}}+\vec{u}'^{n+\nicefrac{1}{2}}) - {\mathrm{Re}^{-1}}\Vert \vec{\omega}^{h,n+\frac{1}{2}}\Vert^2 - \frac{\mathrm{Re}^{-1}}{2}\Vert \vec{\omega}'^{n+\frac{1}{2}}\Vert^2 - \| \sqrt{\tm}\vec{u}'^{n+\nicefrac{1}{2}} \|^2\;.
\end{multline}
\end{lemma}
\begin{proof}
We remark that, by our choice of fine-scale spaces, for any $p^h\in\mathbb{V}_h^3$,
\begin{equation*}
    \hat{p}^h := p^h - \sum_{e\in\mathcal{M}} \left(\int_{\Omega^e} p^h \mathrm{d}V\right) \mathbf{1}_{\Omega^e} \in \mathbb{W}_{h,0}^3\;,
\end{equation*}
where $\mathbf{1}_D$ is the indicator function on $D$.
As a result, for any $v'\in\mathbb{W}_{h,0}^2$,
\begin{equation*}
    (p^h,\div{v'}) = (\hat{p}^h,\div{v'})\;.
\end{equation*}
Then, make the following choices for the test functions in Formulation \ref{method:ns-stabilization}:
\begin{equation*}
    \begin{split}
        \vec{v}^h = \vec{u}^{h,n+\nicefrac{1}{2}},\quad q^h = p^{h,n+\nicefrac{1}{2}},\quad \vec{\tau}^h = \mathrm{Re}^{-1}\vec{\omega}^{h,n+\nicefrac{1}{2}}\;,\\
        \vec{v}' = \vec{u}'^{n+\nicefrac{1}{2}},\quad q' = p'^{n+\nicefrac{1}{2}} + \hat{p}^{h,n+\nicefrac{1}{2}},\quad \vec{\tau}' = \frac{1}{2}\mathrm{Re}^{-1}\vec{\omega}'^{n+\nicefrac{1}{2}}\;.
    \end{split}
\end{equation*}
Summing up all the terms in Formulation \ref{method:ns-stabilization},
and noting that
\begin{multline*}
\left(\frac{\vec{u}^{h,n+1}-\vec{u}^{h,n}}{\Delta t}+\frac{\vec{u}'^{n+1}-\vec{u}'^{n}}{\Delta t},\frac{\vec{u}^{h,n+1}+\vec{u}^{h,n}}{2}\right)+ \left(\frac{\vec{u}^{h,n+1}-\vec{u}^{h,n}}{\Delta t}+\frac{\vec{u}'^{n+1}-\vec{u}'^{n}}{\Delta t},\frac{\vec{u}'^{n+1}+\vec{u}'^{n}}{2}\right)\\=
\frac{1}{2\Delta t}\left(\Vert \vec{u}^{h,n+1} + \vec{u}'^{n+1}\Vert^2 - \Vert \vec{u}^{h,n} + \vec{u}'^{n}\Vert^2 \right)\;,
\end{multline*}
we obtain the energy estimate \eqref{eq:ns-stabilization-energy-estimate}.
\end{proof}

\begin{remark}
    In the modelling section, we chose to discard the terms $(\vec{\omega}^h,\vec{\tau}') - (\vec{u}^h,\curl{\vec{\tau}'})$ from the fine-scale vorticity equation.
    This assumption can be omitted without spoiling the energy stability of the formulation.
    However, it will result in the energy evolution of Lemma \ref{lem:energy-evolution} being transformed from an equality to an upper bound.
    This was a deliberate choice made by the authors. 
\end{remark}

\subsection{General boundary conditions}
While our theoretical analysis will only focus on the case of homogeneous boundary conditions, the method can easily accommodate general boundary conditions of the form presented in \eqref{eq:intro-boundary-conditions}.
We can do this by replacing the coarse-scale spaces with
\begin{align*}
    \mathbb{V}_{h,\Gamma_{\hat{\vec{\omega}}},\hat{\vec{\omega}}}^1&= \mathbb{V}_h^1 \cap H_{\hat{\vec{\omega}}}(\curl{};\Omega,\Gamma_{\hat{\vec{\omega}}})\;,\\
    \mathbb{V}_{h,\Gamma_{\hat{u}},\hat{u}}^2&= \mathbb{V}_h^2 \cap H_{\hat{u}}(\div{};\Omega,\Gamma_{\hat{u}})\;,
\end{align*}
and by altering Formulation \ref{method:ns-stabilization} with the appropriate boundary terms to yield the following semi-discrete formulation.
\begin{methodBox}\label{method:ns-stabilization-gen-bc}
Find $(\vec{\omega}^h,\vec{u}^h,p^h,\vec{\omega}',\vec{u}',p') \in \mathbb{V}_{h,\Gamma_{\hat{\vec{\omega}}},\hat{\vec{\omega}}}^1\times\mathbb{V}_{h,\Gamma_{\hat{u}},\hat{u}}^2\times\mathbb{V}_{h,0}^3\times \mathbb{W}_{h,\vec{0}}^1\times\mathbb{W}_{h,0}^2\times\mathbb{W}_{h,0}^3$, such that, for all $(\vec{\tau}^h,\vec{v}^h,q^h,\vec{\tau}',\vec{v}',q') \in \mathbb{V}_{h,\Gamma_{\hat{\vec{\omega}}},\vec{0}}^1\times\mathbb{V}_{h,\Gamma_{\hat{u}},0}^2\times\mathbb{V}_{h,0}^3\times \mathbb{W}_{h,\vec{0}}^1\times\mathbb{W}_{h,0}^2\times\mathbb{W}_{h,0}^3
$:
\begin{subequations}\label{eq:sd-ns-stabilized-formulation-general-bc}
    \begin{align}
        \nonumber \left(\partial_t\vec{u}^h + \partial_t\vec{u}',\vec{v}^h\right) + (\vec{\omega}^h\times\left[\vec{u}^h+\vec{u}'\right],\vec{v}^h)+&\\
    \label{eq:sd-ns-stabilized-formulation-coarse-momentum} (\vec{\omega}'\times\vec{u}^h,\vec{v}^h) +\mathrm{Re}^{-1}(\curl{\vec{\omega}^h},\vec{v}^h) - (p^h,\div{\vec{v}^h})  
    &= (\vec{f},\vec{v}^h) - (\hat{p},\mathcal{T}_{\perp}\vec{v}^h)_{\Gamma_{\hat{p}}}\;,\\
    (q^h,\div{\vec{u}^h}) &= 0\;,\\
    (\vec{\omega}^h,\vec{\tau}^h) - (\vec{u}^h,\curl{\vec{\tau}^h}) - (\vec{u}',\curl{\vec{\tau}^h})
    &= -(\hat{\vec{u}},\mathcal{T}_{\parallel}\vec{\tau}^h)_{\Gamma_{\hat{\vec{u}}}} \;,\\
    \nonumber \left(\partial_t\vec{u}^h + \partial_t\vec{u}',\vec{v}'\right) + (\vec{\omega}^h\times \vec{u}',\vec{v}')+(\tm \vec{u}',\vec{v}') +& \\
    \label{eq:sd-formulation-fine-scale-momentum} \frac{\mathrm{Re}^{-1}}{2}( \curl{\vec{\omega}'},{\vec{v}'}) - ({p'} ,\div{\vec{v}'})- ({p}^h ,\div{\vec{v}'}) &= (\mathrm{res_M},\vec{v}')\\
    \label{eq:sd-formulation-fine-scale-mass}({q'},\div{\vec{u}'})  &= 0\;,\\
    \label{eq:sd-formulation-fine-scale-vorticity}(\vec{\omega}',\vec{\tau}') - ({\vec{u}'},\curl{\vec{\tau}'})  &= 0\;,
    \end{align}
\end{subequations}
 where the momentum residual is defined as:
\begin{equation}
    \mathrm{res_M} := \vec{f} - \vec{\omega}^{h}\times \vec{u}^{h} - \mathrm{Re}^{-1}\curl{\vec{\omega}^{h}}.
\end{equation}
\end{methodBox}

\section{Stability, uniqueness and convergence of the Oseen equations}
\label{sec:convergence}
To demonstrate the convergence of the proposed method, we analyse the time-dependent linearised Navier--Stokes equations, known as the Oseen equations again in the setting of homogeneous boundary conditions and the same modelling assumptions.
The weak formulation in velocity-vorticity-pressure variables reads: find $(\vec{\omega},\vec{u},p)\in H_0(\curl{}; \Omega,\partial\Omega)\times H_0(\div{};\Omega,\partial \Omega)\times L^2_0(\Omega)$, such that for all $(\vec{\tau},\vec{v},q)\in H_0(\curl{}; \Omega,\partial\Omega)\times H_0(\div{};\Omega,\partial \Omega)\times L^2_0(\Omega)$,
\begin{subequations}\label{eq:weak-osseen}
\begin{align}
    \label{eq:weak-osseen-momentum}(\sigma\,\vec{u},\,\vec{v}) + \Big(\frac{1}{\sqrt{\nu}}\,\vec{\omega}\times\vec{\beta},\,\vec{v}\Big)
    + \sqrt{\nu}\,(\curl{\vec{\omega}},\,\vec{v}) - (p,\,\div{\vec{v}}) &= (\vec{f},\,\vec{v})\;,\\
    (\div{\vec{u}},\,q) &= 0\;,\\
    (\vec{\omega},\,\vec{\tau}) - \sqrt{\nu}\,(\vec{u},\,\curl{\vec{\tau}}) &= 0\;,
\end{align}
\end{subequations}
where $\nu,\sigma>0$ are positive constants and $\vec{\beta}\in L^{\infty}$ is a known solenoidal advecting velocity field.
We assume that the velocity field is nowhere zero to simplify the analysis.
Furthermore, we also assumed that this velocity field can be decomposed into coarse- and fine-scale contributions using some projection, i.e. $\vec{\beta} = \vec{\beta}^h + \vec{\beta}'$ where $\vec{\beta}^h,\vec{\beta}'\in L^{\infty}$ and $\vec{\beta}^h \in \mathbb{V}_{h,0}^2$.
Clearly, in the incompressible Navier--Stokes case, the velocity field has such a coarse- and fine-scale decomposition induced by the Stokes projector, but the analysis can be done in this more general setting. 
The stabilized discrete weak form that we will analyze is presented in Formulation \ref{method:oseen-stabilization} below. 
\begin{methodBox}\label{method:oseen-stabilization}
Find $(\vec{\omega}^h,\vec{u}^h,p^h,\vec{\omega}',\vec{u}',p') \in\mathbb{S}_h$
, such that, for all $(\vec{\tau}^h,\vec{v}^h,q^h,\vec{\tau}',\vec{v}',q') \in\mathbb{S}_{h}$:
\begin{subequations}\label{eq:oseen-stabilized-formulation}
    \begin{align}
        \nonumber \sigma(\vec{u}^h+\vec{u}',\vec{v}^h)+ \frac{1}{\sqrt{\nu}}(\vec{\omega}^h\times\left[\vec{\beta}^h+\vec{\beta}'\right],\vec{v}^h)+ \frac{1}{\sqrt{\nu}}(\vec{\omega}'\times\vec{\beta}^h,\vec{v}^h)+ &\\
    \label{eq:oseen-stabilized-formulation-coarse-momentum} \sqrt{\nu}(\curl{\vec{\omega}^h},\vec{v}^h) - (p^h,\div{\vec{v}^h}) &= (\vec{f},\vec{v}^h)\;,\\
    (q^h,\div{\vec{u}^h}) &= 0\;,\\
    (\vec{\omega}^h,\vec{\tau}^h) - \sqrt{\nu}(\vec{u}^h,\curl{\vec{\tau}^h}) - \sqrt{\nu}(\vec{u}',\curl{\vec{\tau}^h}) &= 0\;,\\
    \nonumber \sigma(\vec{u}',\vec{v}') +  (\tm \vec{u}',\vec{v}') + \sqrt{\frac{\nu}{2}}( \curl{\vec{\omega}'},{\vec{v}'}) -
    (p',\div{\vec{v}'}) + &\\
     \label{eq:oseen-stabilized-formulation-fine-momentum}\left( \sigma \vec{u}^h + \frac{1}{\sqrt{\nu}}\vec{\omega}^h\times\left[\vec{\beta}^h + \vec{\beta}'\right] + \sqrt{\nu}\curl{\vec{\omega}^h},\vec{v}'\right)- (p^h,\div{\vec{v}'}) &= (\vec{f},\vec{v}')\;,\\
    ({q'},\div{\vec{u}'})  &= 0\;,\\
    (\vec{\omega}',\vec{\tau}') - \sqrt{\frac{\nu}{2}}(\curl{\vec{u}'},{\vec{\tau}'})  &= 0\;.
    \end{align}
\end{subequations}
\end{methodBox}
\subsection{Stability and uniqueness}
To demonstrate that Formulation \ref{method:oseen-stabilization} permits a stable and unique solution, we extend the argument presented in \cite{anaya_analysis_2019} to our proposed method.
For this, we require the following inf-sup conditions for the chosen finite element spaces. In particular, when the chosen spaces form a structure-preserving discretization of the de Rham complex, these conditions are automatically satisfied \cite{arnold_finite_2018}.
\begin{assumpBox}\label{ass:inf-sup}
    For the finite element spaces, the following inf-sup conditions hold:
    \begin{align}
        \label{eq:infsupVelPresCoarse}\sup_{\substack{\vec{v}^h\in \mathbb{V}_{h,0}^2\\\vec{v}^h\neq 0}}\frac{\left( q^h,\div{\vec{v}^h} \right)}{\Vert \vec{v}^h \Vert_{H(\div{})}}\geq \beta_1 \Vert q^h\Vert\;,\quad &\forall q^h\in\mathbb{V}_{h,0}^3\;,\\
        \label{eq:infsupVelPresFine}\sup_{\substack{\vec{v}'\in \mathbb{W}_{h,0}^2\\\vec{v}'\neq 0}}\frac{\left(  q',\div{\vec{v}'} \right)}{\Vert \vec{v}' \Vert_{H(\div{})}}\geq \beta_2 \Vert q'\Vert\;,\quad &\forall q'\in\mathbb{W}_{h,0}^3\;.
    \end{align}
\end{assumpBox}
This allows us to define the kernel spaces:
\begin{align}
    H^{\mathrm{ker}}_0(\div{};\Omega) &:= \left\{~\vec{v}\in H_0(\div{};\Omega)~:~(q,\div{\vec{v}})=0~\forall q\in L^2_0(\Omega)~\right\},&\quad \\
    \mathbb{V}_{h,0}^{2,\mathrm{ker}}&:=\mathbb{V}_{h,0}^2 \cap H_0^{\mathrm{ker}}(\div{};\Omega,\partial \Omega)\;,\\
    \mathbb{W}_{h,0}^{2,\mathrm{ker}}&:=\mathbb{W}_{h,0}^2 \cap \bigcup_{e\in\mathcal{M}} H_0^{\mathrm{ker}}(\div{};\Omega^e,\partial \Omega^e)\;.
\end{align}
Using these kernel spaces, we define the following function spaces and a stabilisation norm:
\begin{align}
\mathbb{X} &:= H_{\vec{0}}(\curl{};\Omega,\partial\Omega)\times H_0^{\mathrm{ker}}(\div{};\Omega,\partial\Omega) \times \bigcup_{e\in\mathcal{M}} H_{\vec{0}}(\curl{};\Omega^e,\partial \Omega^e) \times\bigcup_{e\in\mathcal{M}} H_0^{\mathrm{ker}}(\div{};\Omega^e,\partial \Omega^e)\\
    \mathbb{X}_h&:=\mathbb{V}_{h,\vec{0}}^1\times\mathbb{V}_{h,0}^{2,\mathrm{ker}}\times\mathbb{W}_{h,\vec{0}}^1\times\mathbb{W}_{h,0}^{2,\mathrm{ker}} \subset \mathbb{X} \\
    \Vert x^h\Vert_{\mathbb{X}_h}^2 &:= \Vert \vec{\omega}^h\Vert^2 + \nu\Vert \curl{\vec{\omega}^h}\Vert^2 + \Vert\vec{\omega}'\Vert^2+ \nu\Vert \curl{\vec{\omega}'}\Vert^2+ \Vert \vec{u}^h+ \vec{u}'\Vert^2 +  \Vert \sqrt{\tm} \vec{u}' \Vert^2\;.
\end{align}
With these spaces, a reduced version of Formulation \ref{method:oseen-stabilization} is as below. 
\begin{methodBox}\label{method:reduced-oseen-stabilization}
Find $(\vec{\omega}^h,\vec{u}^h,\vec{\omega}',\vec{u}') \in\mathbb{X}_h $, such that, for all $(\vec{\tau}^h,\vec{v}^h,\vec{\tau}',\vec{v}') \in\mathbb{X}_{h}$:
\begin{subequations}
\begin{align}
        \nonumber \sigma(\vec{u}^h+\vec{u}',\vec{v}^h)+\frac{1}{\sqrt{\nu}}(\vec{\omega}^h\times\left[\vec{\beta}^h+\vec{\beta}'\right],\vec{v}^h)+&\\ 
        \frac{1}{\sqrt{\nu}}(\vec{\omega}'\times\vec{\beta}^h,\vec{v}^h) 
        +\sqrt{\nu}(\curl{\vec{\omega}^h},\vec{v}^h) &= (\vec{f},\vec{v}^h)\;,\\
     (\vec{\omega}^h,\vec{\tau}^h) - \sqrt{\nu}(\vec{u}^h,\curl{\vec{\tau}^h}) - \sqrt{\nu}(\vec{u}',\curl{\vec{\tau}^h}) &= 0\;,\\
    \nonumber \sigma(\vec{u}',\vec{v}') + (\tm \vec{u}',\vec{v}') + \sqrt{\frac{\nu}{2}}(\curl{\vec{\omega}'},{\vec{v}'})+ &\\
  \left( \sigma \vec{u}^h+ \frac{1}{\sqrt{\nu}}\vec{\omega}^h\times\left[\vec{\beta}^h + \vec{\beta}'\right] + \sqrt{\nu}\curl{\vec{\omega}^h},\vec{v}'\right) &= (\vec{f},\vec{v}')\;,\\
    \label{eq:reduced-oseen-stabilized-formulation}(\vec{\omega}',\vec{\tau}') - \sqrt{\frac{\nu}{2}}({\vec{u}'},\curl{\vec{\tau}'})  &= 0\;.
\end{align}
\end{subequations}
\end{methodBox}
The equivalence of Formulations \ref{method:oseen-stabilization} and \ref{method:reduced-oseen-stabilization} is demonstrated by the following result.
\begin{lemma}\label{lem:equivalence-full-reduced-form}
    If $(\vec{\omega}^h,\vec{u}^h,p^h,\vec{\omega}',\vec{u}',p')$ is a solution to Formulation \ref{method:oseen-stabilization}, $(\vec{\omega}^h,\vec{u}^h,\vec{\omega}',\vec{u}')$ is a solution to Formulation \ref{method:reduced-oseen-stabilization}. If $(\vec{\omega}^h,\vec{u}^h,\vec{\omega}',\vec{u}')$ is a solution to Formulation \ref{method:reduced-oseen-stabilization}, then there exists $(p^h,p')\in \mathbb{V}_{h,0}^3\times\mathbb{W}_{h,0}^3$ so that $(\vec{\omega}^h,\vec{u}^h,p^h,\vec{\omega}',p',\vec{u}')$ is a solution to Formulation \ref{method:oseen-stabilization}.
\end{lemma}
\begin{proof}
    The first result is evident from from Assumption \ref{ass:inf-sup} and $\div{\mathbb{V}}_{h,0}^2\subseteq \mathbb{V}_{h,0}^3$, so that the solution $(\vec{\omega}^h,\vec{u}^h,p^h,\vec{\omega}',\vec{u}',p')$ to Formulation \ref{method:oseen-stabilization}, satisfies $\vec{u}^h\in\mathbb{V}_{h,0}^{2,\mathrm{ker}}$ and $\vec{u}'\in\mathbb{W}_{h,0}^{2,\mathrm{ker}}$. Likewise, for any $p^h\in\mathbb{V}_h^3$,
\begin{equation*}
    \hat{p}^h := p^h - \sum_{e\in\mathcal{M}} \left(\int_{\Omega^e} p^h \mathrm{d}V\right) \mathbf{1}_{\Omega^e} \in \mathbb{W}_{h,0}^3\;,
\end{equation*}
where $\mathbf{1}_D$ is the indicator function on $D$. As a result, for any $v'\in\mathbb{W}_{h,0}^2$,
\begin{equation*}
    (p^h,\div{v'}) = (\hat{p}^h,\div{v'})\;.
\end{equation*}
    For the converse, let $(\vec{\omega}^h,\vec{u}^h,\vec{\omega}',\vec{u}')$ be the solution to Formulation \ref{method:reduced-oseen-stabilization}.
    Then, we first solve for $p^h$ via equation \eqref{eq:oseen-stabilized-formulation-coarse-momentum}, for which a unique solution exists by Babuška–Lax–Milgram and \eqref{eq:infsupVelPresCoarse}. 
    Likewise, we solve for $p'$ from equation \eqref{eq:oseen-stabilized-formulation-fine-momentum}, whose existence and uniqueness rely on Babuška–Lax–Milgram and \eqref{eq:infsupVelPresFine}.
\end{proof}
Let $\mathcal{A}_{\mathrm{red}}$ be the bilinear form and $ \mathcal{F}_{\mathrm{red}}$ the linear form so that the variational form of Formulation \ref{method:reduced-oseen-stabilization} is equal to finding $x^h\in\mathbb{X}_h$, for which,
\begin{equation}\label{eq:compact-reduced-oseen-stabilization}
    \mathcal{A}_{\mathrm{red}}(x^h,y^h) = \mathcal{F}_{\mathrm{red}}(y^h),\quad \forall y^h\in\mathbb{X}_h\;.
\end{equation}
To show stability and uniqueness of Formulation \ref{method:reduced-oseen-stabilization} and \eqref{eq:compact-reduced-oseen-stabilization}, we make use of the following general result \cite[Theorem 1.2]{gatica_simple_2014}, and we impose Assumption \ref{ass:CFL-conditions} presented below.
\begin{theorem}\label{thm:stability}
    Let $\mathcal{A} : \mathcal{X}\times\mathcal{X}\rightarrow\mathbb{R}$ be a bounded bilinear form and $\mathcal{F}:\mathcal{X}\rightarrow\mathbb{R}$ a bounded functional, both defined on the Hilbert space $(\mathcal{X},\left<\cdot,\cdot\right>_{\mathcal{X}})$. If there exists $\alpha>0$ such that
    \begin{equation}\label{eq:infsup-condition}
        \sup_{y\in\mathcal{X}\backslash\{0\}} \frac{\mathcal{A}(x,y)}{\Vert y \Vert_{\mathcal{X}}} \geq \alpha \Vert x\Vert_{\mathcal{X}}\;,\quad \forall x\in\mathcal{X}\;,
    \end{equation}
    and
    \begin{equation}\label{eq:non-degenerate-condition}
        \sup_{x\in\mathcal{X},y\neq 0} \mathcal{A}(x,y) > 0\;,\quad \forall y\in\mathcal{X}\;,
    \end{equation}
    then there exists a unique solution $x\in\mathcal{X}$ to the problem
    \[ \mathcal{A}(x,y) = \mathcal{F}(y)\;,\quad\forall y\in\mathcal{X}\;.\]
    Furthermore, there exists ${C}>0$ (independent of $x$), such that
    \[\Vert x\Vert_{\mathcal{X}}\leq \frac{1}{{C}}\Vert \mathcal{F}\Vert_{\mathcal{X}'}\;.\]
\end{theorem}
\begin{assumpBox}\label{ass:CFL-conditions}
Assume that,
\begin{subequations}
\begin{align}
    \label{eq:ass:beta-nu-sigma}\frac{\max\left(\Vert \vec{\beta}^h + \vec{\beta}'\Vert^2_{\infty},\Vert \vec{\beta}^h \Vert^2_{\infty}\right)}{\nu \sigma }&\leq \frac{1}{6}\;,\\
    \label{eq:ass:tm-simga}\frac{\tm}{\sigma} &\leq 1\;,\\
    \label{eq:ass:beta-nu-h}\frac{\Vert \vec{\beta}^h\Vert_{\infty}h}{\nu} &\leq \frac{1}{24}\;,\\
    \label{eq:ass:beta-h}\frac{|\vec{\beta}^h(x)|}{h^e} &\leq \tm(x)\;,\quad\forall x\in\Omega^e\;, \forall e\in\mathcal{M}\;,
\end{align}
\end{subequations}
where $h^e$ is the mesh size of $\Omega^e$, which we assume to be smaller than the global mesh size $h := \max_{e\in\mathcal{M}} h^e$.
\end{assumpBox}
\begin{remark}
    Note that \eqref{eq:ass:beta-h} is a slightly weaker assumption compared to the literature \cite{francaStabilizedFiniteElement1992,evans_variational_2020} since we only need this for our theoretical analysis.
    For example, in \cite[Assumption 5]{evans_variational_2020}, the stabilization parameter was assumed to have the following form,
    \begin{equation}
        \tau_{\mathrm{M}}^{-1} = \frac{2\vert \beta \vert}{h^e} \frac{1}{f(\gamma^e)}\;, \quad f(\gamma^e) \leq \min\left(1,\frac{4 \gamma^e}{{C}_{\mathrm{inv}}} \right)\;,
    \end{equation}
    where $\gamma^e := \frac{\vert \beta^h\vert h^e}{2 \nu}$ is the element Peclet number and $f$ is some monotone function.
    Also note that standard choices of the parameter $\tau_{\mathrm{M}}$ satisfy this requirement, see Section \ref{sec:numerics} for instance, where the choice from \cite{evans_variational_2020} has been used.
\end{remark}

Relying on these assumptions, the conditions outlined in Theorem \ref{thm:stability} can be shown for the problem \eqref{eq:compact-reduced-oseen-stabilization} as stated in the following lemma.
\begin{lemma}\label{lem:stability}
    Under Assumptions \ref{ass:CFL-conditions}, there exist constants $C_{\mathrm{cont}},C_{\mathrm{inf-sup}}>0$ independent of $\nu$ and $h^e$, such that $\mathcal{A}_{\mathrm{red}}$ is a bilinear form that satsisfies,
    \begin{subequations}
        \begin{align}
        \label{eq:bounded-form} \mathcal{A}_{\mathrm{red}}(x,y) &\leq C_{\mathrm{cont}}\Vert x\Vert_{\mathbb{X}_h} \Vert y\Vert_{\mathbb{X}_h}\;,&\quad \forall x,y &\in\mathbb{X}\;,\\
\label{eq:stability-infsup}        \sup_{\substack{y^h\in\mathbb{X}_h\\y^h\neq 0}} \frac{\mathcal{A}_{\mathrm{red}}(x^h,y^h)}{\Vert y^h\Vert_{\mathbb{X}_h}} &\geq C_{\mathrm{inf-sup}} \Vert x^h\Vert_{\mathbb{X}_h}\;,&\quad \forall x^h&\in\mathbb{X}_h\;,\\
        \label{eq:stability-nondegenerate} \sup_{\substack{x^h\in\mathbb{X}_h\\x^h\neq 0}} \mathcal{A}_{\mathrm{red}}(x^h,y^h) &> 0\;,&\quad \forall y^h&\in\mathbb{X}_h\;.
    \end{align}
    \end{subequations}
\end{lemma}
\begin{remark}\label{rem:general-reduced-bounded}
    The bilinear form $\mathcal{A}_{\mathrm{red}}$ is also continuous on $\widetilde{\mathbb{X}} \times \widetilde{\mathbb{X}}$ with the same constant of continuity $C_{\mathrm{cont}}$, where $\widetilde{\mathbb{X}}$ is the extended space defined as:
    \[ \widetilde{\mathbb{X}} := H_{\vec{0}}(\curl{};\Omega,\partial\Omega)\times H_0(\div{};\Omega,\partial\Omega) \times \bigcup_{e\in\mathcal{M}} H_{\vec{0}}(\curl{};\Omega^e,\partial \Omega^e) \times\bigcup_{e\in\mathcal{M}} H_0(\div{};\Omega^e,\partial \Omega^e)\;.\]
    This can be shown by considering the same estimates from \ref{app:proof-stability-td}.
\end{remark}
The proof of Lemma \ref{lem:stability} is given in \ref{app:proof-stability-td}.
As a result, from Theorem \ref{thm:stability},  Formulation \ref{method:reduced-oseen-stabilization} is stable and unique, as shown in the following corollary.
\begin{corollary}\label{cor:stability-reduced-method}
    Under Assumption \ref{ass:CFL-conditions}, Formulation \ref{method:reduced-oseen-stabilization} has a unique solution $x^h\in\mathbb{X}_h$, for which
    \begin{equation}
        \Vert x^h\Vert_{\mathbb{X}_h} \leq C_{\mathrm{red}} \Vert \vec{f}\Vert\;,
    \end{equation}
    with ${C_\mathrm{red}}$ independent on $\nu$.
\end{corollary}
\begin{proof}
As a direct result of Theorem \ref{thm:stability} and the Lemma \ref{lem:stability}, we find
    \begin{equation*}
        \Vert x^h\Vert_{\mathbb{X}_h} \leq C \Vert \mathcal{F}_{\mathrm{red}} \Vert_{\mathbb{X}_h'}\;,
    \end{equation*}
with $C$ a constant independent of $\nu$. For $\vec{f}\in L^2(\Omega)$ and any $y^h = (\vec{\tau}^h,\vec{v}^h,\vec{v}',\vec{\tau}')\in\mathbb{X}_h$,
\begin{align*}
    (\vec{f},\vec{v}^h + \vec{v}') \leq \Vert\vec{f} \Vert \Vert \vec{v}^h + \vec{v}' \Vert\leq \Vert \vec{f} \Vert \Vert y^h \Vert_{\mathbb{X}_h}\;.
    \end{align*}
This implies that the dual norm of $\mathcal{F}_{\mathrm{red}}$ is estimated as, 
\begin{equation*}
    \Vert \mathcal{F}_{\mathrm{red}}\Vert_{\mathbb{X}_h'}:= \sup_{\substack{y^h\in\mathbb{X}_h}}\frac{\left<\vec{f},y^h\right>}{\Vert y^h\Vert_{\mathbb{X}_h}} \leq \Vert \vec{f}\Vert\;,
    \end{equation*}
showing the desired estimate.
\end{proof}
We extend this result on the stability of the reduced formulation to Formulation \ref{method:oseen-stabilization}.
\begin{theorem}
    Under Assumption \ref{ass:CFL-conditions}, Formulation \ref{method:oseen-stabilization} has a unique and stable solution $(\vec{\omega}^h,\vec{u}^h,p^h,\vec{\omega}',\vec{u}',p')\in\mathbb{S}_h$, for which
    \begin{equation}
        \Vert (\vec{\omega}^h,\vec{u}^h,\vec{u}',\vec{\omega}') \Vert^2_{\mathbb{X}_h} + \Vert p^h\Vert^2 + \Vert p'\Vert^2   \leq C \Vert \vec{f}\Vert\;,
    \end{equation}
    with $C$ independent on $\nu$. 
\end{theorem}
\begin{proof}
    The reduced solution $(\vec{\omega}^h,\vec{u}^h,\vec{\omega}',\vec{u}') \in\mathbb{X}_h$ is first computed, which is unique and stable by Corollary \ref{cor:stability-reduced-method}.
    By Lemma \ref{lem:equivalence-full-reduced-form}, the reduced solution is expanded to a stable and unique solution for Formulation \ref{method:oseen-stabilization}.
    For the coarse-scale pressure, we estimate
    \begin{align}
        \nonumber \Vert p^h\Vert &\leq\frac{1}{\beta_1} \sup_{\substack{\vec{v}^h\in\mathbb{V}_{h,0}^2}}\frac{|(p^h,\div{\vec{v}^h})|}{\Vert \vec{v}^h\Vert_{H(\div{})}}\\ 
        \nonumber &\leq \frac{1}{\beta_1}\sup_{\substack{\vec{v}^h\in\mathbb{V}_{h,0}^2}}\frac{|(f-\sigma (\vec{u}^h+\vec{u}')-\frac{1}{\sqrt{\nu}}\vec{\omega}^h\times\left[\vec{\beta}^h+\vec{\beta}'\right]-\frac{1}{\sqrt{\nu}}\vec{\omega}'\times\vec{\beta}^h-\sqrt{\nu}\curl{\vec{\omega}^h},{\vec{v}^h})|}{\Vert \vec{v}^h\Vert_{H(\div{})}}\\
        \nonumber &\leq \frac{1}{\beta_1}\left( \Vert \vec{f}\Vert + C_{\mathrm{cont}} \Vert (\vec{\omega}^h,\vec{u}^h,\vec{\omega}',\vec{u}')\Vert_{\mathbb{X}_h}\right) \leq \frac{1 + C_{\mathrm{cont}} C_{\mathrm{red}}}{\beta_1}\Vert \vec{f}\Vert\;,
    \end{align}
    where we used \eqref{eq:infsupVelPresCoarse}, Lemma \ref{lem:stability}, Remark \ref{rem:general-reduced-bounded} and Corollary \ref{cor:stability-reduced-method}. The bound for $p'$ can be shown via a similar argument.
\end{proof}
\subsection{Convergence analysis}
Lastly, we will show that Formulation \ref{method:oseen-stabilization} converges at optimal rates. This follows from the consistency of the formulation, which we demonstrate in the following proposition.
\begin{proposition}\label{lem:consistency}
    Let $(\vec{\omega},\vec{u},p)\in H_{\vec{0}}(\curl{};\Omega,\partial\Omega)\times H_0(\div{};\Omega,\partial \Omega)\times L^2_0(\Omega)$ be the solution to \eqref{eq:weak-osseen}. Then, Formulation \ref{method:reduced-oseen-stabilization} is consistent, i.e.,
    \begin{align}
        \mathcal{A}_{\mathrm{red}}\left((\vec{\omega},\vec{u},0,0),y^h\right) &= \mathcal{F}_{\mathrm{red}}(y^h)\;,\quad\forall y^h\in\mathbb{X}_h\;.
    \end{align}
\end{proposition}
\begin{proof}
    Note that the solution $(\vec{\omega},\vec{u},p)$ satisfies
    \begin{equation}\label{eq:momentum-equation-strong}
        (\sigma\,\vec{u},\,\vec{v}) + \Big(\frac{1}{\sqrt{\nu}}\,\vec{\omega}\times\vec{\beta},\,\vec{v}\Big)
    + \sqrt{\nu}\,(\curl{\vec{\omega}},\,\vec{v}) - (\vec{f},\,\vec{v}) = 0\;,
    \end{equation}
    for all $\vec{v}\in H_0^{\mathrm{ker}}(\div{};\Omega,\partial \Omega)$ as $(p,\div{\vec{v}}) = 0$.
    In particular, this holds for any $\vec{v}\in \mathbb{W}_0^{2,\mathrm{ker}}$ and $\vec{v} \in \mathbb{V}_{h,0}^{2,\mathrm{ker}}$.
    As a result, by uniqueness, $\vec{\omega}'=\vec{0}$ and $\vec{u}'=\vec{0}$ solve the fine-scale problem.
    Then, in the absence of $\vec{\omega}'$ and $\vec{u}'$, the coarse-scale subproblem of Formulation \ref{method:reduced-oseen-stabilization} reduces to
    \begin{subequations}
        \begin{align}
            \sigma(\vec{u}^h,\vec{v}^h) + 
            \frac{1}{\sqrt{\nu}}\left(\vec{\omega}^h\times\vec{\beta},\vec{v}^h\right) + \sqrt{\nu}\left(\curl{\vec{\omega}^h},\vec{v}^h\right) &= \left(\vec{f},\vec{v}^h\right)\;,\\
            \left( \vec{\omega}^h,\vec{\tau}^h\right) - \sqrt{\nu}\left( \vec{u}^h,\curl{\vec{\tau}^h}\right) &= \vec{0}\;.
        \end{align}
    \end{subequations}
    This is the Galerkin discretisation of \eqref{eq:weak-osseen} and is clearly consistent, thus showing the result.
\end{proof}

\begin{lemma}\label{lem:convergence-rate-oseen}
    Let $(\vec{\omega},\vec{u},p)\in H_0(\curl{};\Omega,\partial \Omega)\times H_0(\div{};\Omega,\partial \Omega)\times L^2_0(\Omega)$ be the solution to \eqref{eq:weak-osseen} and let 
    $\vec{x}^h$
    be the unique solution to Formulation \ref{method:reduced-oseen-stabilization} under Assumption \ref{ass:CFL-conditions}. Then, there exits a positive constant $\gamma_{\mathrm{red}}$ independent of $\nu$ and $h$ such that,
    \begin{multline}\label{eq:thereom-convergence-reduced}
    \Vert (\vec{\omega},\vec{u},0,0)-x^h\Vert_{\mathbb{X}_h} \leq  \qquad 
    \inf_{\mathclap{(\vec{\tau}^h,\vec{v}^h)\in\mathbb{V}_{h,\vec{0}}^1\times\mathbb{V}_{h,0}^{2,\mathrm{ker}}}}  \qquad \gamma_{\mathrm{red}}  \left(\Vert \vec{u} - \vec{v}^h\Vert^2 + \Vert \vec{\omega}-\vec{\tau}^h\Vert^2 + \nu\Vert \curl{\vec{\omega}-\vec{\tau}^h}\Vert^2 \right)\;.
\end{multline}
\end{lemma}
\begin{proof}
    Let $y^h\in\mathbb{X}_h$ and define
\begin{equation*}
    z^h =\mathrm{argmax}_{\hat{z}^h\in\mathbb{X}_h} \frac{\mathcal{A}_{\mathrm{red}}(x^h-y^h,\hat{z}^h)}{\Vert \hat{z}^h\Vert_{\mathbb{X}_h}}\;.
\end{equation*}
Then, by Lemma \ref{lem:stability} and Proposition \ref{lem:consistency},
\begin{align*}
    \Vert x^h-y^h\Vert_{\mathbb{X}_h} &\leq \frac{C^{-1}_{\mathrm{inf-sup}}}{\Vert z^h\Vert_{\mathbb{X}_h}}\mathcal{A}_{\mathrm{red}}(x^h-y^h,z^h)\\
     &=\frac{C^{-1}_{\mathrm{inf-sup}}}{\Vert z^h\Vert_{\mathbb{X}_h}}\mathcal{A}_{\mathrm{red}}\left((\vec{\omega},\vec{u},0,0)-y^h,z^h\right)\\
    \quad &\leq C^{-1}_{\mathrm{inf-sup}}C_{\mathrm{cont}}\Vert (\vec{\omega},\vec{u},0,0) - y^h \Vert_{\mathbb{X}_h}\;,
\end{align*}
where $C_{\mathrm{cont}}$ is the continuity constant from Lemma \ref{lem:stability}. Then
\begin{align}
    \nonumber \Vert (\vec{\omega},\vec{u},0,0)-x^h\Vert_{\mathbb{X}_h} &\leq \Vert (\vec{\omega},\vec{u},0,0)-y^h\Vert_{\mathbb{X}_h} + \Vert x^h-y^h\Vert_{\mathbb{X}_h}\\
    \nonumber &\leq (1 + C^{-1}_{\mathrm{inf-sup}}C_{\mathrm{cont}})\Vert (\vec{\omega},\vec{u},0,0)-y^h\Vert_{\mathbb{X}_h}\;.
\end{align}
For the r.h.s., it suffices to approximate the coarse-scale components, thereby proving the desired result.
\end{proof}
\begin{remark}\label{rem:reduced-formulation}
    As the space $\mathbb{V}_{h,0}^{2,\mathrm{ker}} \subset H_0^{\mathrm{ker}}(\div{};\Omega,\partial \Omega)$, the infimum in \eqref{eq:thereom-convergence-reduced} over $(\vec{\omega^h},\vec{u}^h)\in\mathbb{V}_{h,\vec{0}}^1\times \mathbb{V}_{h,0}^{2,\mathrm{ker}}$ can be replaced by $(\vec{\omega^h},\vec{u}^h)\in\mathbb{V}_{h,\vec{0}}^1\times \mathbb{V}_{h,0}^{2}$, see \cite[pages 58,59]{gatica_simple_2014}. Crucially, the resulting constant $\gamma_{\mathrm{red}}$ stays independent of $\nu$ and mesh size $h$.
    As a result, the formulation is independent of the coarse-scale pressure and is thus pressure-robust.
\end{remark}
Similarly, we can estimate the error of the coarse-scale pressure with the following lemma.  
\begin{lemma}\label{lem:convergence-pressure}
Let $(\vec{\omega},\vec{u},p)\in H_0(\curl{};\Omega,\partial \Omega)\times H_0(\div{};\Omega,\partial \Omega)\times L^2_0(\Omega)$ be the solution to \eqref{eq:weak-osseen} and let $(\vec{\omega}^h,\vec{u}^h,p^h,\vec{\omega}',\vec{u}',p')$ be the unique solution to Formulation \ref{method:oseen-stabilization} under Assumption \ref{ass:CFL-conditions}. Then, there exits a positive constant $\gamma$ independent of $\nu$ and $h$ such that,
    \begin{multline}
    \Vert p - p^h \Vert^2 
    \leq
    \quad \inf_{\mathclap{(\vec{\tau}^h,\vec{v}^h,q^h)\in\mathbb{V}_{h,\vec{0}}^1\times\mathbb{V}_{h,0}^{2}\times \mathbb{V}_{h,0}^3}} \quad \gamma_{p} \left(\Vert \vec{u} - \vec{v}^h\Vert^2 + \Vert \vec{\omega}-\vec{\tau}^h\Vert^2 + \nu\Vert \curl{\vec{\omega}-\vec{\tau}^h}\Vert^2 + \Vert p - q^h\Vert^2 \right)\;.
\end{multline}
\end{lemma}

\begin{proof}
From \eqref{eq:weak-osseen-momentum} and \eqref{eq:oseen-stabilized-formulation-coarse-momentum}, the following equality holds
\begin{multline}
 \nonumber \sigma(\vec{u},\vec{v}^h)+ \frac{1}{\sqrt{\nu}}(\vec{\omega}\times\left[\vec{\beta}^h+\vec{\beta}'\right],\vec{v}^h)+ \sqrt{\nu}(\curl{\vec{\omega}},\vec{v}^h) - (p,\div{\vec{v}^h}) = \\
  \sigma(\vec{u}^h+\vec{u}',\vec{v}^h)+ \frac{1}{\sqrt{\nu}}(\vec{\omega}^h\times\left[\vec{\beta}^h+\vec{\beta}'\right],\vec{v}^h)+ \frac{1}{\sqrt{\nu}}(\vec{\omega}'\times\vec{\beta}^h,\vec{v}^h)+ \sqrt{\nu}(\curl{\vec{\omega}^h},\vec{v}^h) - (p^h,\div{\vec{v}^h})\;.
\end{multline}
As a result, we note that for any $\vec{v}^h\in\mathbb{V}_{h,0}^2$,
\begin{align}
    \nonumber &\frac{(p-p^h,\div{\vec{v}^h})}{\Vert \vec{v}^h\Vert_{H(\div{})}} \\
    \nonumber =&\frac{\left|\left( \sigma(\vec{u}-\vec{u}^h)+  \frac{1}{\sqrt{\nu}}\left(\vec{\omega}-\vec{\omega^h}\right)\times\left[\vec{\beta}^h+\vec{\beta}'\right] + \frac{1}{\sqrt{\nu}}\left(\vec{0}-\vec{\omega}'\right)\times\vec{\beta}^h + \sqrt{\nu}\curl{\vec{\omega}-\vec{\omega}^h},\vec{v}^h \right)\right|}{\Vert \vec{v}^h \Vert_{H(\div{};\Omega)}}\\
    \nonumber\leq& C_{\mathrm{cont}}\Vert (\vec{\omega},\vec{u},0,0) - (\vec{\omega}^h,\vec{u}^h,\vec{\omega}',\vec{v}') \Vert_{\mathbb{X}_h}\;,
\end{align}
where we invoked Remark \ref{rem:general-reduced-bounded} in the last inequality.
Then, for any $q^h\in\mathbb{V}_{h,0}^3$ and using \eqref{eq:infsupVelPresCoarse}, 
\begin{align}
    \beta_1 \Vert q^h - p^h \Vert &\leq \sup_{\vec{v}^h\in \mathbb{V}_{h,0}^2} \frac{(q^h-p^h,\div{\vec{v}^h})}{\Vert \vec{v}^h \Vert_{H(\div{})}} \\
    &\leq \sup_{\vec{v}^h\in \mathbb{V}_{h,0}^2} \frac{(q^h-p,\div{\vec{v}^h})}{\Vert \vec{v}^h \Vert_{H(\div{})}} + \sup_{\vec{v}^h\in \mathbb{V}_{h,0}^2} \frac{(p-p^h,\div{\vec{v}^h})}{\Vert \vec{v}^h \Vert_{H(\div{})}} \\
    & \leq \Vert q^h - p \Vert + \sup_{\vec{v}^h\in \mathbb{V}_{h,0}^2} \frac{(p-p^h,\div{\vec{v}^h})}{\Vert \vec{v}^h \Vert_{H(\div{})}} 
\end{align}
so that
\begin{align}
    & \Vert p - p^h \Vert \leq \Vert p - q^h \Vert + \Vert q^h - p^h \Vert\\
    \leq & \left(1+\frac{1}{\beta_1}\right)\Vert p - q^h \Vert + \frac{C_{cont}}{\beta_1}\Vert (\vec{\omega},\vec{u},0,0) - (\vec{\omega}^h,\vec{u}^h,\vec{\omega}',\vec{v}') \Vert_{\mathbb{X}_h}
\end{align}
and the desired result is found upon using Lemma \ref{lem:convergence-rate-oseen} and Remark \ref{rem:reduced-formulation}.
\end{proof}
By Lemma \ref{lem:convergence-rate-oseen} and \ref{lem:convergence-pressure}, the formulation convergences optimally with respect to the spaces $\mathbb{V}_{h,\vec{0}}^1,\mathbb{V}_{h,0}^2$ and $\mathbb{V}_{h,0}^3$.
To obtain explicit convergence rates, we assume the following on the approximation power in the spaces of vorticity and velocity.
\begin{assumpBox}\label{ass:approximation-power}
    There exist integers $k^1,k^2$ and $k^3$, interpolation operators,
    \begin{equation*}
      \Pi^1_0 : H_{\vec{0}}(\curl{};\Omega,\partial\Omega) \rightarrow\mathbb{V}_{h,\vec{0}}^1\;,
      \Pi^2_0 : H_0^{\mathrm{ker}}(\div{};\Omega,\partial \Omega) \rightarrow \mathbb{V}_{h,0}^{2,\mathrm{ker}}\;,
      \Pi_0^3 : L^2_0(\Omega) \rightarrow \mathbb{V}_{h,0}^3\;,  
    \end{equation*}
    and interpolation constant $C_{\mathrm{interp}}$ independent on the global mesh size $h$, such that:
     for every $\vec{\tau} \in H_{\vec{0}}(\curl{};\Omega,\partial\Omega) $ and $0\leq \ell \leq k^1$,
    \begin{subequations}
        \begin{align}
          \vert \vec{\tau} - \Pi^1_0 \vec{\tau}\vert_{H^\ell(\curl{})} &\leq C_{\mathrm{interp}} h^{k^1 - \ell} \vert \vec{\tau} \vert_{H^{k^1}(\curl{})}\;,\;&\forall \tau \in H^{k^1}_{\vec{0}}(\curl{};\Omega,\partial\Omega)\;,\;0 \leq \ell \leq k^1\;,\\
          \vert \vec{v} - \Pi^2_0 \vec{v}\vert_{H^\ell(\div{})} &\leq C_{\mathrm{interp}} h^{k^2 - \ell} \vert \vec{v} \vert_{H^{k^2}(\div{})}\;,\;&\forall \vec{v} \in H^{k^2}_0(\div{};\Omega,\partial \Omega)\;,\;0 \leq \ell \leq k^2\;,\\
            \vert q - \Pi^3_0 q\vert_{H^\ell} &\leq C_{\mathrm{interp}} h^{k^3 - \ell} \vert q \vert_{H^{k^3}}\;,\;&\forall q \in H^{k^3}_0(\Omega)\;,\;0 \leq \ell \leq k^3\;.
        \end{align}
    \end{subequations}
\end{assumpBox}
Note that this is a non-standard assumption, as it requires the existence of commuting projectors (between the coarse-scale velocity and coarse-scale pressure spaces).
These exist for divergence-conforming B-spline discretizations \cite{buffa_isogeometric_2011} and certain divergence conforming finite-element discretizations \cite{guzmanConformingDivergencefreeStokes2013}.
With Assumption \ref{ass:approximation-power}, the total convergence rate of our formulation can be stated as below.
\begin{theorem}\label{thm:optimal-convergence-rate-oseen}
    Let $(\vec{\omega},\vec{u},p)\in H_0(\curl{};\Omega,\partial \Omega)\times H_0(\div{};\Omega,\partial \Omega)\times L^2_0(\Omega)$ be the solution to \eqref{eq:weak-osseen} and let $(\vec{\omega}^h,\vec{u}^h,p^h,\vec{\omega}',\vec{u}',p')$ be the unique solution to Formulation \ref{method:oseen-stabilization} under Assumption \ref{ass:CFL-conditions}.
    Furthermore, with Assumption \ref{ass:approximation-power} in place, let integers $k^1,k^2,k^3$ be such that $\omega \in H^{k^1}(\curl{};\Omega)$, $\vec{u} \in H^{k^2}(\div{};\Omega)$ and $p \in H^{k^3}(\Omega)$.
    Then,
    \begin{multline}\label{eq:thm-optimal-convergence-rate}
    \Vert \vec{u} - \vec{u}^h - \vec{u}'\Vert^2 + \Vert \vec{\omega}-\vec{\omega}^h\Vert^2 + \nu\Vert \curl{\vec{\omega}-\vec{\omega}^h}\Vert^2 + \Vert p - p^h\Vert^2 \leq\\
    C_{\mathrm{approx}}^2 h^{2k^2} \vert \vec{u} \vert_{H^{k^2}(\div{})}^2 + \max(1,\nu)C_{\mathrm{approx}}^2 h^{2k^1} \vert \vec{\omega} \vert_{H^{k^1}(\curl{})}^2 + C_{\mathrm{approx}}^2 h^{2k^3} \vert p \vert_{H^{k^3}} \;.
\end{multline}
where $C_{\mathrm{approx}}$ is independent on $\nu$ and mesh size $h$.
\end{theorem}
\begin{proof}
    This is a direct result of Lemma \ref{lem:convergence-rate-oseen}, Lemma \ref{lem:convergence-pressure} and Assumption \ref{ass:approximation-power} by choosing $\vec{\tau}^h = \Pi_0^1 \vec{\omega}, \vec{v}^h = \Pi_0^2 \vec{u}$ and $q^h = \Pi_0^3 p$.
\end{proof}

\begin{remark}\label{rem:adaptivity}
Note, for a constant $\vec{\beta}^h$, the fine-scale velocity can be bounded as
\[\Vert \vec{u}'\Vert^2 \leq \frac{h}{\vert \vec{\beta}^h\vert}(\tm \vec{u}',\vec{u}') \leq \frac{h}{\vert \vec{\beta}^h\vert}\Vert (\vec{\omega}^h,\vec{u}^h,\vec{\omega}',\vec{u}')\Vert_{\mathbb{X}_h}^2\]
and, as such, will vanish at a faster rate than the coarse-scale error terms.
With this observation, the fine-scale velocity $\vec{u}'$ can be considered an a posteriori error indicator.
\end{remark}

\begin{remark}
Theorem \ref{thm:optimal-convergence-rate-oseen} can be altered to control $\Vert \vec{u}-\vec{u}^h\Vert$ instead of $\Vert \vec{u}-\vec{u}^h-\vec{u}'\Vert$, by removing the term $\sigma (\vec{u}',\vec{v}^h)$ from Formulation \ref{method:oseen-stabilization}.
This change is necessary to be able to estimate the inf-sup constant following our analysis, since, without this change, $|(\vec{v}^h+\vec{v}',\vec{v}^h+ \vec{v}')| \ngeq C \left( \Vert \vec{v}^h\Vert^2 + \Vert \vec{v}'\Vert^2 \right) $.
Note that $\sigma (\vec{u}',\vec{v}^h)$ is equivalent to the $(\partial_t \vec{u}',\vec{v}^h)$ in Formulation \ref{method:sd-ns-stabilization}.
While we did not omit this term in our work and thus obtained the a priori estimates in Theorem \ref{thm:optimal-convergence-rate-oseen}, the numerical results suggest that $\Vert \vec{u}- \vec{u}^h\Vert$ converges optimally nevertheless.
\end{remark}

\section{Numerical Results}\label{sec:numerics}
Formulation \ref{method:ns-stabilization} has been implemented in Nutils, an open source Finite Element package for Python; see \cite{van_zwieten_nutils_2022}.
However, due to computational limitations related to Nutils (we will discuss this further in Section \ref{sec:choice-stabilization-parameter}), we only test two-dimensional problems.
For this, we discretize the 2D de Rham complex:
\begin{equation}
    H^1(\Omega) \xrightarrow[]{\rot{}} H(\div{};\Omega) \xrightarrow[]{\div{}} L^2(\Omega)\;,
\end{equation}
for which we use maximally smooth B-spline spaces. 
Given $\vec{k} = (k_x,k_y)$ and $\vec{k}' = (k'_x,k'_y)$ the chosen spline degree vectors for coarse and fine scales, respectively, we will denote the associated tensor-product B-spline spaces as $\mathbb{B}_{k_x,k_y}(D)$ and $\mathbb{B}_{k'_x,k'_y}(D)$, respectively, over domain $D$.
Then, we pick the finite-dimensional coarse and fine spaces for the different unknowns as:
\begin{subequations}
\begin{align}
    \mathbb{V}_{h}^0 &= \mathbb{B}_{k_x,k_y}(\Omega) \cap H^1_{{\phi}}(\Omega,\Gamma_{\hat{\vec{\omega}}})\;,\\
    \mathbb{V}_{h}^1 & = \left[ \mathbb{B}_{k_x,k_y-1}(\Omega)\times \mathbb{B}_{k_x-1,k_y}(\Omega) \right] \cap H_{\phi}(\div{},\Omega,\Gamma_{\hat{u}})\;,\\
    \mathbb{V}_{h}^2 & = \mathbb{B}_{k_x-1,k_y-1}(\Omega) \cap L^2_0(\Omega)\;,\\
    \mathbb{W}_{h,{0}}^0 &=  \cup_{e\in\mathcal{M}} \left( \mathbb{B}_{k'_x,k'_y}(\Omega^e) \cap H^1_{{0}}(\Omega^e,\partial \Omega^e)\right)\;,\\
    \mathbb{W}_{h,0}^1 & =  \cup_{e\in\mathcal{M}} \left( \left[\mathbb{B}_{k'_x,k'_y-1}(\Omega^e)\times \mathbb{B}_{k'_x-1,k'_y}(\Omega^e)\right] \cap H_{0}(\div{},\Omega^e,\partial\Omega^e) \right)\;,\\
    \mathbb{W}_{h,0}^2 & = \cup_{e\in\mathcal{M}} \left( \mathbb{B}_{k'_x-1,k'_y-1}(\Omega^e) \cap L^2_0(\Omega^e)\right)\;,
\end{align}
\end{subequations}
where $\mathcal{M}$ is a tensor-product mesh of equally sized mesh elements (chosen for simplicity).
In Section \ref{sec:choice-stabilization-parameter}, we will first discuss the choice of the stabilization parameter $\tau_{\mathrm{M}}$ and how it affects implementation complexity.
In Section \ref{sec:NR-convergence-rate}, we numerically validate the convergence rate of our formulation, as well as the choice of fine-scale degree $\vec{k}'$.
Lastly, we qualitatively and quantitatively compare the stabilized and unstabilized solutions.
To do this, Section \ref{sec:NR-subgrid-analysis} considers the unsteady shear layer roll-up problem in both viscous and inviscid cases.
Following that, in Section \ref{sec:NR-comparison}, we consider the steady-state lid-driven cavity problem.

\subsection{Choices of stabilization parameter}
\label{sec:choice-stabilization-parameter}
For Formulation \ref{method:ns-stabilization}, the equations are linearized via Picard iterations. For this choice, we choose the stabilization parameter based on \cite{codina_stabilized_2002},
\begin{equation}\label{eq:tm-central}
    \tm := \sqrt{C_1\frac{\vec{u}^{h,\PicIt,n+\frac{1}{2}}\cdot\vec{u}^{h,\PicIt,n+\frac{1}{2}}}{h^2} + C_2\frac{\mathrm{Re}^{-2}}{4 h^4}}\;,
\end{equation}
where $C_1>0$ and $C_2>0$ are parameters to be chosen independent of mesh size $h$. 
While these can be fitted to data, we instead base our choice on \cite{evans_variational_2020}, and we pick:
\begin{equation}
    C_1 := \max\left(k_x,k_y\right)^2\;,C_2 := \max\left(k_x',k_y'\right)^4\;.
\end{equation}
With this choice, we note that $\tau_{\mathrm{M}}$ does not depend on the next Picard iterate $\PicIt+1$, leading to a linear fine-scale model that can be eliminated from the formulation via Schur's complement.
As a result, Formulation \ref{method:ns-stabilization} has the same number of unknown degrees of freedom as the unstabilized formulation.
Moreover, since the fine-scale problems are defined elementwise, the additional cost of associated assembly and solution scales linearly with the number of elements, and the fine-scale problem can be solved in an embarrassingly parallel fashion.
Hence, the cost of solving the linear system associated with the coarse-scale problem will dominate in general.

In certain situations, linearizations via Newton's method are preferred over Picard iterates, as they converge faster.
Unfortunately, with consistent Newton linearization, the fine-scale problem \eqref{eq:sd-formulation-fine-scale-momentum}--\eqref{eq:sd-formulation-fine-scale-vorticity} cannot be eliminated due to the presence of the term $(\vec{\omega}^h\times\vec{u}',\vec{v}')$, and when using the stabilization parameter \eqref{eq:tm-central}, which depends on the next time step $n+1$.
For such cases, we propose an alternative stabilization parameter, 
\begin{equation}\label{eq:tm-upwind}
    \tm := \sqrt{C_1\frac{\vec{u}^{h,n}\cdot\vec{u}^{h,n}}{h^2} + C_2\frac{\mathrm{Re}^{-2}}{4 h^4}}\;,
\end{equation}
and an alteration to Formulation \ref{method:sd-ns-stabilization}, where the non-linear convection term in the fine-scale equation is temporally discretized as
\begin{equation}\label{eq:non-linear-fine-scale-convection-upwind}
(\vec{\omega}^{h,n}\times \vec{u}'^{n+\nicefrac{1}{2}},\vec{v}')\;.
\end{equation}
As a result, the fine-scale problem becomes linear again and can be eliminated even when using Newton's method.
As \eqref{eq:tm-upwind} and \eqref{eq:non-linear-fine-scale-convection-upwind} only depend on the previous solutions $\vec{u}^{h,n},\vec{\omega}^{h,n}$, we will refer to this stabilization as {semi-CN} stabilization, and the use of \eqref{eq:tm-central} as full-CN stabilization.

\subsection{Optimal convergence rates}
\label{sec:NR-convergence-rate}
To numerically demonstrate that the optimal convergence rates, as suggested by Theorem \ref{thm:optimal-convergence-rate-oseen} for the Oseen problem, are also achieved for the incompressible Navier--Stokes equations, we perform a numerical study with the steady-state regularized lid-driven cavity flow. 
We also perform a study for the Taylor--Green vortex flow from \cite{evans_divergence-free_2011}.
\subsubsection{Steady regularized lid-driven cavity flow}
We start by investigating the steady-state version of the incompressible Navier--Stokes equations. 
Here, Formulation \ref{method:ns-stabilization-gen-bc} is altered by removing the unsteady term
\[  \left(\partial_t \vec{u}^{h} + \partial_t \vec{u}',\vec{v}^h+\vec{v}'\right)\;. \]
We remark that our analysis of the Oseen equation relies on the time parameter $\sigma>0$, which plays the role of the reciprocal of the time-step size $\Delta t$.
This analysis thus does not extend to the steady-state case; however, a similar analysis can be performed for the case $\sigma = 0$ by modifying Assumptions \ref{ass:CFL-conditions} to place further restrictions on the elementwise Reynolds number.

\begin{figure}[htp]
    \begin{subfigure}{0.32\textwidth}
        \centering
        \caption{pressure}
        \newcommand{\logLogSlopeTriangle}[5]
{

    \pgfplotsextra
    {
        \pgfkeysgetvalue{/pgfplots/xmin}{\xmin}
        \pgfkeysgetvalue{/pgfplots/xmax}{\xmax}
        \pgfkeysgetvalue{/pgfplots/ymin}{\ymin}
        \pgfkeysgetvalue{/pgfplots/ymax}{\ymax}

        \pgfmathsetmacro{\xArel}{#1}
        \pgfmathsetmacro{\yArel}{#3}
        \pgfmathsetmacro{\xBrel}{#1-#2}
        \pgfmathsetmacro{\yBrel}{\yArel}
        \pgfmathsetmacro{\xCrel}{\xArel}

        \pgfmathsetmacro{\lnxB}{\xmin*(1-(#1-#2))+\xmax*(#1-#2)} 
        \pgfmathsetmacro{\lnxA}{\xmin*(1-#1)+\xmax*#1} 
        \pgfmathsetmacro{\lnyA}{\ymin*(1-#3)+\ymax*#3} 
        \pgfmathsetmacro{\lnyC}{\lnyA+#4*(\lnxA-\lnxB)}
        \pgfmathsetmacro{\yCrel}{\lnyC-\ymin)/(\ymax-\ymin)} 

        \coordinate (A) at (rel axis cs:\xArel,\yArel);
        \coordinate (B) at (rel axis cs:\xBrel,\yBrel);
        \coordinate (C) at (rel axis cs:\xCrel,\yCrel);

        \draw[#5]   (A)-- node[pos=0.5,anchor=south] {1}
                    (B)-- 
                    (C)-- node[pos=0.5,anchor=west] {#4}
                    cycle;
    }
}
\begin{tikzpicture}
\begin{loglogaxis}[%
  xmin=1, xmax=80,
  ymin=10^(-7), ymax=5,
  width=\textwidth,
  height=7cm,
  ticklabel style={font=\tiny}]
    \addplot[GREEN, mark = x] table [GREEN, x=nelems, y=VMS k 1 k’ 2, col sep=comma, mark = x] {Images/convergence_results_p.csv};
    \addplot[GREEN, mark = +] table [GREEN, x=nelems, y=VMS k 1 k’ 3, col sep=comma, mark = +] {Images/convergence_results_p.csv};
    \addplot[GREEN, mark = o] table [GREEN, x=nelems, y=VMS k 1 k’ 4, col sep=comma, mark = o] {Images/convergence_results_p.csv};
    \addplot[GREEN ,style=dashed] table [GREEN, x=nelems, y=No VMS k 1 k’ 1, col sep=comma, style = dashed] {Images/convergence_results_p.csv};

    \addplot[ORANGE] table [ORANGE, x=nelems, y=VMS k 2 k’ 2, col sep=comma] {Images/convergence_results_p.csv};
    \addplot[ORANGE, mark = x] table [ORANGE, x=nelems, y=VMS k 2 k’ 3, col sep=comma] {Images/convergence_results_p.csv};
    \addplot[ORANGE, mark = +] table [ORANGE, x=nelems, y=VMS k 2 k’ 4, col sep=comma] {Images/convergence_results_p.csv};
    \addplot[ORANGE, mark = o] table [ORANGE, x=nelems, y=VMS k 2 k’ 5, col sep=comma] {Images/convergence_results_p.csv};
    \addplot[ORANGE ,style=dashed] table [ORANGE, x=nelems, y=No VMS k 2 k’ 2, col sep=comma] {Images/convergence_results_p.csv};

    \addplot[MAGENTA] table [MAGENTA, x=nelems, y=VMS k 3 k’ 3, col sep=comma] {Images/convergence_results_p.csv};
    \addplot[MAGENTA, mark = x] table [MAGENTA, x=nelems, y=VMS k 3 k’ 4, col sep=comma] {Images/convergence_results_p.csv};
    \addplot[MAGENTA, mark = +] table [MAGENTA, x=nelems, y=VMS k 3 k’ 5, col sep=comma] {Images/convergence_results_p.csv};
    \addplot[MAGENTA, mark = o] table [MAGENTA, x=nelems, y=VMS k 3 k’ 6, col sep=comma] {Images/convergence_results_p.csv};
    \addplot[MAGENTA,style=dashed] table [MAGENTA, x=nelems, y=No VMS k 3 k’ 3, col sep=comma] {Images/convergence_results_p.csv};

    \logLogSlopeTriangle{0.8}{0.15}{0.78}{-1}{GREEN};
    \logLogSlopeTriangle{0.8}{0.15}{0.55}{-2}{ORANGE};
    \logLogSlopeTriangle{0.8}{0.15}{0.37}{-3}{MAGENTA};
    
\end{loglogaxis}
\end{tikzpicture}%
    \end{subfigure}
    \hfill
    \begin{subfigure}{0.32\textwidth}
        \centering
        \caption{velocity}
        \newcommand{\logLogSlopeTriangle}[5]
{

    \pgfplotsextra
    {
        \pgfkeysgetvalue{/pgfplots/xmin}{\xmin}
        \pgfkeysgetvalue{/pgfplots/xmax}{\xmax}
        \pgfkeysgetvalue{/pgfplots/ymin}{\ymin}
        \pgfkeysgetvalue{/pgfplots/ymax}{\ymax}

        \pgfmathsetmacro{\xArel}{#1}
        \pgfmathsetmacro{\yArel}{#3}
        \pgfmathsetmacro{\xBrel}{#1-#2}
        \pgfmathsetmacro{\yBrel}{\yArel}
        \pgfmathsetmacro{\xCrel}{\xArel}

        \pgfmathsetmacro{\lnxB}{\xmin*(1-(#1-#2))+\xmax*(#1-#2)} 
        \pgfmathsetmacro{\lnxA}{\xmin*(1-#1)+\xmax*#1} 
        \pgfmathsetmacro{\lnyA}{\ymin*(1-#3)+\ymax*#3} 
        \pgfmathsetmacro{\lnyC}{\lnyA+#4*(\lnxA-\lnxB)}
        \pgfmathsetmacro{\yCrel}{\lnyC-\ymin)/(\ymax-\ymin)} 

        \coordinate (A) at (rel axis cs:\xArel,\yArel);
        \coordinate (B) at (rel axis cs:\xBrel,\yBrel);
        \coordinate (C) at (rel axis cs:\xCrel,\yCrel);

        \draw[#5]   (A)-- node[pos=0.5,anchor=south] {1}
                    (B)-- 
                    (C)-- node[pos=0.5,anchor=west] {#4}
                    cycle;
    }
}

\begin{tikzpicture}
\begin{loglogaxis}[%
  xmin=1, xmax=80,
  ymin=10^(-8), ymax=0.07,
  width=\textwidth,
  height=7cm,
  ticklabel style={font=\tiny}]
    \addplot[GREEN, mark = x] table [x=nelems, y=VMS k 1 k’ 2, col sep=comma] {Images/convergence_results_u.csv};
    \addplot[GREEN, mark = +] table [x=nelems, y=VMS k 1 k’ 3, col sep=comma] {Images/convergence_results_u.csv};
    \addplot[GREEN, mark = o] table [x=nelems, y=VMS k 1 k’ 4, col sep=comma] {Images/convergence_results_u.csv};
    \addplot[GREEN, style = dashed] table [x=nelems, y=No VMS k 1 k’ 1, col sep=comma] {Images/convergence_results_u.csv};

    \addplot[ORANGE] table [x=nelems, y=VMS k 2 k’ 2, col sep=comma] {Images/convergence_results_u.csv};
    \addplot[ORANGE, mark = x] table [x=nelems, y=VMS k 2 k’ 3, col sep=comma] {Images/convergence_results_u.csv};
    \addplot[ORANGE, mark = +] table [x=nelems, y=VMS k 2 k’ 4, col sep=comma] {Images/convergence_results_u.csv};
    \addplot[ORANGE, mark = o] table [x=nelems, y=VMS k 2 k’ 5, col sep=comma] {Images/convergence_results_u.csv};
    \addplot[ORANGE ,style=dashed] table [x=nelems, y=No VMS k 2 k’ 2, col sep=comma] {Images/convergence_results_u.csv};

    \addplot[MAGENTA] table [x=nelems, y=VMS k 3 k’ 3, col sep=comma] {Images/convergence_results_u.csv};
    \addplot[MAGENTA, mark = x] table [x=nelems, y=VMS k 3 k’ 4, col sep=comma] {Images/convergence_results_u.csv};
    \addplot[MAGENTA, mark = +] table [x=nelems, y=VMS k 3 k’ 5, col sep=comma] {Images/convergence_results_u.csv};
    \addplot[MAGENTA, mark = o] table [x=nelems, y=VMS k 3 k’ 6, col sep=comma] {Images/convergence_results_u.csv};
    \addplot[MAGENTA,style=dashed] table [x=nelems, y=No VMS k 3 k’ 3, col sep=comma] {Images/convergence_results_u.csv};

    \logLogSlopeTriangle{0.8}{0.15}{0.85}{-1}{GREEN};
    \logLogSlopeTriangle{0.8}{0.15}{0.67}{-2}{ORANGE};
    \logLogSlopeTriangle{0.8}{0.15}{0.45}{-3}{MAGENTA};
\end{loglogaxis}
\end{tikzpicture}%
    \end{subfigure}
    \hfill
    \begin{subfigure}{0.32\textwidth}
        \centering
        \caption{vorticity}

\newcommand{\logLogSlopeTriangle}[5]
{

    \pgfplotsextra
    {
        \pgfkeysgetvalue{/pgfplots/xmin}{\xmin}
        \pgfkeysgetvalue{/pgfplots/xmax}{\xmax}
        \pgfkeysgetvalue{/pgfplots/ymin}{\ymin}
        \pgfkeysgetvalue{/pgfplots/ymax}{\ymax}

        \pgfmathsetmacro{\xArel}{#1}
        \pgfmathsetmacro{\yArel}{#3}
        \pgfmathsetmacro{\xBrel}{#1-#2}
        \pgfmathsetmacro{\yBrel}{\yArel}
        \pgfmathsetmacro{\xCrel}{\xArel}

        \pgfmathsetmacro{\lnxB}{\xmin*(1-(#1-#2))+\xmax*(#1-#2)} 
        \pgfmathsetmacro{\lnxA}{\xmin*(1-#1)+\xmax*#1} 
        \pgfmathsetmacro{\lnyA}{\ymin*(1-#3)+\ymax*#3} 
        \pgfmathsetmacro{\lnyC}{\lnyA+#4*(\lnxA-\lnxB)}
        \pgfmathsetmacro{\yCrel}{\lnyC-\ymin)/(\ymax-\ymin)} 

        \coordinate (A) at (rel axis cs:\xArel,\yArel);
        \coordinate (B) at (rel axis cs:\xBrel,\yBrel);
        \coordinate (C) at (rel axis cs:\xCrel,\yCrel);

        \draw[#5]   (A)-- node[pos=0.5,anchor=south] {1}
                    (B)-- 
                    (C)-- node[pos=0.5,anchor=west] {#4}
                    cycle;
    }
}

\begin{tikzpicture}[spy using outlines=
{circle, magnification=8, connect spies}]
\begin{loglogaxis}[%
  xmin=1, xmax=80,
  ymin=10^(-9), ymax=0.5,
  width=\textwidth,
  height=7cm,
  ticklabel style={font=\tiny}]
    \addplot[GREEN, mark = x] table [x=nelems, y=VMS k 1 k’ 2, col sep=comma] {Images/convergence_results_w.csv};
    \addplot[GREEN, mark = +] table [x=nelems, y=VMS k 1 k’ 3, col sep=comma] {Images/convergence_results_w.csv};
    \addplot[GREEN, mark = o] table [x=nelems, y=VMS k 1 k’ 4, col sep=comma] {Images/convergence_results_w.csv};
    \addplot[GREEN, style = dashed] table [x=nelems, y=No VMS k 1 k’ 1, col sep=comma] {Images/convergence_results_w.csv};

    \addplot[ORANGE] table [x=nelems, y=VMS k 2 k’ 2, col sep=comma] {Images/convergence_results_w.csv};
    \addplot[ORANGE, mark = x] table [x=nelems, y=VMS k 2 k’ 3, col sep=comma] {Images/convergence_results_w.csv};
    \addplot[ORANGE, mark = +] table [x=nelems, y=VMS k 2 k’ 4, col sep=comma] {Images/convergence_results_w.csv};
    \addplot[ORANGE, mark = o] table [x=nelems, y=VMS k 2 k’ 5, col sep=comma] {Images/convergence_results_w.csv};
    \addplot[ORANGE ,style= dashed] table [x=nelems, y=No VMS k 2 k’ 2, col sep=comma] {Images/convergence_results_w.csv};

    \addplot[MAGENTA] table [x=nelems, y=VMS k 3 k’ 3, col sep=comma] {Images/convergence_results_w.csv};
    \addplot[MAGENTA, mark = x] table [x=nelems, y=VMS k 3 k’ 4, col sep=comma] {Images/convergence_results_w.csv};
    \addplot[MAGENTA, mark = +] table [x=nelems, y=VMS k 3 k’ 5, col sep=comma] {Images/convergence_results_w.csv};
    \addplot[MAGENTA, mark = o] table [x=nelems, y=VMS k 3 k’ 6, col sep=comma] {Images/convergence_results_w.csv};
    \addplot[MAGENTA, style = dashed] table [x=nelems, y=No VMS k 3 k’ 3, col sep=comma] {Images/convergence_results_w.csv};
    \coordinate (spypoint)     at (0.7,-3);
    \coordinate (magnifyglass) at (1,-16); 

    \logLogSlopeTriangle{0.8}{0.15}{0.77}{-2}{GREEN};
    \logLogSlopeTriangle{0.8}{0.15}{0.6}{-3}{ORANGE};
    \logLogSlopeTriangle{0.8}{0.15}{0.4}{-4}{MAGENTA};
\end{loglogaxis}
\spy [size=2.5cm] on (spypoint)
    in node[fill=white] at (magnifyglass); 

\end{tikzpicture}%
    \end{subfigure}\\
    \begin{subfigure}{\textwidth}
        \hspace{0.25\textwidth}
        \begin{tikzpicture}[scale=0.5]
        \pgflowlevelsynccm
        \draw (-0.5,0) rectangle ++ (\textwidth-0.5,2.25);
        \draw[fill=MAGENTA] (0.25,0.25) rectangle ++(0.5,0.25);
        \node at (1.25,0.375) {$k=3$};
        \draw[fill=ORANGE] (0.25,1) rectangle ++(0.5,0.25);
        \node at (1.25,1.125) {$k=2$};
        \draw[fill=GREEN] (0.25,1.75) rectangle ++(0.5,0.25);
        \node at (1.25,1.875) {$k=1$};
        
        \draw[gray, dashed] (3,1.875) -- ++(0.5,0);
        \node[text width=4.75cm,align=flush left]  at (6.5,1.875) {no stabilization};
        \draw[gray] (3,1.125) -- ++(0.5,0);
        \node[text width=4.75cm,align=flush left]  at (6.5,1.125) {stabilization, $k'=k$};
        \draw[gray] (3,0.375) -- ++(0.5,0);
        \node[text width=4.75cm,align=flush left] at (6.5,0.375) {stabilization, $k'=k + 1$};
        \node[gray] at (3.25,0.375) {$\times$};
    
        \draw[gray] (9,1.875) -- ++(0.5,0);
        \node[gray] at (9.25,1.875) {$+$};
        \node[text width=4.75cm,align=flush left]  at (12.5,1.875) {stabilization, $k'=k+2$};
        \draw[gray] (9,1.125) -- ++(0.5,0);
        \draw[gray] (9.25,1.125) circle (0.075);
        \node[text width=4.75cm,align=flush left]  at (12.5,1.125) {stabilization, $k'=k+3$};
        \end{tikzpicture}
    \end{subfigure}
    \caption{The $L^2$ error convergence rates of the steady unregularized lid-driven cavity flow are shown above. This is performed for coarse-scale degrees $\vec{k} = (k,k)$ where $k=1,2,3$, with different degrees indicated by different colours. The different choices of fine-scale degree $\vec{k}' = (k',k')$ are instead indicated by different markers.
    Additionally, results without stabilization are shown with dashed lines.
    For a given coarse-scale degree, the lines overlap, indicating that the stabilized solution does not spoil the accuracy or the convergence rate.
    As later tests show, the stabilized solution does, however, improve the quality of the solution for under-resolved meshes.
    }
    \label{fig:NR-convergence-rate-lid-driven}
\end{figure}
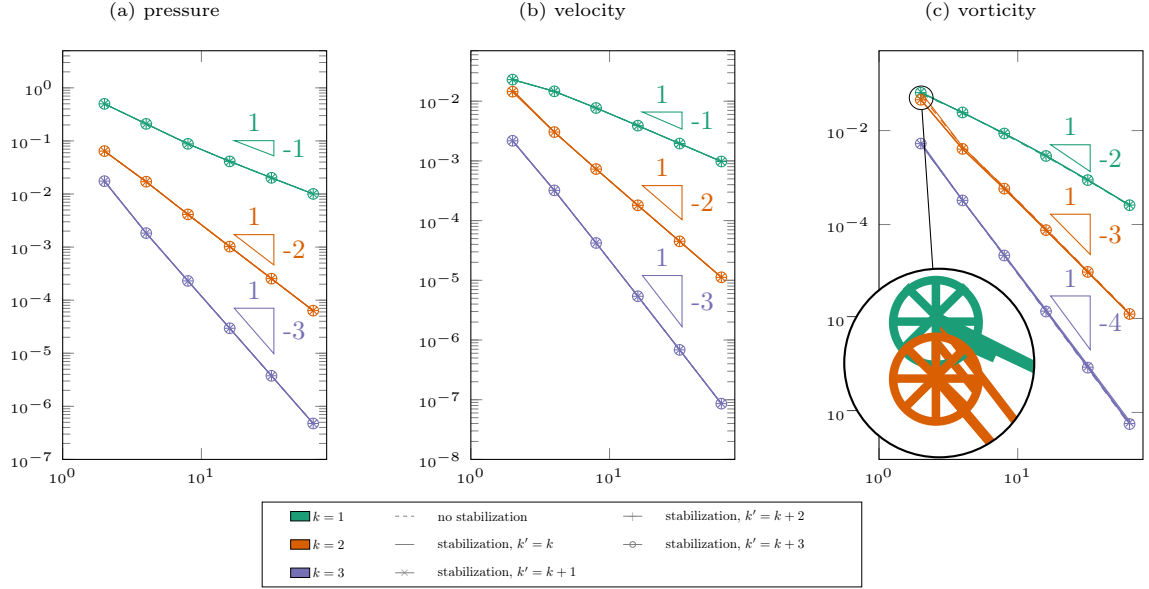

To analyze the steady-state case, we test convergence via a manufactured solution: the steady, regularized lid-driven cavity flow, initially proposed by \cite{shih_effects_1989}.
We pick the solution to be 
\begin{equation}
    \vec{u} = \begin{bmatrix}
        ( x^4 - 2 x^3 + x^2 ) ( 4 y^3 - 2 y)\\
        - ( 4 x^3 - 6 x^2 + 2 x ) ( y^4 - y^2 )
    \end{bmatrix}\;,~\omega = \mathrm{rot}(\vec{u})= \partial_x u_y - \partial_y u_x\;,~p= \sin(x)\sin(y)\;,
\end{equation}
and enforce appropriate Dirichlet boundary conditions on the velocity (by choosing the boundaries $\Gamma_{\hat{\vec{u}}} = \Gamma_{\hat{u}} = \partial \Omega$).
In Figure \ref{fig:NR-convergence-rate-lid-driven}, the convergence rates measured in $L^2$ norm can be seen for Reynolds number $\mathrm{Re} = 1000$.
Here, we investigate the effects of choosing different coarse-scale degrees $\vec{k} = (k,k), k = 1,2,3$, and of different fine-scale degrees $\vec{k}' = (k',k')$.
Note that since we require bubble functions for fine scales, we always have $k' \geq 2$.

We expect optimal convergence rates of $k+1$ for the vorticity and $k$ for the velocity and pressure solutions, and the numerical solutions confirm this; only the vorticity solution seems to show a slightly longer pre-asymptotic range for $k=1$.
In addition, these plots suggest that the choice of fine-scale degree $k'$ does not affect the convergence rate much.
Only for $k=2$ can we see some minor deviations, but for $k' \geq k + 1$ this has already stabilized.

\subsubsection{Unsteady two-dimensional Taylor--Green vortex flow}
We now look at the convergence by solving the unsteady two-dimensional Taylor--Green vortex, whose solution is analytically known \cite{evans_divergence-free_2011}. Here, we solve over the periodic domain $\Omega = [0,2\pi]^2$ by choosing the initial velocity as
\begin{equation}
    \vec{u}_{\text{initial}}(x,y) = \begin{bmatrix} \sin(x) \cos(y) \\
        -\cos(x) \sin(y)
    \end{bmatrix}\;.
\end{equation}
These equations (in rotational form) evolve as:
\begin{align}
    \omega(x,y,t) &= 2 \sin(x)\sin(y) e^{\frac{-2t}{\mathrm{Re}}}\;,\\
    \vec{u}(x,y,t) &= \begin{bmatrix} \sin(x) \cos(y) \\
        -\cos(x) \sin(y)
    \end{bmatrix}e^{\frac{-2 t}{\mathrm{Re}}}\;,\\
    p(x,y,t) &= \frac{\cos(2x) + \cos(2y)}{4}e^{\frac{-4t}{Re}} + \frac{\vec{u}(x,y,t)\cdot \vec{u}(x,y,t)}{2}\;.
\end{align}
Note that the pressure solution differs from that in \cite{evans_divergence-free_2011} by a factor $\frac{\vec{u}\cdot\vec{u}}{2}$ as, in the velocity-vorticity-pressure formulation, we solve for the total pressure $p = P + \frac{\vec{u}\cdot\vec{u}}{2}$, instead of the static pressure $P$.
We solve this problem for coarse-scale spline degrees $\vec{k}=(1,1),(2,2),(3,3)$, and fine-scale spline degrees $\vec{k}' = (k+1,k+1)$, and we choose the step size as the largest 
$\Delta t \leq \min\left(h^{\frac{k+1}{2}},\frac{h^2}{4}\mathrm{Re}\right)$ (taken from \cite{evans_divergence-free_2011}) that perfectly divides $T=1$. 

Note that, in our simulations, we initialize the fine-scale fields as zero, i.e., $\omega'^0 = 0,\vec{u}'^0 = \vec{0},p'^0 = 0$.
This is motivated by the fact that, in general, the analytical initial conditions are not known, making it impossible to initialize the fine-scale fields consistently.
Moreover, even if the analytical fields were known, initializing the fine-scale fields consistently would require an initial projection of the analytical solution onto the coarse- and fine-scale spaces, while preserving the divergence-free constraint for the velocities, which could be computationally expensive.
In the case of the Taylor--Green vortex, however, this choice has interesting consequences, as we will see next.

In Table \ref{tab:NR-Taylor--Green}, we show the $L^2$ errors for the coarse-scale velocity solution at time $T=1$ for both the unstabilized and stabilized formulations.
We observe optimal convergence rates of approximately $k$ for the velocity solution.
However, as we see from the provided $L^2$ norms of the fine-scale quantities, the fine-scale velocity and vorticity fields are always of the order of machine precision. As a result, the stabilized and unstabilized solutions are (almost) identical.
This seems to be a particular feature of the Taylor--Green vortex and our choice to initialize the fine-scale fields at zero.

For the Taylor--Green vortex solution, it holds that $\partial_t \vec{u}(t) + \mathrm{Re}^{-1}\curl{\omega(t)} = 0$, so that $\nabla p(t) = -\omega(t)\times \vec{u}(t)$.
Interestingly, {at each time-step $n$,} our stabilized discrete solution seems to mimic this behaviour when the fine-scale fields are initialized at zero:
\[ \left( \frac{\vec{u}^{h,n+1} - \vec{u}^{h,n}}{\Delta t},\vec{v}'\right) + \mathrm{Re}^{-1} \left( \curl{\omega^{h,n+\frac{1}{2}}},\vec{v}'\right)=0\;,\quad \forall \vec{v}'\in\mathbb{W}_h^1\;, \]
and 
\[ \left( \omega^{h,n+\frac{1}{2}}\times \vec{u}^{h,n+\frac{1}{2}},\vec{v}'\right) - \left( p^{h,n+\frac{1}{2}},\vec{v}'\right) - \left( p'^{n+\frac{1}{2}},\vec{v}'\right)=0\;,\quad \forall \vec{v}'\in\mathbb{W}_h^1\;. \]
{
As a result, the forcing term of the fine-scale problem ($(\mathrm{res_{M}},\vec{v}') = 0$) vanishes, and $\vec{u}'^{n+\frac{1}{2}} = \vec{0}$ and $\omega'^{n+\frac{1}{2}} = 0$ solve the fine-scale problem.
}

\begin{table}[htp]
\centering
\caption{The convergence of the unstabilized and VMS stabilized formulation for the Taylor--Green vortex benchmark. Observe that the error in the velocity is (almost) identical, due to the vanishing fine-scale fields. The calculated convergence rate is $0.99$ for $k=1$, $2.09$ for $k=2$ and $3.17$ for $k=3$.}
\label{tab:NR-Taylor--Green}
\begin{tabular}{rrrrrr}
\hline
\multicolumn{1}{l}{} & \multicolumn{5}{c}{Degree 1}                                                                                 \\
                     & \multicolumn{1}{c}{No stabilization}   & \multicolumn{4}{c}{VMS stabilization}                               \\ \hline
\# Elements          & $\Vert \vec{u}^h - \vec{u}\Vert_{L^2}$ & $\Vert \vec{u}^h - \vec{u}\Vert_{L^2}$ & $\Vert\vec{u}'\Vert_{L^2}$  & $\Vert \omega'\Vert_{L^2}$ & $\Vert p'\Vert_{L^2}$ \\ \hline
$4 \times 4$         & 1.92738413169754                       & 1.92738413169754                        & 6.04E-17  & 2.56E-16 & 1.69E+00  \\
$8 \times 8$         & 0.979519064239168                     & 0.979519064239168                     & 4.99E-16  &  1.20E-14 & 1.34E-01  \\
$16 \times 16$       & 0.492644572600598                      & 0.492644572600598                      & 7.48E-17  &   7.12E-16 & 5.12E-03
  \\ \hline \hline
\multicolumn{1}{l}{} & \multicolumn{5}{c}{Degree 2}                                                                                 \\
                     & \multicolumn{1}{c}{No stabilization}   & \multicolumn{4}{c}{VMS stabilization}                               \\  \hline
\# Elements          & $\Vert \vec{u}^h - \vec{u}\Vert_{L^2}$ & $\Vert \vec{u}^h - \vec{u}\Vert_{L^2}$ & $\Vert\vec{u}'\Vert_{L^2}$  & $\Vert \omega'\Vert_{L^2}$  & $\Vert p'\Vert_{L^2}$ \\ \hline
$4 \times 4$         & 0.547347069541981                      & 0.547347069541981                      & 3.90E-14  &   1.99E-14 & 1.43E-01 \\
$8 \times 8$         & 0.108380117382897                      & 0.108380117382897                      & 6.82E-14   &   4.07E-13 & 1.61E-02\\
$16 \times 16$       & 0.0255313734834109                     & 0.0255313734834109                     & 3.68E-13   &  4.46E-12 & 2.06E-03 \\ \hline \hline
\multicolumn{1}{l}{} & \multicolumn{5}{c}{Degree 3}                                                                                 \\
                     & \multicolumn{1}{c}{No stabilization}   & \multicolumn{4}{c}{VMS stabilization}                               \\ \hline
\# Elements          & $\Vert \vec{u}^h - \vec{u}\Vert_{L^2}$ & $\Vert \vec{u}^h - \vec{u}\Vert_{L^2}$ & $\Vert\vec{u}'\Vert_{L^2}$  & $\Vert\omega'\Vert_{L^2}$  & $\Vert p'\Vert_{L^2}$ \\ \hline
$4 \times 4$         & 0.174148833469405                      & 0.174148833469405                       & 2.21E-14  &    6.94E-14 & 1.33E+00\\
$8 \times 8$         & 0.0142579039840667                     & 0.0142579039840667                     & 8.71E-14   &  4.81E-13 & 4.02E-03 \\
$16 \times 16$       & 0.00158134529847338                    & 0.00158134529847338                    & 3.77E-14  &   6.08E-13 & 7.00E-05 \\ \hline
\end{tabular}
\end{table}

\subsection{Two-Dimensional shear layer roll-up problem}
\label{sec:NR-subgrid-analysis}
We investigate the effectiveness of the introduced stabilization on an unsteady problem: the two-dimensional shear layer roll-up problem.
Consider the periodic domain $\Omega = [0,2\pi]^2$, with initial conditions $\vec{u}_{\text{initial}} = [u_{x,\text{initial}} v_{y,\text{initial}}]^T$:
\begin{equation}\label{eq:NR-shear-layer-initial-cond}
    {u}_{x,\text{initial}}(x,y) = \begin{cases}
        \tanh{\left(\frac{y - \frac{\pi}{2}}{\delta}\right)} & y \leq \pi\;,\\
        \tanh{\left(\frac{\frac{3\pi}{2}- y}{\delta}\right)} & \text{else}\;,\\
    \end{cases}\;,\quad u_{y,\text{initial}} = \epsilon \sin(x)\;,
\end{equation}
where we choose $\delta = \frac{\pi}{15}$ and $\epsilon = 0.05$.
We will solve this problem for two cases: the inviscid case, where no dissipation is expected, and the viscous case, where kinetic energy dissipates.
For both cases, we investigate the kinetic energy of the coarse scale, fine scale, and the combined scales, defined as
\begin{equation}
    K_h(t) := \frac{1}{2}\Vert \vec{u}^h(t) \Vert^2\;,\quad K'(t) := \frac{1}{2}\Vert \vec{u}'(t) \Vert^2\;,\quad K(t) := \frac{1}{2}\Vert \vec{u}^h(t) +\vec{u}'(t) \Vert^2\;.
\end{equation}

\subsubsection{Inviscid Case}
We consider the inviscid case ($\mathrm{Re} = \infty$), where the kinetic energy should be conserved.
The energy evolution at the discrete level is governed by two effects: the ability of the initial coarse-scale space to represent the initial condition \eqref{eq:NR-shear-layer-initial-cond}, and the dissipation induced by the fine-scales.
In general, both effects will lead to a discrepancy between the continuous and discrete kinetic energies.
In particular, our formulation will not conserve energy as dissipation will occur through stabilization in the fine-scale governing equations; also compare Lemma \ref{lem:energy-evolution}.
With mesh refinement, we expect both effects to diminish.

\begin{figure}[htp]
    \centering
    \begin{subfigure}[b]{0.3\textwidth}
        \includegraphics[width=\textwidth]{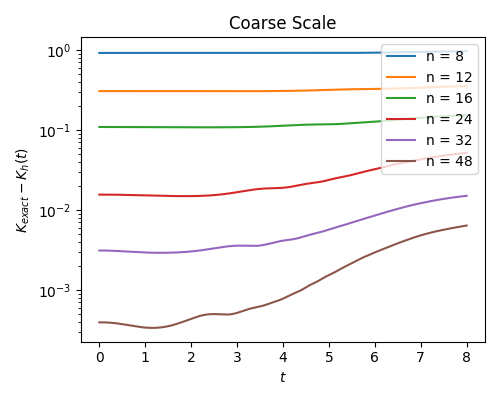}
        \caption{}
    \end{subfigure}
    \hfill
    \begin{subfigure}[b]{0.3\textwidth}
        \includegraphics[width=\textwidth]{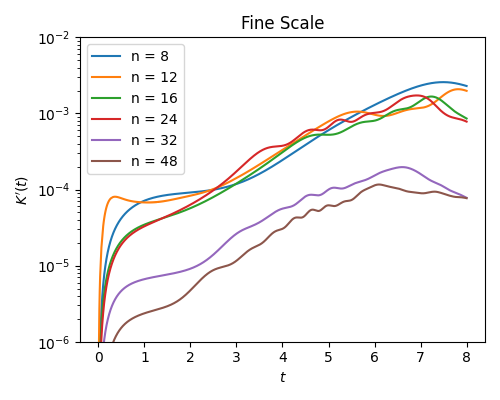}
        \caption{}
    \end{subfigure}
    \hfill
    \begin{subfigure}[b]{0.3\textwidth}
        \includegraphics[width=\textwidth]{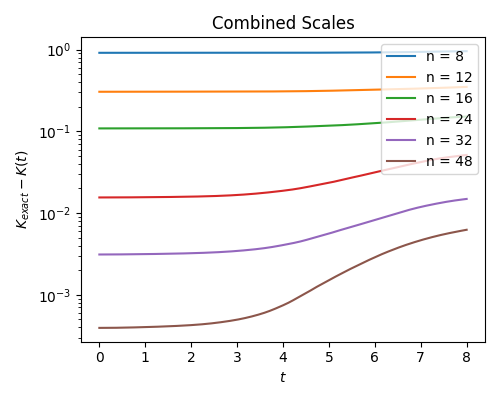}
        \caption{}
    \end{subfigure}
    \caption{For the inviscid case, the above plots show the changes in the coarse scale, fine scale, and total kinetic energies. The results for the mesh sizes $n=8,12,16,24,32,48$ with time step $\Delta t = 0.001$, with and without stabilization, are shown for degrees $k=2$ and $k'=3$.}
    \label{fig:inviscid-case-result-kinetic-energy}
\end{figure}

See Figure \ref{fig:inviscid-case-result-kinetic-energy}, where we solve the problem for various mesh element sizes with time step $\Delta t= 0.001$, and compare the kinetic energy to the reference value calculated from \eqref{eq:NR-shear-layer-initial-cond}.
Here, in line with the expectations mentioned above, we observe that, as the mesh is refined, the discrepancy between the initial energy (at $t=0.0$) at the continuous and discrete levels decreases, as the space can better approximate \eqref{eq:NR-shear-layer-initial-cond}.
Moreover, as the mesh is refined, the dissipation induced by the fine-scales is reduced, which can be seen a) by the lower amount of kinetic energy stored in the fine-scales, and b) the fact that the total kinetic energy stays closer to $K_{exact}$ as it evolves.

\subsubsection{Viscid case, $\mathrm{Re}=1600$.}
\begin{figure}
    \centering
    \includegraphics[width=0.3\linewidth]{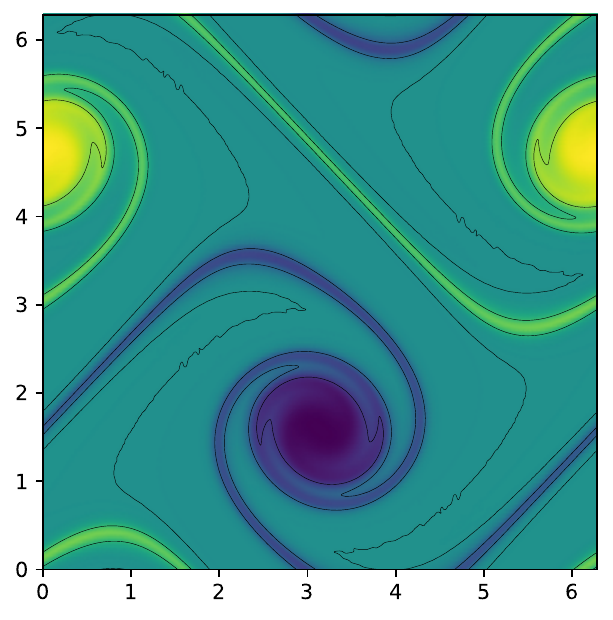}
    \caption{The vorticity reference solution for the viscid shear layer roll up ($\mathrm{Re} = 1600$), solved with cubic B-splines on a mesh of $128\times128$ elements. The contour lines are given by $-6,-5,\dots,5,6$.}
    \label{fig:viscid-case-reference-solution}
\end{figure}
In the viscous case, the exact solution also dissipates energy, so there is no analytical reference.
Instead, we solve the problem over a fine mesh ($n=128$) and take this as the benchmark for comparisons.
For stabilization, we will also investigate the implicit and semi-implicit formulations as discussed in Section \ref{sec:choice-stabilization-parameter}.
For this, all simulations use a time step $\Delta t = 0.01$.

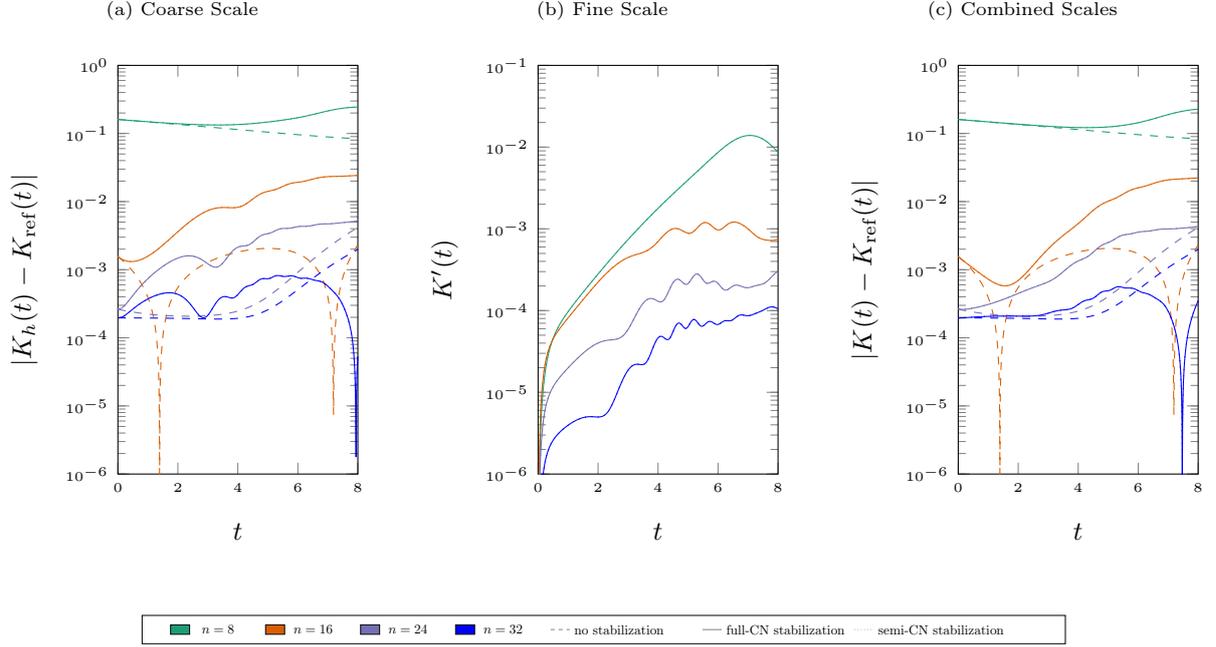
\begin{figure}[htp]
    \centering
    \begin{subfigure}{0.3\textwidth}
        \caption{Coarse Scale}
        \label{fig:NR-viscid-Re1600-coarse}
        \newcommand{\logLogSlopeTriangle}[5]
{

    \pgfplotsextra
    {
        \pgfkeysgetvalue{/pgfplots/xmin}{\xmin}
        \pgfkeysgetvalue{/pgfplots/xmax}{\xmax}
        \pgfkeysgetvalue{/pgfplots/ymin}{\ymin}
        \pgfkeysgetvalue{/pgfplots/ymax}{\ymax}

        \pgfmathsetmacro{\xArel}{#1}
        \pgfmathsetmacro{\yArel}{#3}
        \pgfmathsetmacro{\xBrel}{#1-#2}
        \pgfmathsetmacro{\yBrel}{\yArel}
        \pgfmathsetmacro{\xCrel}{\xArel}

        \pgfmathsetmacro{\lnxB}{\xmin*(1-(#1-#2))+\xmax*(#1-#2)} 
        \pgfmathsetmacro{\lnxA}{\xmin*(1-#1)+\xmax*#1} 
        \pgfmathsetmacro{\lnyA}{\ymin*(1-#3)+\ymax*#3} 
        \pgfmathsetmacro{\lnyC}{\lnyA+#4*(\lnxA-\lnxB)}
        \pgfmathsetmacro{\yCrel}{\lnyC-\ymin)/(\ymax-\ymin)} 

        \coordinate (A) at (rel axis cs:\xArel,\yArel);
        \coordinate (B) at (rel axis cs:\xBrel,\yBrel);
        \coordinate (C) at (rel axis cs:\xCrel,\yCrel);

        \draw[#5]   (A)-- node[pos=0.5,anchor=south] {1}
                    (B)-- 
                    (C)-- node[pos=0.5,anchor=west] {#4}
                    cycle;
    }
}

\begin{tikzpicture}[spy using outlines=
{circle, magnification=8, connect spies}]
\begin{semilogyaxis}[%
  xmin=0, xmax=8,
  ymin=10^(-6), ymax=1,
  width=\textwidth,
  height=7cm,
  xlabel=$t$,
  ylabel=$\left|K_h(t) - K_{\mathrm{ref}}(t)\right|$,
  ticklabel style={font=\tiny}]
    \addplot[GREEN, style = dashed] table [x=time, y=8_False_full, col sep=comma] {Images/ViscidAnalysis/coarse_data_abs.csv};
    \addplot[GREEN] table [x=time, y=8_True_full, col sep=comma] {Images/ViscidAnalysis/coarse_data_abs.csv};
    \addplot[GREEN, style = dotted] table [x=time, y=8_True_semi, col sep=comma] {Images/ViscidAnalysis/coarse_data_abs.csv};

    \addplot[ORANGE, style = dashed] table [x=time, y=16_False_full, col sep=comma] {Images/ViscidAnalysis/coarse_data_abs.csv};
    \addplot[ORANGE] table [x=time, y=16_True_full, col sep=comma] {Images/ViscidAnalysis/coarse_data_abs.csv};
    \addplot[ORANGE, style = dotted] table [x=time, y=16_True_semi, col sep=comma] {Images/ViscidAnalysis/coarse_data_abs.csv};

    \addplot[MAGENTA, style = dashed] table [x=time, y=24_False_full, col sep=comma] {Images/ViscidAnalysis/coarse_data_abs.csv};
    \addplot[MAGENTA] table [x=time, y=24_True_full, col sep=comma] {Images/ViscidAnalysis/coarse_data_abs.csv};
    \addplot[MAGENTA, style = dotted] table [x=time, y=24_True_semi, col sep=comma] {Images/ViscidAnalysis/coarse_data_abs.csv};

    \addplot[blue, style = dashed] table [x=time, y=32_False_full, col sep=comma] {Images/ViscidAnalysis/coarse_data_abs.csv};
    \addplot[blue] table [x=time, y=32_True_full, col sep=comma] {Images/ViscidAnalysis/coarse_data_abs.csv};
    \addplot[blue, style = dotted] table [x=time, y=32_True_semi, col sep=comma] {Images/ViscidAnalysis/coarse_data_abs.csv};

\end{semilogyaxis}

\end{tikzpicture}%
    \end{subfigure}
    \hfill
    \begin{subfigure}{0.3\textwidth}
        \caption{Fine Scale}
        \label{fig:NR-viscid-Re1600-fine}
        \newcommand{\logLogSlopeTriangle}[5]
{

    \pgfplotsextra
    {
        \pgfkeysgetvalue{/pgfplots/xmin}{\xmin}
        \pgfkeysgetvalue{/pgfplots/xmax}{\xmax}
        \pgfkeysgetvalue{/pgfplots/ymin}{\ymin}
        \pgfkeysgetvalue{/pgfplots/ymax}{\ymax}

        \pgfmathsetmacro{\xArel}{#1}
        \pgfmathsetmacro{\yArel}{#3}
        \pgfmathsetmacro{\xBrel}{#1-#2}
        \pgfmathsetmacro{\yBrel}{\yArel}
        \pgfmathsetmacro{\xCrel}{\xArel}

        \pgfmathsetmacro{\lnxB}{\xmin*(1-(#1-#2))+\xmax*(#1-#2)} 
        \pgfmathsetmacro{\lnxA}{\xmin*(1-#1)+\xmax*#1} 
        \pgfmathsetmacro{\lnyA}{\ymin*(1-#3)+\ymax*#3} 
        \pgfmathsetmacro{\lnyC}{\lnyA+#4*(\lnxA-\lnxB)}
        \pgfmathsetmacro{\yCrel}{\lnyC-\ymin)/(\ymax-\ymin)} 

        \coordinate (A) at (rel axis cs:\xArel,\yArel);
        \coordinate (B) at (rel axis cs:\xBrel,\yBrel);
        \coordinate (C) at (rel axis cs:\xCrel,\yCrel);

        \draw[#5]   (A)-- node[pos=0.5,anchor=south] {1}
                    (B)-- 
                    (C)-- node[pos=0.5,anchor=west] {#4}
                    cycle;
    }
}

\begin{tikzpicture}[spy using outlines=
{circle, magnification=80, connect spies}]
\begin{semilogyaxis}[%
  at = {(0cm,0cm)},
  xmin=0, xmax=8,
  ymin=10^(-6), ymax=0.1,
  width=\textwidth,
  height=7cm,
    xlabel=$t$,
  ylabel=$K'(t)$,
  ticklabel style={font=\tiny}]
    \addplot[GREEN] table [x=time, y=8_True_full, col sep=comma] {Images/ViscidAnalysis/fine_data.csv};
    \addplot[GREEN, style = dotted] table [x=time, y=8_True_semi, col sep=comma] {Images/ViscidAnalysis/fine_data.csv};

    \addplot[ORANGE] table [x=time, y=16_True_full, col sep=comma] {Images/ViscidAnalysis/fine_data.csv};
    \addplot[ORANGE, style = dotted] table [x=time, y=16_True_semi, col sep=comma] {Images/ViscidAnalysis/fine_data.csv};
    
    \addplot[MAGENTA] table [x=time, y=24_True_full, col sep=comma] {Images/ViscidAnalysis/fine_data.csv};
    \addplot[MAGENTA, style = dotted] table [x=time, y=24_True_semi, col sep=comma] {Images/ViscidAnalysis/fine_data.csv};

    \addplot[blue] table [x=time, y=32_True_full, col sep=comma] {Images/ViscidAnalysis/fine_data.csv};
    \addplot[blue, style = dotted] table [x=time, y=32_True_semi, col sep=comma] {Images/ViscidAnalysis/fine_data.csv};

\end{semilogyaxis}


\end{tikzpicture}%
    \end{subfigure}
    \hfill
    \begin{subfigure}{0.3\textwidth}
        \caption{Combined Scales}
        \label{fig:NR-viscid-Re1600-total}
        \newcommand{\logLogSlopeTriangle}[5]
{

    \pgfplotsextra
    {
        \pgfkeysgetvalue{/pgfplots/xmin}{\xmin}
        \pgfkeysgetvalue{/pgfplots/xmax}{\xmax}
        \pgfkeysgetvalue{/pgfplots/ymin}{\ymin}
        \pgfkeysgetvalue{/pgfplots/ymax}{\ymax}

        \pgfmathsetmacro{\xArel}{#1}
        \pgfmathsetmacro{\yArel}{#3}
        \pgfmathsetmacro{\xBrel}{#1-#2}
        \pgfmathsetmacro{\yBrel}{\yArel}
        \pgfmathsetmacro{\xCrel}{\xArel}

        \pgfmathsetmacro{\lnxB}{\xmin*(1-(#1-#2))+\xmax*(#1-#2)} 
        \pgfmathsetmacro{\lnxA}{\xmin*(1-#1)+\xmax*#1} 
        \pgfmathsetmacro{\lnyA}{\ymin*(1-#3)+\ymax*#3} 
        \pgfmathsetmacro{\lnyC}{\lnyA+#4*(\lnxA-\lnxB)}
        \pgfmathsetmacro{\yCrel}{\lnyC-\ymin)/(\ymax-\ymin)} 

        \coordinate (A) at (rel axis cs:\xArel,\yArel);
        \coordinate (B) at (rel axis cs:\xBrel,\yBrel);
        \coordinate (C) at (rel axis cs:\xCrel,\yCrel);

        \draw[#5]   (A)-- node[pos=0.5,anchor=south] {1}
                    (B)-- 
                    (C)-- node[pos=0.5,anchor=west] {#4}
                    cycle;
    }
}

\begin{tikzpicture}[spy using outlines=
{circle, magnification=8, connect spies}]
\begin{semilogyaxis}[%
  xmin=0, xmax=8,
  ymin=10^(-6), ymax=1,
  width=\textwidth,
  height=7cm,
    xlabel=$t$,
  ylabel=$\left|K(t) - K_{\mathrm{ref}}(t)\right|$,
  ticklabel style={font=\tiny}]
    \addplot[GREEN, style = dashed] table [x=time, y=8_False_full, col sep=comma] {Images/ViscidAnalysis/total_data_abs.csv};
    \addplot[GREEN] table [x=time, y=8_True_full, col sep=comma] {Images/ViscidAnalysis/total_data_abs.csv};
    \addplot[GREEN, style = dotted] table [x=time, y=8_True_semi, col sep=comma] {Images/ViscidAnalysis/total_data_abs.csv};

    \addplot[ORANGE, style = dashed] table [x=time, y=16_False_full, col sep=comma] {Images/ViscidAnalysis/total_data_abs.csv};
    \addplot[ORANGE] table [x=time, y=16_True_full, col sep=comma] {Images/ViscidAnalysis/total_data_abs.csv};
    \addplot[ORANGE, style = dotted] table [x=time, y=16_True_semi, col sep=comma] {Images/ViscidAnalysis/total_data_abs.csv};

    \addplot[MAGENTA, style = dashed] table [x=time, y=24_False_full, col sep=comma] {Images/ViscidAnalysis/total_data_abs.csv};
    \addplot[MAGENTA] table [x=time, y=24_True_full, col sep=comma] {Images/ViscidAnalysis/total_data_abs.csv};
  ymin=-0.01, ymax=0.01,
    \addplot[MAGENTA, style = dotted] table [x=time, y=24_True_semi, col sep=comma] {Images/ViscidAnalysis/total_data_abs.csv};

    \addplot[blue, style = dashed] table [x=time, y=32_False_full, col sep=comma] {Images/ViscidAnalysis/total_data_abs.csv};
    \addplot[blue] table [x=time, y=32_True_full, col sep=comma] {Images/ViscidAnalysis/total_data_abs.csv};
    \addplot[blue, style = dotted] table [x=time, y=32_True_semi, col sep=comma] {Images/ViscidAnalysis/total_data_abs.csv};

\end{semilogyaxis}

\end{tikzpicture}%
    \end{subfigure}\\
    \begin{subfigure}{\textwidth}
        \hspace{0.125\textwidth}
        \begin{tikzpicture}[scale=0.5]
        \pgflowlevelsynccm
        \draw (-0.5,0) rectangle ++ (\textwidth+8.5cm,0.75);
        \draw[fill=GREEN] (0.25,0.25) rectangle ++(0.5,0.25);
        \node at (1.5,0.375) {$n=8$};
        \draw[fill=ORANGE] (2.75,0.25) rectangle ++(0.5,0.25);
        \node at (4,0.375) {$n=16$};
        \draw[fill=MAGENTA] (5.25,0.25) rectangle ++(0.5,0.25);
        \node at (6.5,0.375) {$n=24$};
        \draw[fill=blue] (7.75,0.25) rectangle ++(0.5,0.25);
        \node at (9,0.375) {$n=32$};
        
        \draw[gray, dashed] (10.25,0.375) -- ++(0.5,0);
        \node[text width=4.75cm,align=flush left]  at (13.25,0.375) {no stabilization};
        \draw[gray] (14.25,0.375) -- ++(0.5,0);
        \node[text width=4.75cm,align=flush left]  at (17.25,0.375) {full-CN stabilization};
        \draw[gray,dotted] (18.25,0.375) -- ++(0.5,0);
        \node[text width=4.75cm,align=flush left]  at (21.25,0.375) {semi-CN stabilization};
        \end{tikzpicture}
    \end{subfigure}
    \caption{The differences of the kinetic energies compared to the reference ($n=128$) in the viscous case for $\mathrm{Re}=1600$ are shown.
    The results are for the mesh sizes $n=8,16,24,32$ with time step $\Delta t = 0.01$, degrees $ k=3$ and $ k'=4$, and for full-CN and semi-CN stabilization choices. Note that the semi-CN and full-CN stabilization results are so close that they overlap.}
    \label{fig:NR-viscid-Re1600}
\end{figure}

\begin{figure}[htp]
    \centering
    \begin{subfigure}[b]{0.24\textwidth}
        \includegraphics[width=\textwidth]{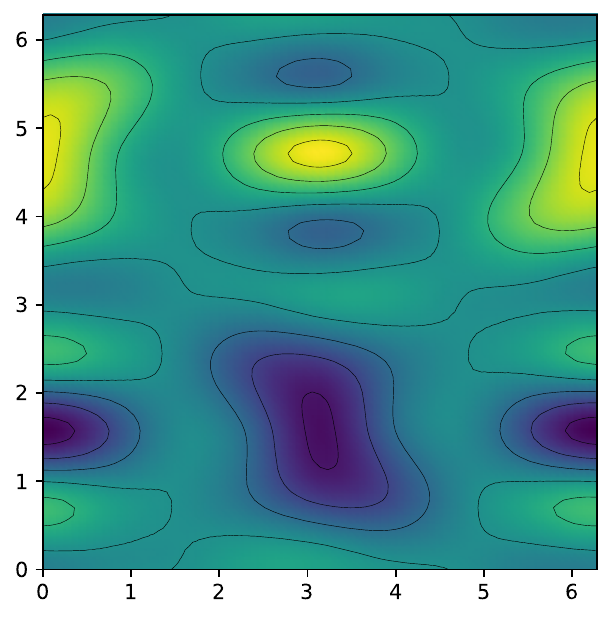}
        \caption{$n=8$, No stabilization}
        \label{fig:NR-viscid-Re1600-vor-8-false}
    \end{subfigure}
    \hfill
    \begin{subfigure}[b]{0.24\textwidth}
        \includegraphics[width=\textwidth]{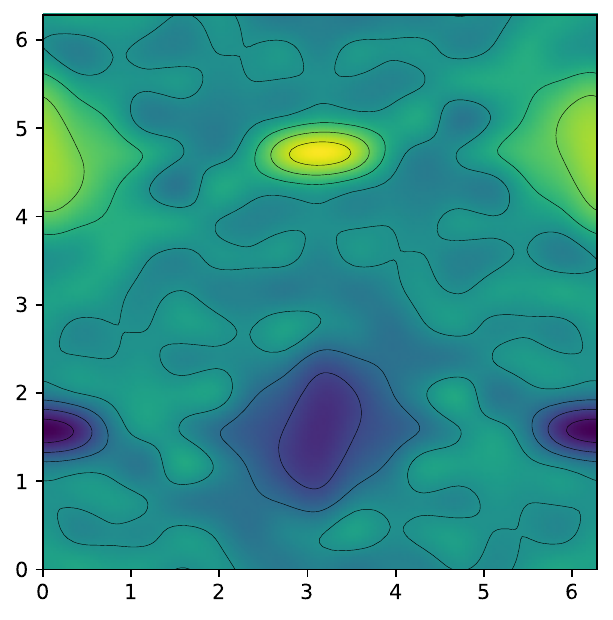}
        \caption{$n=16$, No stabilization}
        \label{fig:NR-viscid-Re1600-vor-16-false}
    \end{subfigure}
    \hfill
    \begin{subfigure}[b]{0.24\textwidth}
        \includegraphics[width=\textwidth]{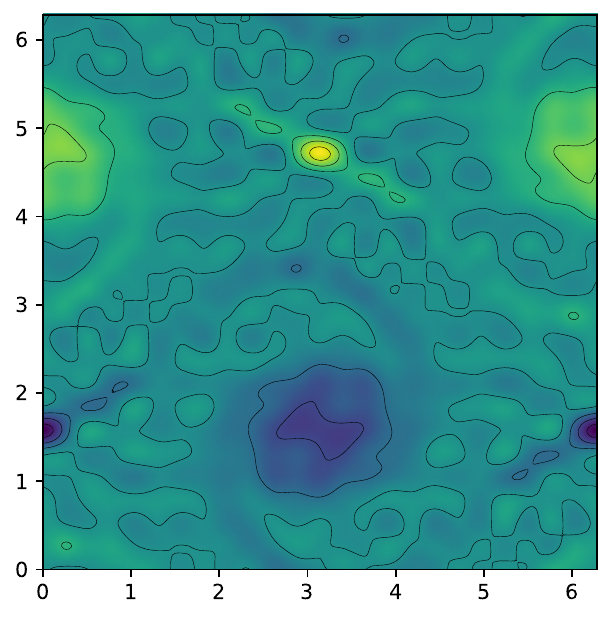}
        \caption{$n=24$, No stabilization}
        \label{fig:NR-viscid-Re1600-vor-24-false}
    \end{subfigure}
    \hfill
    \begin{subfigure}[b]{0.24\textwidth}
        \includegraphics[width=\textwidth]{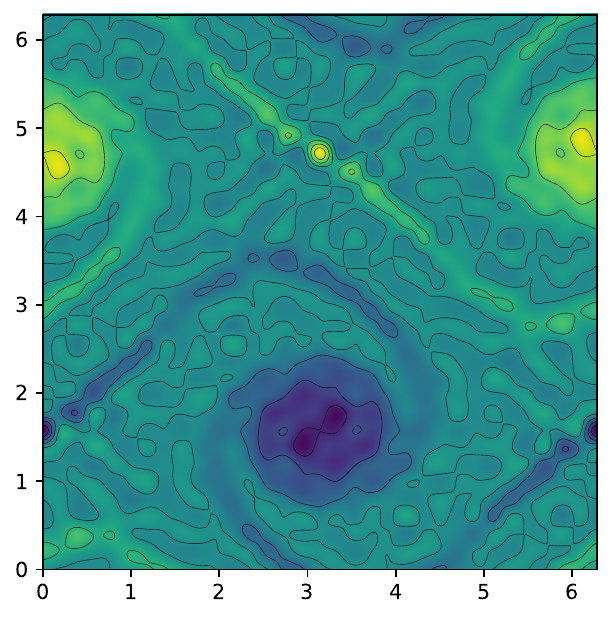}
        \caption{$n=32$, No stabilization}
        \label{fig:NR-viscid-Re1600-vor-32-false}
    \end{subfigure}\\
    \begin{subfigure}[b]{0.24\textwidth}
        \includegraphics[width=\textwidth]{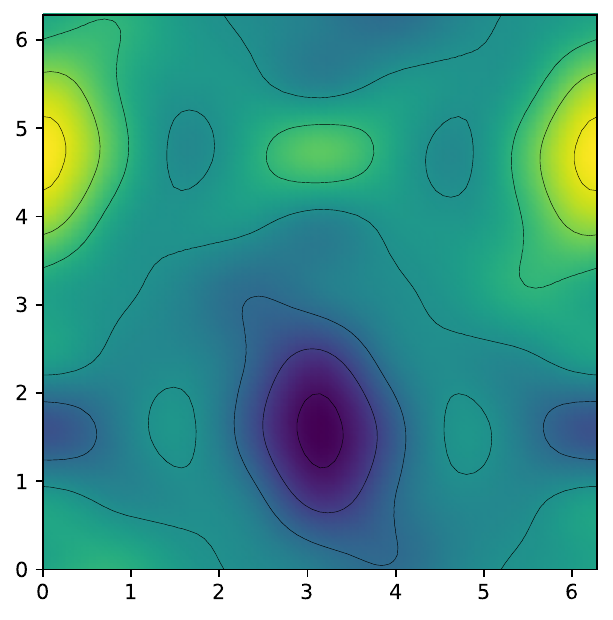}
        \caption{$n=8$, Stabilization}
        \label{fig:NR-viscid-Re1600-vor-8-true}
    \end{subfigure}
    \hfill
    \begin{subfigure}[b]{0.24\textwidth}
        \includegraphics[width=\textwidth]{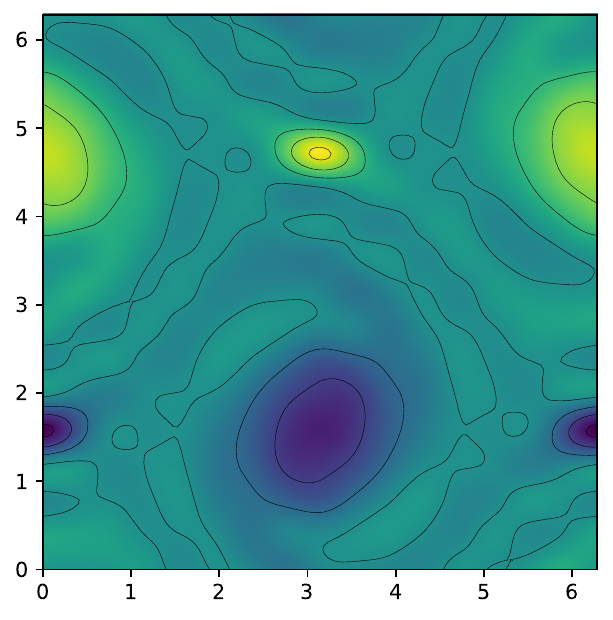}
        \caption{$n=16$, Stabilization}
        \label{fig:NR-viscid-Re1600-vor-16-true}
    \end{subfigure}
    \hfill
    \begin{subfigure}[b]{0.24\textwidth}
        \includegraphics[width=\textwidth]{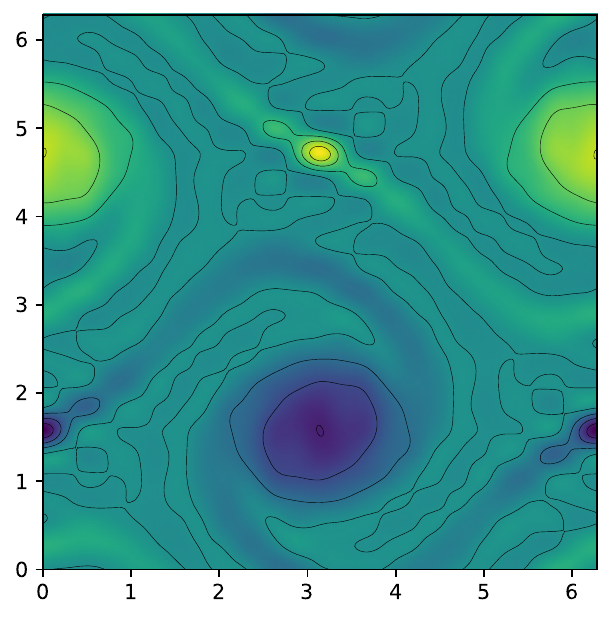}
        \caption{$n=24$, Stabilization}
        \label{fig:NR-viscid-Re1600-vor-24-true}
    \end{subfigure}
    \hfill
    \begin{subfigure}[b]{0.24\textwidth}
        \includegraphics[width=\textwidth]{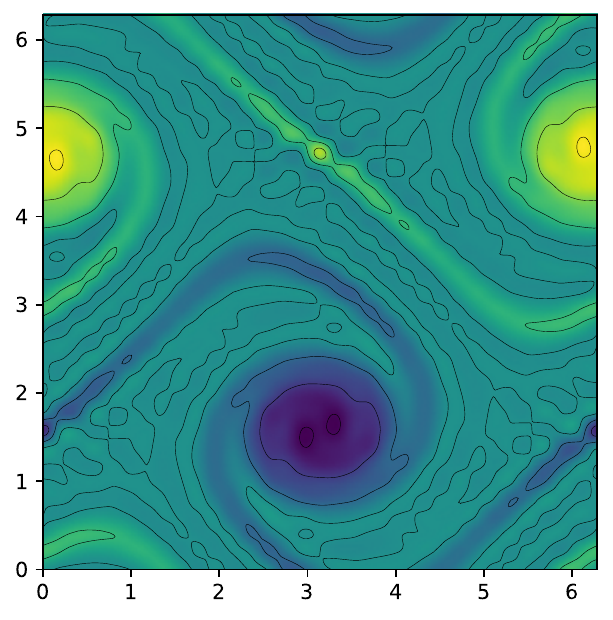}
        \caption{$n=32$, Stabilization}
        \label{fig:NR-viscid-Re1600-vor-32-true}
    \end{subfigure}
    \caption{The vorticity at $T=8.0$ for the viscous case ($\mathrm{Re}=1600$). Observe that the stabilized solutions are a better qualitative approximation of the high-resolution solution. The contour levels are given by $-6, -5, \dots, 5, 6$.}
    \label{fig:viscid-case-result-kinetic-energy}
\end{figure}

The results of kinetic energy evolution are shown in Figure \ref{fig:NR-viscid-Re1600}, where we compare the different computed kinetic energies with the reference energy $K_{\mathrm{ref}}(t)$.
As expected, except for the $n=32$ case, the kinetic energy of the unstabilized simulation is slightly closer to the reference solution, indicating more dissipation in the stabilized cases.
(As the mesh is refined, both get closer to the reference, as expected.)
Note that the semi-CN and full-CN results are indistinguishable, showing that both formulations are almost identical.

Remarkably, despite the slightly increased dissipation, the stabilized solutions are qualitatively much closer to the reference in all cases, see Figure \ref{fig:NR-viscid-Re1600}.
For example, comparing the $n=8$ vorticity fields in Figure \ref{fig:NR-viscid-Re1600-vor-8-true} and Figure \ref{fig:NR-viscid-Re1600-vor-8-false} with the reference solution in Figure \ref{fig:viscid-case-reference-solution}, we note that the unstabilized solution shows many small additional vortices.
Instead, the stabilized solution accurately predicts the correct number of large vortices and their location.
For $n=16,24,32$, the unstabilized solution again predicts the correct number of vortices, but does not resolve the `arms' spiraling from the vortices as accurately as the stabilized solutions.
In addition, the stabilized solutions are less oscillatory compared to the unstabilized solutions.

\subsection{Two-dimensional lid-driven cavity problem}
\label{sec:NR-comparison}
\begin{figure}[htp]
    \centering
    \begin{subfigure}{0.3\textwidth}
        \includegraphics[width = \textwidth]{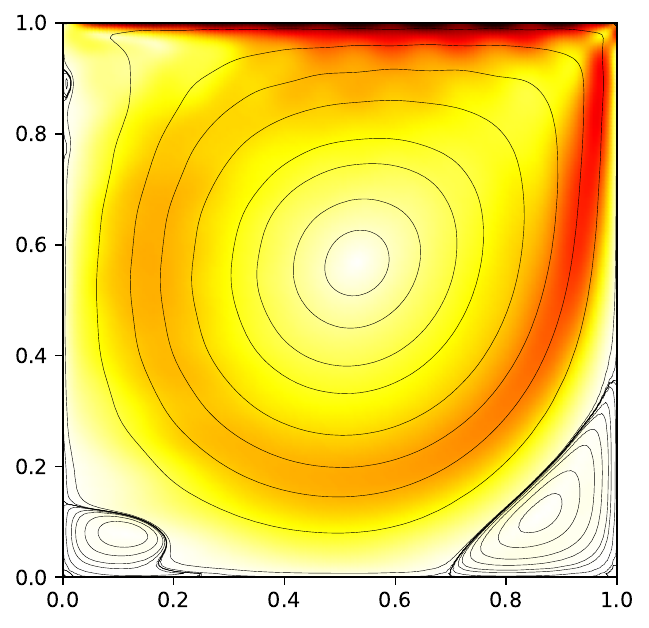}
    \end{subfigure}
    \hfill
    \begin{subfigure}{0.3\textwidth}
        \includegraphics[width = \textwidth]{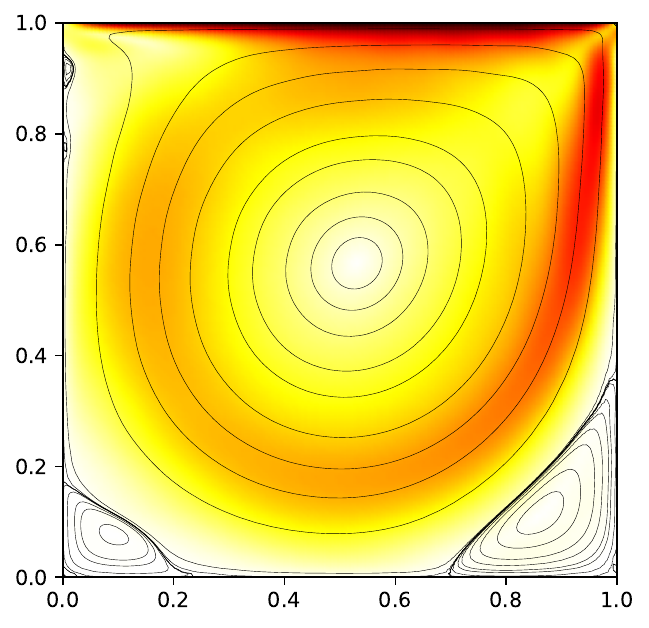}
    \end{subfigure}
    \hfill
    \begin{subfigure}{0.3\textwidth}
        \includegraphics[width = \textwidth]{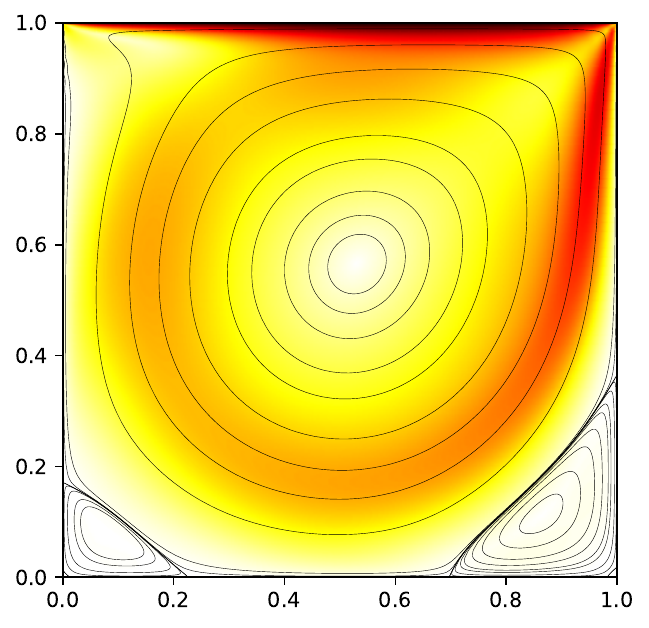}
    \end{subfigure}
    \begin{subfigure}{0.3\textwidth}
        \includegraphics[width = \textwidth]{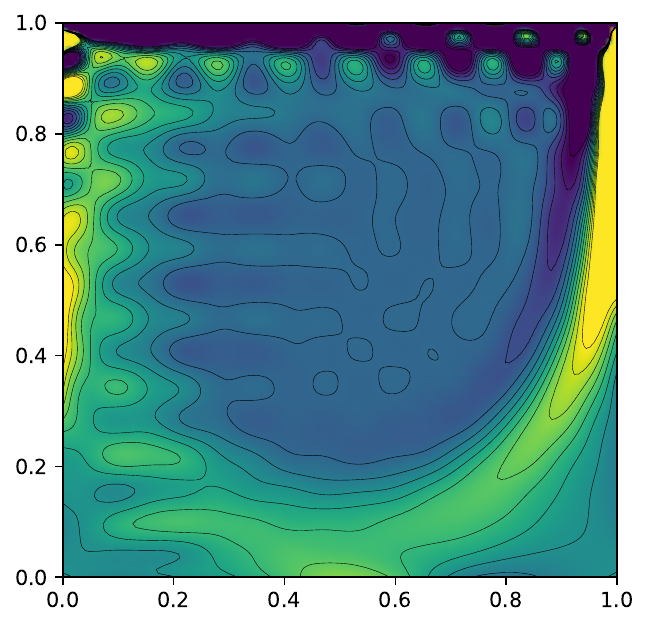}
        \caption{$16\times16$, no stabilization}
    \end{subfigure}
    \hfill
    \begin{subfigure}{0.3\textwidth}
        \includegraphics[width = \textwidth]{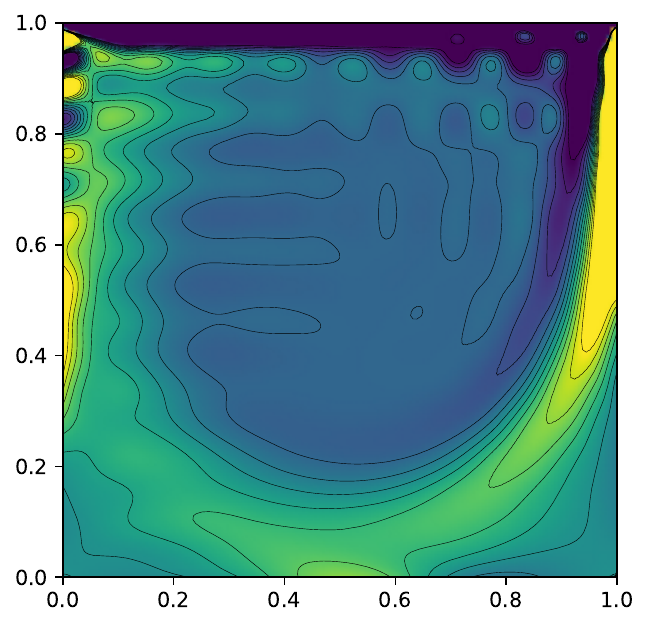}
        \caption{$16\times16$, stabilization}
    \end{subfigure}
    \hfill
    \begin{subfigure}{0.3\textwidth}
        \includegraphics[width = \textwidth]{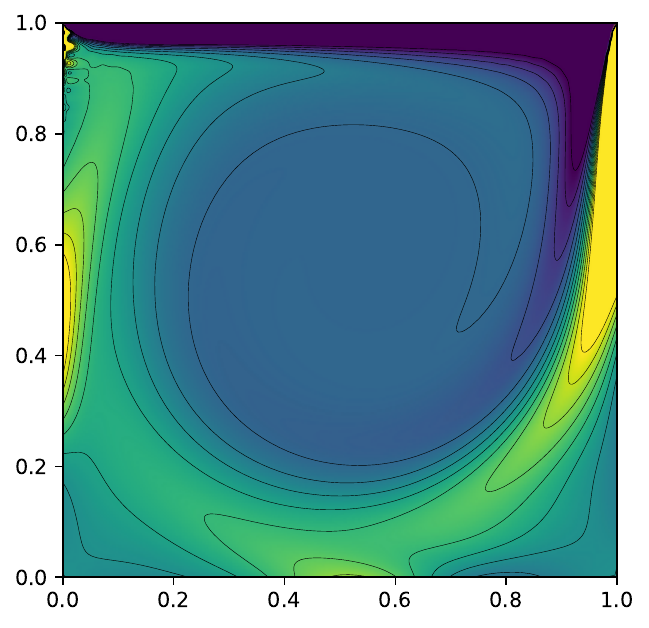}
        \caption{$64\times64$, no stabilization}
    \end{subfigure}
    \caption{The velocity and vorticity solutions of the two-dimensional steady-state lid-driven cavity flow. Here we use $\vec{k} = (4,4)$ and $\vec{k}'=(5,5)$. In the top figures, the magnitude of the velocity is shown in color, and the contour lines depict the streamfunction. The height of the contour lines is taken from \cite[table 7]{botellaBenchmarkSpectralResults1998}. The bottom figure depicts the vorticity in colour and contour lines, with height levels given by $-6, -5, \dots, 5, 6$.  }
    \label{fig:NR-qualitative}
\end{figure}

Lastly, we compare the solutions obtained with and without stabilization for the two-dimensional lid-driven cavity problem.
We solve this problem over the unit square $\Omega=[0,1]^2$, with Dirichlet boundary condition $\vec{u} = [u_x, u_y]$, where:
\begin{align}
    u_x(x,y) &= \begin{cases}
        1 & \text{if }y=1\;,\\
        0 &\text{else}\;
    \end{cases}\\
    u_y(x,y) &= 0\;.
\end{align}
and $\Gamma_{\hat{\vec{u}}} = \Gamma_{\hat{u}} = \partial\Omega$.
We choose a Reynolds number $\mathrm{Re} = 1000$ and bi-variation spline degree $\vec{k} = (4,4)$. 
In Figure \ref{fig:NR-qualitative}, the velocity (with streamlines), and vorticity solutions are shown for a low-resolution mesh $16\times16$ (stabilization and no stabilization) and a high-resolution mesh $64\times 64$ (no stabilization). 
Here, we again observe that the stabilized low-resolution solution more closely resembles the high-resolution solution than the unstabilized low-resolution simulation, better capturing the corner vortices and being less oscillatory.

\section{Conclusion}
We have introduced a variational multiscale stabilization formulation of the vorticity-velocity-pressure formulation of the incompressible Navier--Stokes equations, where both the coarse- and fine-scale discretization conform to the structure of the continuous de Rham complex.
Our proposed formulation is residual-based and triggers the coarse and fine decomposition using the Stokes projector \eqref{eq:stokes-projector}.
The fully discrete formulation uses modelling assumptions similar to those from the variational multiscale literature \cite{van_opstal_isogeometric_2017,evans_variational_2020}, and simplifies the fine-scale equations by using a Darcy--Stokes model.
Finally, the fine-scale problem is further simplified by decoupling the element problems and using bubble spaces for their discretization.

The resulting discrete formulation provides several desirable properties.
We show how the fine-scale equations can be solved element-wise in parallel, with both Picard and Newton-type linearizations.
We prove an explicit energy evolution for the fully discrete problem, thus also showing energetic stability.
Both the computed coarse- and fine-scale velocities are pointwise divergence-free.
For an Oseen problem related to our formulation, we prove stability, pressure-robustness, and optimal a priori estimates.
Finally, we numerically validated our approach for several incompressible Navier--Stokes benchmark problems in two dimensions.
The numerical results demonstrate that the proposed stabilization leads to improved qualitative results compared to the unstabilized formulation on coarse meshes, with no impact on its quantitative accuracy and convergence.

For future research, we aim to extend this formulation to the incompressible magnetohydrodynamics equations by supplementing it with Maxwell's equations. 
Additionally, we want to explore the use of the fine-scale solutions to drive adaptive refinement in combination with adaptive discretizations of the de Rham complex \cite{shepherd2024locally,dijkstra2025macro}; see Remark \ref{rem:adaptivity}.
Using a fine-scale solution for adaptivity has been done before; see, for example, \cite{irisarri_posteriori_2020}, and it has the benefit of requiring no additional computational cost -- the a posteriori error estimate is embedded in the VMS formulation.
Finally, our current formulation provides stabilization across the entire domain and utilizes the same bubble spaces for the fine-scale discretization everywhere.
However, since our fine-scale problems are localized, both of these aspects could be adapted based on local flow features (for instance, with local degree or mesh adaptivity, or turning off stabilization in certain regions), leading to a more efficient approach toward stabilization.
We will explore these aspects in future research.

\section{Acknowledgements}
The research of KD is supported by the Peter Paul Peterich / DIAM Fast Track Scholarship.
The research of DT was partially supported by project number 202.150 awarded through the Veni research programme by the Dutch Research Council (NWO).

\bibliography{references_fixed}{}
\bibliographystyle{plain}
\appendix

\section{Proof of Lemma \ref{lem:stability}}
\label{app:proof-stability-td}
The proof of Lemma \ref{lem:stability} is shown in three parts. 
We start by noting that the boundedness of $\mathcal{A}_{\mathrm{red}}$ is trivial from its definition, leaving the proofs of conditions \eqref{eq:stability-infsup} and \eqref{eq:stability-nondegenerate}.

\subsection{Boundedness \eqref{eq:bounded-form}}
Let $x,y \in \mathbb{X}$ where $x = (\vec{\omega}^h,\vec{u}^h,\vec{\omega}',\vec{u}')$ and $y = (\vec{\tau}^h,\vec{v}^h,\vec{\tau}',\vec{v}')$ and estimate all but one of the terms of \eqref{method:reduced-oseen-stabilization} straightforwardly by
\begin{subequations}
    \begin{align}
        \vert\sigma (\vec{u}^h + \vec{u}',\vec{v}^h + \vec{v}')\vert &\leq \sigma \Vert \vec{u}^h +\vec{u}'\Vert \Vert \vec{v}^h + \vec{v}'\Vert\;,\\
        \vert\frac{1}{\sqrt{\nu}}(\vec{\omega}^h\times\left[\vec{\beta}^h+\vec{\beta}'\right],\vec{v}^h+\vec{v}')\vert &\leq \frac{2\Vert \vec{\beta}^h+\vec{\beta}'\Vert_{\infty} }{\sqrt{\nu}}\Vert \vec{\omega}^h\Vert \Vert \vec{v}^h +\vec{v}' \Vert\;,\\ 
        \vert\sqrt{\nu}(\curl{\vec{\omega}^h},\vec{v}^h + \vec{v}') \vert&\leq \sqrt{\nu}\Vert\curl{\vec{\omega}^h}\Vert \Vert\vec{v}^h + \vec{v}'\Vert\;, \\
     \vert(\vec{\omega}^h,\vec{\tau}^h)\vert &\leq \Vert \vec{\omega}^h\Vert \Vert\vec{\tau}^h\Vert\;, \\
     \vert\sqrt{\nu}(\vec{u}^h + \vec{u}',\curl{\vec{\tau}^h})\vert &\leq \Vert \vec{u}^h + \vec{u}'\Vert \sqrt{\nu}\Vert\curl{\vec{\tau}^h}\Vert\;, \\
    \vert(\tm \vec{u}',\vec{v}')\vert&\leq (\tm \vec{u}',\vec{u}') (\tm \vec{v}',\vec{v}')\;,\\
    \vert\sqrt{\frac{\nu}{2}}(\curl{\vec{\omega}'},{\vec{v}'})\vert &\leq \frac{1}{\sqrt{2}}\sqrt{\nu}\Vert \curl{\vec{\omega}'}\Vert_\mathcal{M} \Vert \tau_{\mathrm{M}}\Vert_{\infty} (\tm \vec{v}',\vec{v}')\;, \\
  \vert\left( \sqrt{\nu}\curl{\vec{\omega}^h},\vec{v}'\right)\vert&\leq \sqrt{\nu}\Vert \curl{\vec{\omega}^h}\Vert_\mathcal{M} \Vert \tau_{\mathrm{M}}\Vert_{\infty} (\tm \vec{v}',\vec{v}')\;, \\
  \vert(\vec{\omega}',\vec{\tau}')\vert &\leq \Vert \vec{\omega}'\Vert\Vert \vec{\tau}'\Vert\;, \\
  \vert \sqrt{\frac{\nu}{2}}({\vec{u}'},\curl{\vec{\tau}'})\vert  &\leq \frac{1}{\sqrt{2}}\sqrt{\nu}\Vert \curl{\vec{\tau}'}\Vert_\mathcal{M} \Vert \tau_{\mathrm{M}}\Vert_{\infty} (\tm \vec{u}',\vec{u}')\;.
    \end{align}
\end{subequations}
Note that $\Vert \tau_{\mathrm{M}}\Vert \leq \frac{h^e}{\Vert \vec{\beta}^h\Vert_{\infty}} < \infty$.
The only remaining term to estimate is $\vert\frac{1}{\sqrt{\nu}}(\vec{\omega}'\times\vec{\beta}^h,\vec{v}^h) \vert$, which relies on the following results.
Because of condition \eqref{eq:ass:beta-h}, for any vector \(\vec{u} \in \left[L^2(\Omega^e)\right]^d\), it holds that
\begin{align}
\nonumber \int_{{\Omega}^e} |\vec{u}\times\vec{\beta}^h|^2 dV & \leq \sum_{i=1,j=1}^{d,d}\int_{{\Omega}^e} | u_i \beta^h_j |^2 dV\leq \sum_{i=1,j=1}^{d,d}\Vert\vec{\beta}^h\Vert_{L^\infty(\Omega^e)}\int_{\Omega^e} |\vec{\beta}^h| |u_i|^2 dV =\\
    &d h^e \Vert \vec{\beta}^h\Vert_{L^\infty(\Omega^e)}\int_{\Omega^e} \frac{|\beta^h|}{h^e}|\vec{u}|^2 dV \leq d h^e \Vert \vec{\beta}^h\Vert_{L^\infty(\Omega^e)} (\tm \vec{u},\vec{u})_{L^2(\Omega^e)}^2 \;.
\end{align}
Hence,  we estimate
\begin{align}
    \nonumber \vert\frac{1}{\sqrt{\nu}}(\vec{\omega}'\times\vec{\beta}^h,\vec{v}^h) \vert&\leq \vert\frac{1}{\sqrt{\nu}}(\vec{\omega}'\times\vec{\beta}^h,\vec{v}^h+\vec{v}') \vert + \vert\frac{1}{\sqrt{\nu}}(\vec{\omega}'\times\vec{\beta}^h,\vec{v}') \vert\\
    &\leq \frac{2\Vert \vec{\beta}^h\Vert_{\infty} }{\sqrt{\nu}}\Vert \vec{\omega}'\Vert \Vert \vec{v}^h +\vec{v}' \Vert + \sqrt{\frac{d h^e\Vert\vec{\beta}^h\Vert_{\infty}}{\nu}}\Vert \vec{\omega}'\Vert(\tm \vec{v}',\vec{v}')\;.
\end{align}
As a result, for any $\vec{x},\vec{y}\in\mathbb{X}$, we have that
\[ \mathcal{A}_{\mathrm{red}}(x,y) \leq C_{\mathrm{cont}} \Vert \vec{x}\Vert_{\mathbb{X}_h} \Vert \vec{y}\Vert_{\mathbb{X}_h} \]
with $C_{\mathrm{cont}} $ independent on $\nu$. In fact, by Assumption \ref{ass:CFL-conditions}, it is estimated by
\begin{equation}
    C_{\mathrm{cont}} \leq \max\left( \sigma, \frac{\sqrt{\sigma}}{3},\frac{h^e}{\Vert\vec{\beta}^h\Vert_{\infty}}, \sqrt{\frac{d}{24}} \right) \leq \max\left( \sigma, \frac{\sqrt{\sigma}}{3},\frac{h_{\mathrm{max}}}{\Vert\vec{\beta}^h\Vert_{\infty}}, \sqrt{\frac{d}{24}} \right)\;,
\end{equation}
for some maximal elementsize $h_{\mathrm{max}}$ considered.

\subsection{Inf-sup stability \eqref{eq:stability-infsup}}
Given $x^h = ({\vec{\omega}}^h,\vec{u}^h,\vec{\omega}',\vec{u}')\in\mathbb{X}_h$, we pick $y^h = (\vec{\omega}^h, \vec{u}^h + \hat{c}\sqrt{\nu}\curl{\vec{\omega}^h}, \vec{\omega}',\vec{u}' + \hat{c}\sqrt{\nu}\curl{\vec{\omega}'})$. Then:
\begin{align}
\nonumber &\mathcal{A}_{\mathrm{red}}(x^h,y^h) \\
    \nonumber =& \sigma (\vec{u}^h+\vec{u}',\vec{u}^h) + \frac{1}{\sqrt{\nu}}(\vec{\omega}^h\times\left[\beta^h+\beta'\right],\vec{u}^h) + \frac{1}{\sqrt{\nu}}(\vec{\omega}'\times\left[\beta^h\right],\vec{u}^h) + \sqrt{\nu}(\curl{\vec{\omega}^h},\vec{u}^h) + \\
    \nonumber &\hat{c}\sqrt{\nu}\sigma (\vec{u}^h+\vec{u}',\curl{\vec{\omega}^h}) + \hat{c}(\vec{\omega}^h\times\left[\beta^h+\beta'\right],\curl{\vec{\omega}^h}) + \\
    \nonumber &\hat{c}(\vec{\omega}'\times\left[\beta^h\right],\curl{\vec{\omega}^h}) + \hat{c}\nu (\curl{\vec{\omega}^h},\curl{\vec{\omega}^h}) +\\
    \nonumber &(\vec{\omega}^h,\vec{\omega}^h) - \sqrt{\nu}(\vec{u}^h+\vec{u}',\curl{\vec{\omega}^h})+\\
    \nonumber &\sigma (\vec{u}^h+\vec{u}',\vec{u}') + (\tau_M^{-1} \vec{u}',\vec{u}') + \sqrt{\frac{\nu}{2}}(\curl{\vec{\omega}'},\vec{u}') + \frac{1}{\sqrt{\nu}}\left( \vec{\omega}^h\times[\beta^h+\beta'],\vec{u}' \right) +\\
    \nonumber &\sqrt{\nu}(\curl{\vec{\omega}^h},\vec{u}' )+
    \hat{c}\sqrt{\nu}\sigma (\vec{u}^h+\vec{u}',\curl{\vec{\omega}'}) + \hat{c}\sqrt{\nu}(\tau_M^{-1} \vec{u}',\curl{\vec{\omega}'}) + \\
    \nonumber &\hat{c}\nu\sqrt{\frac{1}{2}}(\curl{\vec{\omega}'},\curl{\vec{\omega}'}) + \hat{c}\left( \vec{\omega}^h\times[\beta^h+\beta'],\curl{\vec{\omega}'} \right) + \hat{c}\nu(\curl{\vec{\omega}^h},\curl{\vec{\omega}'} )+\\
    \nonumber &(\vec{\omega}',\vec{\omega}') - \sqrt{\frac{\nu}{2}}(\vec{u}',\curl{\vec{\omega}'})\\
    = \label{eq:to-estimate-infsup-sigma-formulation}&\sigma \Vert \vec{u}^h+\vec{u}'\Vert^2 + \Vert \vec{\omega}^h\Vert^2 + \hat{c}\nu\Vert\curl{\vec{\omega}^h}\Vert^2 + (\tau_M^{-1}\vec{u}',\vec{u}') + \frac{\hat{c}\nu}{\sqrt{2}}\Vert \curl{\vec{\omega}'}\Vert^2+\Vert \vec{\omega}'\Vert^2 +\\
    \nonumber &\frac{1}{\sqrt{\nu}}(\vec{\omega}^h\times\left[\beta^h+\beta'\right],\vec{u}^h) + \frac{1}{\sqrt{\nu}}(\vec{\omega}'\times\left[\beta^h\right],\vec{u}^h)+\hat{c}\sqrt{\nu}\sigma (\vec{u}^h+\vec{u}',\curl{\vec{\omega}^h}) + \\
    \nonumber &\hat{c}(\vec{\omega}^h\times\left[\beta^h+\beta'\right],\curl{\vec{\omega}^h})+\hat{c}(\vec{\omega}'\times\left[\beta^h\right],\curl{\vec{\omega}^h})+ \frac{1}{\sqrt{\nu}}\left( \vec{\omega}^h\times[\beta^h+\beta'],\vec{u}' \right) +\\
    \nonumber &\hat{c}\sqrt{\nu}\sigma (\vec{u}^h+\vec{u}',\curl{\vec{\omega}'}) + \hat{c}\sqrt{\nu}(\tau_M^{-1} \vec{u}',\curl{\vec{\omega}'}) + \hat{c}\left( \vec{\omega}^h\times[\beta^h+\beta'],\curl{\vec{\omega}'} \right) + \\
    \nonumber  &\hat{c}\nu(\curl{\vec{\omega}^h},\curl{\vec{\omega}'} )\;.
\end{align}
To estimate \eqref{eq:to-estimate-infsup-sigma-formulation}, we make use of the following estimates. The first (and most complex) estimate is given by:
\begin{align*}
    \frac{1}{\sqrt{\nu}}(\vec{\omega}'\times\beta^h,\vec{u}^h) &= \frac{1}{\sqrt{\nu}}(\vec{\omega}'\times\beta^h,\vec{u}^h + \vec{u}') - \frac{1}{\sqrt{\nu}}(\vec{\omega}'\times\beta^h,\vec{u}')\;.
\end{align*}
Where these terms are estimated by    
\begin{align*}
    |\frac{1}{\sqrt{\nu}}(\vec{\omega}'\times\beta^h,\vec{u}^h + \vec{u}')| &\geq -\frac{4 \Vert \beta^h\Vert^2_{\infty}}{\nu \sigma}\Vert \vec{\omega}' \Vert^2 - \frac{1}{8}\sigma \Vert \vec{u}^h + \vec{u}'\Vert^2\;,\\
    |\frac{1}{\sqrt{\nu}}(\vec{\omega}'\times\beta^h,\vec{u}')| &= \sum_{e\in\mathcal{M}} |\frac{1}{\sqrt{\nu}}(\vec{\omega}'\times\beta^h,\vec{u}')_{\Omega^e}|  \geq -d\mathrm{Re}_h \Vert\vec{\omega}'\Vert^2 - (\frac{1}{4}\tm\vec{u}',\vec{u}')\;.
\end{align*}
Here, $\mathrm{Re}_h := \max_{\Omega^e} \frac{\Vert \beta^h\Vert_{L^\infty(\Omega^e)} h^e}{\nu}$, is the maximal elementwise Reynolds number and $d$ is the number of dimensions (so $d=2$ or $d=3$).
For this last equation, we made use of the elementwise description of the integral to apply the following inequality over each element $\Omega^e$ and that $\frac{|\vec{\beta}^h|}{h}\leq \tm$:
\begin{align*}
\int_{{\Omega}^e} |\vec{u}\times\vec{\beta}|^2 dV & \leq \sum_{ij}\int_{{\Omega}^e} | u_i \beta_j |^2 dV\leq \sum_{ij}\Vert\vec{\beta}\Vert_{L^\infty(\Omega^e)}\int_{\Omega^e} |\vec{\beta}| |u_i|^2 dV =\\
    &d h^e \Vert \vec{\beta}\Vert_{L^\infty(\Omega^e)}\int_{\Omega^e} \frac{|\beta|}{h^e}|\vec{u}|^2 dV \leq d h^e \Vert \vec{\beta}\Vert_{L^\infty(\Omega^e)} (\tm \vec{u},\vec{u})_{L^2(\Omega^e)}^2 \;.
\end{align*}
So that we can estimate 
\begin{align*}
    |\left( (\vec{\vec{\omega}} \times \vec{\beta}) , \vec{u}\right)_{L^2(\Omega^e)}| &=|\left( (\vec{u} \times \vec{\beta}) , \vec{\vec{\omega}}\right)_{L^2(\Omega^e)}|  \leq  \Vert\vec{u} \times \vec{\beta}\Vert_{L^2({\Omega}^e)} \Vert\vec{\vec{\omega}}\Vert_{L^2({\Omega}^e)}\\
 &\leq \frac{d h^e \Vert \vec{\beta}\Vert_{L^\infty(\Omega^e)}}{\nu}\Vert \vec{\vec{\omega}}\Vert_{L^2(\Omega^e)}^2 + (\frac{1}{4}\tm\vec{u},\vec{u})_{L^2(\Omega^e)}\;.
\end{align*}
By summing over all mesh elements, we find the desired estimate.
The second non-trivial term is
\begin{align}
    \nonumber |\hat{c}\sqrt{\nu}(\tau_M^{-1} \vec{u}',\curl{\vec{\omega}'}) | &\geq -\frac{\hat{c}^2\nu \sigma}{2}\Vert\curl{\vec{\omega}'}\Vert^2  - \frac{1}{2\sigma}(\tm\vec{u}',\tau_M^{-1}\vec{u}')\\
    \label{eq:appendix:tmucurlu} &\geq -\frac{\hat{c}^2\nu \sigma}{2}\Vert\curl{\vec{\omega}'}\Vert^2  - \frac{1}{2}(\tm\vec{u}',\vec{u}')\;.
\end{align}
Where we used Assumption \ref{ass:CFL-conditions}. The remainder of the terms are estimated via similar estimates as in \cite{anaya_analysis_2019}, which result in the estimates:
\begin{subequations}
\begin{align*}
    \frac{1}{\sqrt{\nu}}(\vec{\omega}^h\times\left[\beta^h+\beta'\right],\vec{u}^h+\vec{u}') &\geq -\frac{4\Vert \beta^h + \beta'\Vert_{\infty}^2}{\nu\sigma}\Vert \vec{\omega}^h\Vert^2 - \frac{\sigma}{8}\Vert \vec{u}^h + \vec{u}'\Vert^2 \\
    \hat{c}\sqrt{\nu}\sigma (\vec{u}^h+\vec{u}',\curl{\vec{\omega}^h}) &\geq - \frac{\sigma}{8}\Vert \vec{u}^h+\vec{u}' \Vert^2 - 2\hat{c}^2\nu\sigma\Vert \curl{\vec{\omega}^h}\Vert^2 \\
    \hat{c}(\vec{\omega}^h\times\left[\beta^h+\beta'\right],\curl{\vec{\omega}^h})&\geq -\frac{\Vert \beta^h+\beta'\Vert^2}{\nu\sigma}\Vert \vec{\omega}^h\Vert^2 - \frac{1}{2}\hat{c}^2\nu\sigma\Vert \curl{\vec{\omega}^h}\Vert^2 \\
    \hat{c}(\vec{\omega}'\times\left[\beta^h\right],\curl{\vec{\omega}^h}) &\geq -\frac{\Vert \beta^h\Vert^2}{\nu\sigma}\Vert \vec{\omega}' \Vert^2 - \frac{1}{2}\hat{c}^2\nu\sigma\Vert \curl{\vec{\omega}^h}\Vert^2\\
    \hat{c}\sqrt{\nu}\sigma (\vec{u}^h+\vec{u}',\curl{\vec{\omega}'}) &\geq - \frac{\sigma}{8} \Vert \vec{u}^h+\vec{u}'\Vert^2- 2\hat{c}^2\nu\sigma\Vert\curl{\vec{\omega}'}\Vert^2\\
     \hat{c}\left( \vec{\omega}^h\times[\beta^h+\beta'],\curl{\vec{\omega}'} \right) &\geq - \frac{\Vert \beta^h+\beta'\Vert_{\infty}^2}{\sigma\nu}\Vert \vec{\omega}^h\Vert^2 - \frac{1}{2}\hat{c}^2\nu \sigma \Vert \curl{\vec{\omega}'}\Vert^2 \\ 
     \hat{c}\nu(\curl{\vec{\omega}^h},\curl{\vec{\omega}'} ) & \geq - \frac{\nu\hat{c}}{2} \Vert \curl{\vec{\omega}^h}\Vert^2 - \frac{\nu\hat{c}}{2} \Vert \curl{\vec{\omega}'}\Vert^2
\end{align*}
\end{subequations}
Then, \eqref{eq:to-estimate-infsup-sigma-formulation} can be estimated by using Assumption \ref{ass:CFL-conditions} and choosing $\hat{c} =  \frac{1}{7\sigma}$:
\begin{align}
    \nonumber \eqref{eq:to-estimate-infsup-sigma-formulation} \geq& \frac{\sigma}{2}\Vert \vec{u}^h + \vec{u}'\Vert^2 +\left[1 - 6\frac{\Vert \beta^h + \beta\Vert^2_{L^\infty}}{\sigma \nu}\right]\Vert\vec{\omega}^h\Vert^2+ \frac{1}{14 \sigma} \nu\Vert \curl{\vec{\omega}^h}\Vert^2\\
    \nonumber &+\frac{1}{4}(\tau_M^{-1}\vec{u}',\vec{u}')+ \left[1 - 5\frac{\Vert\beta^h+\beta'\Vert^2_{L^{\infty}}}{\nu\sigma} - d \frac{\Vert \beta^h\Vert_{L^{\infty}}h}{\nu}\right] \Vert\vec{\omega}'\Vert^2 \\
    &+ \frac{1}{14 \sigma} \nu \Vert\curl{\vec{\omega}'}\Vert^2 \geq C_1 \Vert x^h\Vert_{\mathbb{X}_h}^2\;,
\end{align}
where
\begin{equation}
    C_1 := \min \left(\frac{\sigma}{2},1-6\frac{\Vert \vec{\beta}^h+\vec{\beta}'\Vert_{L^\infty}^2 }{\sigma \nu},\frac{1}{14 \sigma} ,\frac{1}{4}, 1 - 5\frac{\Vert\beta^h+\beta'\Vert^2_{L^{\infty}}}{\nu\sigma} - d \frac{\Vert \beta^h\Vert_{L^{\infty}}h}{\nu}\right)\;.
\end{equation}
Then, with \eqref{eq:appendix:tmucurlu} and \eqref{eq:ass:tm-simga}, we estimate
\begin{align}
    \nonumber \Vert \vec{u}^h + \hat{c}\sqrt{\nu}\curl{\vec{\omega}^h} + \vec{u}' + \hat{c}\sqrt{\nu}\curl{\vec{\omega}'}\Vert^2 + (\tau_M^{-1}(\vec{u}' + \hat{c}\sqrt{\nu}\curl{\vec{\omega}'}),\vec{u}' + \hat{c}\sqrt{\nu}\curl{\vec{\omega}'})&\leq\\
    \nonumber \Vert \vec{u}^h + \vec{u}'\Vert^2 + \hat{c}^2\nu \Vert\curl{\vec{\omega}^h}\Vert^2 + \hat{c}^2\nu \Vert\curl{\vec{\omega}'}\Vert^2  +&\\
    \nonumber(\tm \vec{u}',\vec{u}') + 2 (\tm \vec{u}',\hat{c}\sqrt{\nu}\curl{\omega'}) + \hat{c}^2\nu(\tm\curl{\omega'},\curl{\omega'}) &\leq\\
    \Vert \vec{u}^h + \vec{u}'\Vert^2 + \hat{c}^2\nu \Vert\curl{\vec{\omega}^h}\Vert^2 + \hat{c}\left(\hat{c} + \hat{c}\sigma + \frac{1}{7}\right)\nu \Vert\curl{\vec{\omega}'}\Vert^2 + \frac{3}{2}(\tm \vec{u}',\vec{u'})&=\\
    \Vert \vec{u}^h + \vec{u}'\Vert^2 + \frac{1}{49\sigma^2}\nu \Vert\curl{\vec{\omega}^h}\Vert^2 + \frac{1}{7\sigma}\left(\frac{1}{7\sigma} + \frac{2}{7}\right)\nu \Vert\curl{\vec{\omega}'}\Vert^2 + \frac{3}{2}(\tm \vec{u}',\vec{u'})&\leq\\
    C_2^2 \left(\Vert \vec{u}^h + \vec{u}'\Vert^2 + \nu \Vert\curl{\vec{\omega}^h}\Vert^2 + \nu \Vert\curl{\vec{\omega}'}\Vert^2 + (\tm \vec{u}',\vec{u'}) \right)\;,&
\end{align}
where
\begin{equation}
    C_2^2 := \max\left(\frac{1}{49\sigma^2},\frac{1}{7\sigma}\left(\frac{1}{7\sigma} + \frac{2}{7}\right),\frac{3}{2}  \right)\;.
\end{equation}
Hence, we have that
\begin{equation}
    \Vert (\vec{\omega}^h,\vec{u}^h+\hat{c}\sqrt{\nu}\curl{\vec{\omega}^h},\vec{\omega}',\vec{u}'+\hat{c}\sqrt{\nu}\curl{\vec{\omega}'})\Vert_{\mathbb{X}_h}\leq C_2\Vert (\vec{\omega}^h,\vec{u}^h,\vec{\omega}',\vec{u}')\Vert_{\mathbb{X}_h}
\end{equation}
Showing the inf-sup condition \eqref{eq:stability-infsup} with constant $C_{\mathrm{inf-sup}}=C_1C_2^{-1}$, where the constant is independent of $\nu$ and mesh size $h^e$.\\

\subsection{Non-degeneracy \eqref{eq:stability-nondegenerate}}
We pick the same functions for the target and test functions. This results in:
\begin{align*}
    &\sigma\Vert \vec{u}^h+\vec{u}'\Vert^2 + \Vert \vec{\omega}^h\Vert^2 + \Vert \vec{\omega}'\Vert^2 + (\tau_M^{-1}\vec{u}',\vec{u}') + \frac{1}{\sqrt{\nu}}(\vec{\omega}^h\times\left[\beta^h+\beta'\right],\vec{u}^h+\vec{u}') + \frac{1}{\sqrt{\nu}}(\vec{\omega}'\times\beta^h,\vec{u}^h)\\
    \geq& \frac{7}{8}\sigma\Vert \vec{u}^h+\vec{u}'\Vert^2 + \Vert\vec{\omega}^h\Vert^2 + \left[ 1- 4\frac{\Vert\beta^h+\beta'\Vert_{\infty}^2}{\sigma \nu} - d\frac{\Vert \beta^h\Vert_{\infty}h}{\nu} \right]\Vert\vec{\omega}'\Vert^2 + \frac{3}{4}(\tau_M^{-1}\vec{u}',\vec{u}')
\end{align*}
Here, by estimating the resulting terms as in the inf-sup case and using the same conditions, we show condition \eqref{eq:stability-nondegenerate}.
\end{document}